\documentclass[a4paper,11pt,reqno]{article}

\usepackage[utf8]{inputenc}
\usepackage[T1]{fontenc}
\usepackage{lmodern}
\usepackage[english]{babel}
\usepackage{microtype}
\usepackage{csquotes}

\usepackage[backend=biber,hyperref=true,style=alphabetic,giveninits=true]{biblatex}
\usepackage{biblatex-mr}
\usepackage{tikz}

\definecolor{midnightblue}{rgb}{0.1, 0.1, 0.44}
\usepackage[colorlinks,allcolors=midnightblue]{hyperref}
\AtBeginDocument{}
\AtBeginDocument{}

\usepackage[a4paper,hmargin=2.5cm,bottom=3cm,top=3cm,footskip=3\baselineskip,
	marginparwidth=2.4cm,marginparsep=0.05cm]{geometry}

\usepackage{caption}
\usepackage{subcaption}

\usepackage{bm}

\usepackage{fsmath}

\usepackage[frozencache]{minted}

\newcommand{\loc}{\mathrm{loc}}
\DeclareMathOperator{\TV}{TV}

\newcommand{\dx}{\d x}

\newcommand{\dt}{\d t}

\newcommand{\II}{I\!I}
\newcommand{\III}{I\!I\!I}
\newcommand{\IV}{I\!V}

\begin{document}

%


\begin{center}
\vspace{5mm}

{\Large\bf%
Entropy solutions of non-local scalar conservation laws with congestion via deterministic particle method}

\vspace{5mm}

{\bf by}

\vspace{5mm}

{\bf Emanuela Radici \& Federico Stra}\\[1mm]
\noindent  {\it Institute of Mathematics, EPFL, CH-1015 Lausanne, Switzerland}\\[1mm]
email: \ {\tt emanuela.radici@epfl.ch}, \hspace*{1mm} {\tt federico.stra@epfl.ch}
\\[8mm]
\end{center}

\begin{abstract}
We develop deterministic particle schemes to solve non-local scalar conservation laws with congestion. We show that the discrete approximations converge to the unique entropy solution with an explicit rate of convergence under more general assumptions that the existing literature: the velocity fields are less regular (in particular the interaction force can have a discontinuity at the origin) with no prescribed attractive/repulsive regime and the mobility can have unbounded support. We complement our results with some numerical simulations, among which we show the applicability of the schemes to the multi-species setting.
\end{abstract}

\tableofcontents

\section{Introduction}

In this article we study the scalar conservation law
\begin{equation}\label{eq:conservation-law}
\partial_t\rho(t,x) + \div_x\bigl[\rho(t,x) v\bigl(\rho(t,x)\bigr)\bigl(V(t,x)-(W'*\rho)(t,x)\bigr)\bigr] = 0
\end{equation}
where the transport vector field depends non-locally on the solution itself via a convolution and $v:[0,\infty)\to[0,\infty)$ is a non increasing function to model congestion.
We are interested in the Cauchy problem with initial condition $\rho_0$ being a probability measure in $\Prob(\setR) \cap L^\infty(\setR) \cap BV(\setR)$ with compact support.
In the sequel we adopt the shorthand notations $m(\rho)=\rho v(\rho)$ for the mobility and $U(t,x)=V(t,x)-(W'*\rho)(t,x)$ for the uncongested velocity field. Moreover, for convenience we denote $W'=\partial_x W$.

The case $V=1, W=0$ has been studied for the first time in \cite{DiFrancescoRosini} with a similar approach. The case $V=0$ and with a regular attractive potential $W$ has been studied in \cite{DiFrancescoFagioliRadici}. The case $W=0$ and some symmetry assumptions on $V$ has been studied in \cite{DiFrancescoStivaletta}. In the recent article \cite{FagioliTse} the authors consider a general equation of the form \eqref{eq:conservation-law} and study the gradient flow structure with respect to a variation of the Wasserstein distance. A more detailed comparison between our result and theirs will be discussed after we present our main theorem and the assumptions on the mobility and the potentials.

Since weak solutions of \eqref{eq:conservation-law}, i.e.\ functions $\rho$ satisfying
\begin{equation}\label{eq:continuity-equation}
\int_0^\infty \int_\setR \bigl[
	\rho\partial_t\phi + m(\rho)U\partial_x\phi \bigr] \dx \dt = 0
\qquad \forall \phi\in C^\infty_c\bigl((0,\infty)\times\setR\bigr),
\end{equation}
are not unique in general, we are interested in studying a more restrictive class of solutions, \emph{entropy solutions} introduced by \cite{Kruzkov}, for which we can show uniqueness.

\begin{definition}[Entropy Solution]\label{def:entropy-solution}
A function $\rho:[0,\infty)\times\setR\to\setR$ is an \emph{entropy solution} to \eqref{eq:conservation-law} if
\begin{equation}\label{eq:entropy-condition}
\int_0^\infty \int_\setR \left\{\abs{\rho-c}\partial_t\phi
	+ \sign(\rho-c)\left[\bigl(m(\rho)-m(c)\bigr) U\partial_x\phi
	- m(c)\partial_x U \phi \right]\right\} \dx \dt \geq 0
\end{equation}
for every constant $c\geq0$ and non-negative test function $\phi\in C^\infty_c\bigl((0,\infty)\times\setR;[0,\infty)\bigr)$.
\end{definition}

The definition is adapted from \cite{volpert}, where the inequality is required only for the absolute value instead of all entropy-flux pairs $(\eta,q)$.
In addition, notice that in our problem the vector field $U$ depends non-locally on the solution itself, therefore the uniqueness for problem \eqref{eq:conservation-law} does not follow immediately from the results by Kru\v{z}kov.


The article revolves around the study of a Lagrangian discretization of the PDE \eqref{eq:conservation-law}, in which the solution $\rho$ is approximated by a piecewise constant density associated to particles evolving according to an explicit system of ODEs.
We construct two similar deterministic particle approximations (integrated and sampled), which differ only in the way the interaction term is computed. All the results contained in the article apply equally well to both schemes.


Such particle approximations enjoy good a priori estimates, compactness properties and a suitable approximate entropy inequality \eqref{eq:discrete-entropy}. By passing to the limit as the number of particles grows to infinity we can obtain the existence of an entropy solution for problem \eqref{eq:conservation-law}.


In addition, we demonstrate the practical feasibility of this particle discretization as a numerical scheme for computing good approximations of the density and the underlying flow.

In comparison to approaches that regularize the equation, such as Godunov or the Finite Elements Method, our scheme does not incur in a loss of resolution of the entropic shocks, which are instead always resolved to a sharp discontinuity.

We now introduce the above mentioned particle discretization, with the two alternative schemes (integrated and sampled), and then proceed to state our main results.

\subsection{Deterministic particle schemes}\label{sec:particle-method}

Given a diffuse probability $\rho\in\Prob(\setR)$, we discretize it by means of $N+1$ particles $(x_0,x_1,\dots,x_N)$ in increasing order with the property that $\rho\bigl((x_{i-1},x_i)\bigr)=1/N$ for every $i\in\{1,\dots,N\}$, from which we can reconstruct a piecewise constant density
\begin{align*}
\bar\rho &= \sum_{i=1}^N \rho_i \bm{1}_{(x_{i-1},x_i)}, &
\rho_i &= \frac{1}{N(x_i-x_{i-1})},
\end{align*}
which approximates the original probability $\rho$. For convenience, we define also the external densities $\rho_0=\rho_{N+1}=0$.

We then evolve these particles by making them flow with the transport vector field of the conservation law \eqref{eq:conservation-law}. Two terms require particular attention: the interaction $W'*\rho$ and the congestion $v(\rho)$.

The most obvious way to treat the interaction is to compute exactly the convolution with the measure $\bar\rho$ determined by the particles
\[
\begin{split}
(W'*\bar\rho)(x)
&= \sum_{i=1}^N \int_{x_{i-1}}^{x_i} W'(x-y)\rho_i \d y \\
&= \sum_{i=1}^N \rho_i[W(x-x_{i-1})-W(x-x_i)] \\
&= \sum_{i=0}^{N-1} (\rho_{i+1} - \rho_i)W(x-x_i) .
\end{split}
\]
We refer to this approach as \emph{integrated interaction}. An alternative is to sample the convolution at the points $(x_0,\dots,x_N)$
\[
(W'*\rho)(x)
\simeq \left(W'*\frac1{N+1}\sum_{i=0}^N\delta_{x_i}\right)(x)
= \frac1{N+1}\sum_{i=1}^N W'(x-x_i) .
\]
When $x$ is one of the points $x_i$ we have to interpret $W'(0)=0$. We will not make use of this form of interaction. A third option, suitable only when we want to compute the convolution at one of the points $x_i$, is to sample the convolution at the other points $x_j$
\[
(W'*\rho)(x_i)
\simeq \left(W'*\frac1N\sum_{j\neq i}\delta_{x_j}\right)(x_i)
= \frac1N\sum_{j\neq i} W'(x_i-x_j) .
\]
We refer to this approach as \emph{sampled interaction}. It coincides with the second approach multiplied by a factor $(N+1)/N$. It is equivalent to computing $W'*\dot\rho(x_i)$ where
\[
\dot\rho = \frac1N\sum_{i=0}^N\delta_{x_i}
\]
and we interpret $W'(0)=0$. Notice that $\dot\rho$ is not a probability because it has mass $(N+1)/N$, but when we use to sample the interaction in correspondence of a particles, then it is equivalent to using a probability because the particle itself does not count.

Regarding the congestion term $v(\rho)$, we need a way to decide which density to use in the computation. At each point $x_i$ we have the two densities on the left and right side
\begin{align*}
\lim_{x\to x_i^-}\bar\rho(x) &= \rho_i, &
\lim_{x\to x_i^+}\bar\rho(x) &= \rho_{i+1}.
\end{align*}
For the computation of the congestion, we decide to use the \emph{downwind} density, i.e.\ the density in the direction of motion, toward which the particle is traveling. This means that if the particle $x_i$ has to move to the right, meaning that the transport vector field is non negative, we use the right density $\rho_{i+1}$ in the computation of $v(\rho)$; if the particle has to move left, we use $\rho_i$ instead.


In conclusion, we can formulate the two following systems of ordinary differential equations describing the evolution of the particles $(x_0,\dots,x_N)$:
\begin{itemize}
\item the model with integrated interaction
\begin{align}
\tag{ODE$_I$}\label{eq:ode-integrated-interaction}
&\left\{\begin{aligned}
x_i'(t) &= v_i(t) \bar U_i(t), \\
\bar U_i(t) &= V\bigl(t,x_i(t)\bigr) - (W'*\bar\rho)\bigl(t,x_i(t)\bigr) \\
	&= V\bigl(t,x_i(t)\bigr) - \sum_{j=0}^N(\rho_{j+1}(t) - \rho_j(t)) W\bigl(t,x_i(t)-x_j(t)\bigr), \\
v_i(t) &= \begin{cases}
	v\bigl(\rho_i(t)\bigr), & \text{if } \bar U_i(t) < 0, \\
	v\bigl(\rho_{i+1}(t)\bigr), & \text{if } \bar U_i(t) \geq 0, \end{cases}
\end{aligned}\right. \\
\intertext{\item and the model with sampled interaction}
\tag{ODE$_S$}\label{eq:ode-sampled-interaction}
&\left\{\begin{aligned}
x_i'(t) &= v_i(t) \dot U_i(t), \\
\dot U_i(t) &= V\bigl(t,x_i(t)\bigr) - (W'*\dot\rho)\bigl(t,x_i(t)\bigr) \\
	&= V\bigl(t,x_i(t)\bigr) - \frac1N\sum_{j\neq i}W'\bigl(t,x_i(t)-x_j(t)\bigr), \\
v_i(t) &= \begin{cases}
	v\bigl(\rho_i(t)\bigr), & \text{if } \dot U_i(t) < 0, \\
	v\bigl(\rho_{i+1}(t)\bigr), & \text{if } \dot U_i(t) \geq 0. \end{cases}
\end{aligned}\right.
\end{align}
\end{itemize}

The only difference between the two models is how we compute the non-local interaction, so that the resulting vector field is either $\bar U=V-W'*\bar\rho$ or $\dot U=V-W'*\dot\rho$.

In order to work with the Cauchy problem, the last thing that remains to be specified is the initial condition.
Given $\rho_0\in\Prob(\setR)\cap L^\infty(\setR)\cap BV(\setR)$ with bound $\rho_0\leq R_0$ and compact support $\supp(\rho_0)\subseteq[-S_0,S_0]$, we take as initial condition for the systems of ODEs \eqref{eq:ode-integrated-interaction} and \eqref{eq:ode-sampled-interaction} a family of sorted particles $X_0^N=(x_{0,0}^N,x_{1,0}^N,\dots,x_{N,0}^N)$ such that the corresponding piecewise constant $\bar\rho_0^N$ satisfies for every $N\in\setN$
\begin{equation}\label{eq:initial-conditions}
\bar\rho_0^N \leq R_0, \qquad
\supp(\bar\rho_0^N) \subseteq [-S_0,S_0], \qquad
TV(\bar\rho_0^N) \leq B_0, \qquad
\bar\rho_0^N \to \rho_0 \text{ in $L^1(\setR)$}.
\end{equation}

Notice that in particular, given $\rho_0\in\Prob(\setR)\cap L^\infty(\setR)\cap BV(\setR)$ with bound $\rho_0\leq R_0$ and compact support $\supp(\rho_0)\subseteq[-S_0,S_0]$ and given $N\in\setN$, an initial condition $X_0^N=(x_{0,0}^N,x_{1,0}^N,\dots,x_{N,0}^N)$ satisfying \eqref{eq:initial-conditions} always exists, for instance, one can consider particles $x^N_{0,0}$ and $x^N_{N,0}$ such that $[x^N_{0,0},x^N_{N,0}]$ is the smallest interval containing $\supp(\rho_0)$, i.e.\ it is its convex hull, and then consider intermediate particles $x^N_{i,0}$ for $i=1,\dots,N-1$ with the property that $\rho([x^N_{i-1,0},x^N_{i,0}])=1/N$. In fact, clearly the support and the density bound are satisfied. Moreover, letting $\eta_\eps\in C^\infty_c\bigl([-\eps,\eps];[0,\infty)\bigr)$ be normalized mollifiers and $\rho_\eps=\rho_0*\eta_\eps$, then $\bar\rho_\eps^N\to\bar\rho_0^N$ in $L^1$ for $\eps\to0$ and
\[
\TV(\rho_\eps) = \norm{D\rho_0*\eta_\eps}_1 \leq \norm{D\rho_0}_1\norm{\eta_\eps}
\leq \norm{D\rho_0}_1 = \TV(\rho_0).
\]
Since $\rho_\eps$ is continuous, by the intermediate value theorem it is easy to check that $\TV(\bar\rho_\eps^N)\leq\TV(\rho_\eps)\leq\TV(\rho_0)$ for every $N$. Therefore by lower-semicontinuity we get for every $N$
\[
\TV(\bar\rho_0^N) \leq \liminf_{\eps\to0} \TV(\bar\rho_\eps^N) \leq \TV(\rho_0) \leq B_0.
\]
Finally, we need to show that $\bar\rho^N\to\rho_0$ in $L^1(\setR)$. Since $\bar\rho^N$ is bounded in $BV(\setR)$, we have that $\bar\rho^N\to\sigma$ in $L^1(\setR)$ for some $\sigma\in L^1(\setR)$. On the other hand,
\[
W_1(\bar\rho^N,\rho_0) \leq \sum_{i=1}^N \frac1N (x_i-x_{i-1}) = \frac{2S_0}N,
\]
hence $\bar\rho^N\weakto\rho_0$. By the uniqueness of the limit we deduce that $\sigma$ must be equal to $\rho_0$, and therefore $\bar\rho^N\to\rho_0$ in $L^1(\setR)$.

\subsection{Main results}\label{sec:main-results}

We present here the main result of the article. We are able to show the existence and uniqueness of the entropy solution for the equation \eqref{eq:conservation-law}, and prove that both particle schemes converge to this entropy solution.

\begin{theorem}\label{thm existence of entropy solution}
Let $v$, $V$ and $W$ satisfy the assumptions \autoref{as:assumptions} and let $\rho_0\in\Prob(\setR)\cap L^\infty(\setR)\cap BV(\setR)$ be compactly supported.
For every $N\in\setN$, let $\bar\rho^N$ be the piecewise constant density associated to the particles solving either \eqref{eq:ode-integrated-interaction} or \eqref{eq:ode-sampled-interaction}, with initial datum $\bar\rho_0^N$ satisfying \eqref{eq:initial-conditions} for some constants $R_0,S_0,B_0\in(0,\infty)$.

Then $\bar\rho^N$ converges in $L^1_\loc\bigl([0,\infty)\times\setR\bigr)$ to an entropic solution
\[
\rho\in 
L^\infty_\loc\bigl([0,\infty)\times\setR\bigr)
\cap BV_\loc\bigl([0,\infty)\times\setR\bigr)
\cap C\bigl([0,\infty);L^1(\setR)\bigr)
\]
of \eqref{eq:conservation-law} with initial datum $\rho_0$ (according to \autoref{def:entropy-solution}), with the property that for every $t\in[0,\infty)$
\begin{equation}\label{eq:solution-estimates}
\rho(t) \leq R(t), \qquad
\supp\bigl(\rho(t)\bigr) \subseteq [-S(t),S(t)], \qquad
\TV\bigl(\rho(t)\bigr) \leq B(t),
\end{equation}
for some increasing functions $R,S,B:[0,\infty)\to[0,\infty)$.
Moreover, the entropy solution of \eqref{eq:conservation-law} is unique.
Finally, the convergence $\bar\rho^N\to\rho$ in $L^1$ is quantitative with the rate
\begin{equation}\label{eq:convergence-rate}
\norm{\rho(t)-\bar\rho^N(t)}_{L^1(\setR)}
\leq \left( \norm{\rho(0)-\bar\rho^N(0)}_{L^1(\setR)} + \frac{4H(t)}N \right)
	\exp\bigl(t A(t)\bigr) ,
\end{equation}
for some increasing functions $A,H:[0,\infty)\to[0,\infty)$.
\end{theorem}


We now present the assumptions on the vector fields and the mobility for which we are able to obtain \autoref{thm existence of entropy solution}.

\begin{assumptions}\label{as:assumptions}
We make the following assumptions on the mobility $v$ and the fields $V$ and $W$.
There exist three non-negative non-decreasing functions $F\in L^\infty_\loc\bigl([0,\infty)\bigr)$, $G\in L^\infty_\loc\bigl([0,\infty)\bigr)$ and $\lambda:[0,\infty)\to[0,\infty)$ non-decreasing with $\int_y^\infty \frac{\d x}{\lambda(x)} = \infty$ for every $y>0$, such that:\footnote{With the notation $H\in L^p_\loc\bigl([0,\infty)\bigr)$ we mean that $H\in L^p\bigl([0,R)\bigr)$ for every $R>0$.}
\begin{itemize}
\item $v\in W^{1,\infty}_\loc\bigl([0,\infty);[0,\infty)\bigr)$ is a non-increasing function satisfying
\begin{equation}\label{eq:assumption-v}
\qquad
\norm{v'}_{L^\infty([0,r])} \leq G(r);
\end{equation}
\item either $W$ is repulsive at small scales, i.e.\ there exists $h>0$ such that
\begin{equation}\label{eq:assumption-W-repulsive}
\sign(x)W'(t,x)\leq0, \qquad \forall x\in[-h,h],
\end{equation}
or the congestion decays sufficiently fast so that
\begin{equation}\label{eq:assumption-v-decay}
\int_1^\infty \frac{\d r}{r^2v(r)}=\infty;
\end{equation}
\item $V(t,\plchldr)\in W^{2,1}_\loc(\setR)$, with the estimates
\begin{align}\label{eq:assumption-V}
\norm{V(t,\plchldr)}_{W^{1,\infty}([-x,x])}&\leq F(t)G(x) &
\norm{V''(t,\plchldr)}_{L^1([-x,x])}&\leq F(t)G(x) ;
\end{align}
\item $W(t,\plchldr)\in W^{1,\infty}_\loc(\setR)$,
$W'(t,\plchldr)\in BV_\loc(\setR)$ with distributional derivative $D_x W'(t,\plchldr) = W''(t,\plchldr)\leb^1 + w(t)\delta_0$, where $W''(t,\plchldr)\in W^{1,1}_\loc\bigl((-\infty,0]\bigr)\cap W^{1,1}_\loc\bigl([0,\infty)\bigr)$, with the estimates\footnote{By $W'''$ we denote the weak derivative of $W''$ in $(-\infty,0)\cup(0,\infty)$.}
\begin{equation}\label{eq:assumptions-W}
\begin{aligned}
\norm{W'(t,\plchldr)}_{L^\infty([-x,x])} &\leq F(t)G(x) , &
\abs{w(t)} &\leq F(t) , \\
\norm{W''(t,\plchldr)}_{L^\infty([-x,x])} &\leq F(t)G(x) , &
\norm{W'''(t,\plchldr)}_{L^1([-x,x])} &\leq F(t)G(x) ;
\end{aligned}
\end{equation}
\item for every $t>0$ and $x>0$ we have
\begin{equation}\label{eq:mild-growth}
V_+(t,x)\leq F(t)\lambda(x), \quad
V_-(t,-x)\leq F(t)\lambda(x), \quad
\norm{W'(t,\plchldr)}_{L^\infty([-2x,2x])}\leq F(t)\lambda(x) .
\end{equation}
\end{itemize}
\end{assumptions}

\begin{remark}
Notice that the assumption \eqref{eq:assumption-v} on $v$ is satisfied for instance if $v(r)\leq \log(r)/r$. A particular case is when $v(r)=0$ for $r$ larger than some constant $r_{\max}>0$. On the opposite side, another case of particular interest is $v=1$, as long as the interaction $W$ is repulsive at small scales.
Moreover, they hypothesis \eqref{eq:mild-growth} can be slightly generalized to
\[
\norm{W'(t,\plchldr)_+}_{L^\infty([-2x,0])}\leq F(t)\lambda(x), \qquad
\norm{W'(t,\plchldr)_-}_{L^\infty([0,2x])}\leq F(t)\lambda(x),
\]
allowing to work with interaction potentials with a strong long-distance attraction.
\end{remark}

This study has been motivated by \cite{DiFrancescoFagioliRadici}, which introduced a deterministic particle scheme to solve \eqref{eq:conservation-law} with $V=0$ and a regular attractive non-local interaction potential $W$. The main novelty of our scheme is the way $v$ intervenes in the ODEs for the particles. Namely, instead of weighing each individual interaction term with the congestion evaluated in the appropriate direction, we first compute the total velocity field $U=V-W*\rho$ and then multiply it by the congestion evaluated in the direction resulting from its sign.

This scheme appear to be more natural from the modelization point of view, in the sense that the agents first feel the velocity field that they want to follow and only afterwards they look at the congestion ahead.


In the recent \cite{FagioliTse} the authors investigate the gradient flow structure of \eqref{eq:conservation-law} by introducing the same scheme \eqref{eq:ode-sampled-interaction}. As a byproduct of their study, they are able to show that the many particle limit of $\bar\rho^N$ is an entropy solution of \eqref{eq:conservation-law} under more restrictive hypothesis than \autoref{as:assumptions}. In particular, they require that $v$ vanishes above a certain threshold, $V$ and $W'$ have linear growth and most importantly either $W\in W^{3,\infty}(\setR)$ or it has the specific form $W(x)=\pm\abs{x}$.

In the present article we introduce concurrently \eqref{eq:ode-sampled-interaction} and for the first time the more regular \eqref{eq:ode-integrated-interaction}, which at least from the numerical point of view allows to treat the multi-species problem.
Both schemes allow us to treat asymmetric and time-dependent potentials which are singular at the origin and to employ congestion terms $v$ which do not vanish above a threshold.
In addition we can recover an explicit convergence rate for $\bar\rho^N\to\rho$.

\subsection{Outline of the article}

In \autoref{sec:a-priori-estimates}, after some preliminary lemmas which compare the two velocity fields $\bar U$ and $\dot U$ and estimate their Lipschitz norm and discrete second derivative, we provide some a priori estimates for the piecewise constant densities $\bar\rho^N$ which allow us to deduce uniform bounds for the support (\autoref{lem:supp-bound}), the density (\autoref{lem:max-princ}) and the total variation (\autoref{prop:total-variation}). Thanks to a compactness argument developed in \cite{RossiSavare} which generalizes Aubin-Lion's lemma, we then deduce the convergence in $L^1$ of $\bar\rho^N$ (\autoref{thm:compactness}). Furthermore, we show that such a limit $\rho$ is continuous in time with respect to the spacial $L^1$ norm (\autoref{cor:continuity}).

In \autoref{sec:entropy-solution} we begin by showing an approximate entropy inequality for the piecewise constant densities $\bar\rho^N$ (\autoref{prop:discrete-entropy-condition}). We pass this inequality to the limit with \autoref{prop:entropic-limit}, deducing that the density $\rho$ is an entropy solution. Finally, we provide a generalization of the stability theorem of \cite{KarlsenRisebro} allowing for time-dependent potentials and approximate entropy condition (\autoref{thm:KR}), and use it to deduce both the uniqueness of the entropy solution and the convergence rate of the particle scheme.

In \autoref{sec:numerics} we complement the theoretical results with an implementation of the numerical schemes and demonstrate its applicability to some interesting examples.

\subsection{Future perspectives}

The most natural extension of this work would be to relax the \autoref{as:assumptions} in order to allow interactions which are more singular in the origin, for instance Riesz potentials $W(x)=\pm\frac{1}{\abs{x}^\alpha}$ or Lennard-Jones potentials $W(x)=\pm\left(\frac{1}{\abs{x}^\alpha}-\frac{1}{\abs{x}^\beta}\right)$ for $\alpha,\beta>0$. Among the technical difficulties there is establishing a bound for the growth of the total variation.

Another possible generalization consists in adding a diffusion term in the conservation law \eqref{eq:conservation-law}.
In the setting with a non-linear mobility this has already been investigated in \cite{FagioliRadiciDiffusion}, for the choice $V=0$ and $W$ a regular attractive non-local interaction, working in a fixed bounded interval $[0,L]$, with a strictly positive initial datum.
A similar study has been conducted in \cite{FagioliRadiciOpinion} in the setting of \emph{opinion dynamics}, which is a variation in which the congestion $v(\rho)$ is replaced by a function of space $v(x)$. Finally, in \cite{DaneriRadici} the authors consider the diffusive conservation law without congestion, $V=0$ and a symmetric interaction $W$ and show the convergence of the particle scheme in the unbounded domain $\setR$ by exploiting the gradient flow structure of the equation.
In this case it is not natural to derive a total variation estimate for the solution because boundedness in $BV$ is not expected in general.
A possible future objective would be to study the equation \eqref{eq:conservation-law} with a diffusion term $\Delta A(\rho)$ with $A$ increasing under our more general \autoref{as:assumptions} in the entire domain $\setR$. A delicate point is to prove the discrete entropy inequality when the supports are not equi-bounded.

Another direction for future development would be to consider a system of conservation laws
\begin{equation}\label{eq:multi-species}
\partial_t\rho_i + \div_x\Biggl[\rho_i v_i(\rho)
	\Biggl(V_i - \sum_{j=1}^S W'_{i,j}*\rho_j \Biggr)\Biggr] = 0,
\qquad \forall i \in \{1,\dots,S\},
\end{equation}
where $\rho=(\rho_1,\dots,\rho_S)$ are $S$ different species subject to external fields $(V_1,\dots,V_S)$, mutual interactions $(W'_{i,j})_{i,j=1}^S$ and congestions $(v_1,\dots,v_S)$. While from the computational point of view this extension does not pose particular challenges, we are not able to prove however any convergence of the numerical scheme for this case, leaving it for future research. We provide some examples in the last section.

Leaving the world of conservation laws, it would be interesting to add a source term to the right hand side of the equation \eqref{eq:conservation-law} and generalize this particle approach to the case where mass is not conserved.

Finally, having in mind applications to traffic flow, instead of working on the real line one could pose the conservation law as a PDE on the edges of a planar graph, with boundary rules on how the traffic flows at the vertices.

\section{Theoretical analysis}

In what follows we discuss the many particle limit for both the models. Whenever the difference between $\bar U$ and $\dot U$ does not affect the proof of a statement, we directly use the unified notation where $U$ indicates both the sampled and the integrated fields at the same time.

\subsection{A priori estimates and compactness}\label{sec:a-priori-estimates}


\begin{lemma}\label{lem:discrete-derivatives}
Let $X=(x_0,\dots,x_N)$ be sorted particles and let $U_i$ be either $\bar U_i$ or $\dot U_i$ as defined in \eqref{eq:ode-integrated-interaction} or \eqref{eq:ode-sampled-interaction} respectively. Then
\begin{align}
\abs{U_i-U_{i-1}}
&\leq [C_1+C_2\rho_i] (x_i-x_{i-1})
	\label{eq:lipschitz-bound} \\
\abs{\dot U_{i+1}-2\dot U_i+\dot U_{i-1}}
&\leq \label{eq:second-derivative-bound-sampled}
	C_1 \frac{\abs{\rho_{i+1}-\rho_i}}{N\rho_i\rho_{i+1}}
	+ \frac1N \left(\frac1{\rho_i}+\frac1{\rho_{i+1}}\right) \left(
		\frac{C_3}N + \frac12\int_{x_{i-1}}^{x_{i+1}} \abs{V''(x)}\d x
		\right)
	\spliteq\notag
	+ \frac1{N^2\rho_i} \sum_{j\neq i, i\pm 1}
		\int_{x_{i-1}}^{x_{i+1}} \abs{W'''(y-x_j)} \d y \\
\abs{\bar U_{i+1}-2\bar U_i+\bar U_{i-1}}
&\leq \label{eq:second-derivative-bound-integrated}
	C_1 \frac{\abs{\rho_{i+1}-\rho_i}}{N\rho_i\rho_{i+1}}
	+ \frac1N \left(\frac1{\rho_i}+\frac1{\rho_{i+1}}\right) \left(
		\frac{C_3}N + \frac12\int_{x_{i-1}}^{x_{i+1}} \abs{V''(x)}\d x
		\right)
	\spliteq\notag
	+ \frac1{N\rho_i}
		\int_{x_{i-1}}^{x_{i+1}} \abs{W'''}*\bar\rho(y) \d y .
\end{align}
where
\begin{gather*}
C_1 = \norm{V'}_{L^\infty([x_0,x_N])} + \norm{W''}_{L^\infty([x_0-x_N,x_N-x_0])},
\qquad
C_2 = 2 \norm{W'}_{L^\infty([x_0-x_N,x_N-x_0])}, \\
C_3 = 3 \norm{W''}_{L^\infty([x_0-x_N,x_N-x_0])} .
\end{gather*}
\end{lemma}

\begin{proof}
Let us start with \eqref{eq:lipschitz-bound} first. The field $V$ itself is locally Lipschitz, so $\abs{V_i-V_{i-1}}\leq\norm{V'}_{L^\infty([x_0,x_N])} (x_i-x_{i-1})$. Regarding the interaction term, we treat the sampled and the interaction version separately. For the sampled one we have
\[
\begin{split}
&\hspace{-.25cm}\abs{W'*\dot\rho(x_i)-W'*\dot\rho(x_{i-1})} \\
&\leq \frac1N \sum_{j\neq i-1,i} \abs{W'(x_i-x_j)-W'(x_{i-1}-x_j)}
	+ \frac1N\abs{W'(x_i-x_{i-1})-W'(x_{i-1}-x_i)} \\
&\leq \frac1N \sum_{j\neq i-1,i} \norm{W''}_{L^\infty([x_0-x_N,x_N-x_0])}(x_i-x_{i-1})
	+2 \norm{W'}_{L^\infty([x_0-x_N,x_N-x_0])} \rho_i (x_i-x_{i-1}) \\
&\leq \bigl[\norm{W''}_{L^\infty([x_0-x_N,x_N-x_0])}
	+2 \norm{W'}_{L^\infty([x_0-x_N,x_N-x_0])} \rho_i\bigr] (x_i-x_{i-1}).
\end{split}
\]
We were able to use the Lipschitz estimate for $W'$ because for every $j\neq i-1,i$ the two terms $x_i-x_j$ and $x_{i-1}-x_j$ have the same sign.
For the integrated interaction on the other hand, we have the estimate
\[
\begin{split}
&\hspace{-.25cm}\abs{W'*\bar\rho(x_i)-W'*\bar\rho(x_{i-1})} \\
&\leq \sum_{j\neq i} \int_{x_{j-1}}^{x_j} \abs{W'(x_i-x)-W'(x_{i-1}-x)} \rho_j \d x
	+ \int_{x_{i-1}}^{x_i} \abs{W'(x_i-x)-W'(x_{i-1}-x)} \rho_i \d x \\
&\leq \norm{W''}_{L^\infty([x_0-x_N,x_N-x_0])}(x_i-x_{i-1})
	+ \frac2N \norm{W'}_{L^\infty([x_0-x_N,x_N-x_0])} \\
&\leq \bigl[\norm{W''}_{L^\infty([x_0-x_N,x_N-x_0])}
	+2 \norm{W'}_{L^\infty([x_0-x_N,x_N-x_0])} \rho_i\bigr] (x_i-x_{i-1}).
\end{split}
\]

Let us now move on to discrete second derivative estimates.
We deal first with the external potential $V$.
\[
\begin{split}
V(x_{i+1})-2V(x_i)+V(x_{i-1})
&= V'(x_i)(x_{i+1}-x_i) + \frac12\int_{x_i}^{x_{i+1}} V''(x)(x_{i+1}-x)\d x
	\spliteq
	+ V'(x_i)(x_{i-1}-x_i) + \frac12\int_{x_i}^{x_{i-1}} V''(x)(x_{i-1}-x)\d x \\
&= V'(x_i)(x_{i+1}-2x_i+x_{i-1})
	\spliteq
	+ \frac12\int_{x_i}^{x_{i+1}} V''(x)(x_{i+1}-x)\d x
	+ \frac12\int_{x_i}^{x_{i-1}} V''(x)(x_{i-1}-x)\d x ,
\end{split}
\]
hence
\begin{equation}\label{eq:second-derivative-V}
\begin{split}
&\hspace{-.5cm} \abs{V(x_{i+1})-2V(x_i)+V(x_{i-1})} \\
&\leq \frac{\abs{\rho_{i+1}-\rho_i}}{N\rho_i\rho_{i+1}} \norm{V'}_{L^\infty([x_0,x_N])}
	+ \frac1{2N\rho_i}\int_{x_{i-1}}^{x_i} \abs{V''(x)}\d x
	+ \frac1{2N\rho_{i+1}}\int_{x_i}^{x_{i+1}} \abs{V''(x)}\d x \\
&\leq \frac{\abs{\rho_{i+1}-\rho_i}}{N\rho_i\rho_{i+1}} \norm{V'}_{L^\infty([x_0,x_N])}
	+ \frac1{2N} \left(\frac1{\rho_i}+\frac1{\rho_{i+1}}\right)
		\int_{x_{i-1}}^{x_{i+1}} \abs{V''(x)}\d x
\end{split}
\end{equation}

Next, we have to deal with the interaction term. Here we have to distinguish between the sampled and integrated scheme since the terms $W'*\dot\rho$ and $W'*\bar\rho$ are treated slightly differently.

In the setting of the sampled particle scheme, recalling that $W$ is in $W^{3,1}_\loc$ on both $(0,\infty)$ and $(-\infty,0)$, we can compute
\[
\begin{split}
&\abs{W'*\dot\rho (x_{i+1}) - 2W'*\dot\rho(x_i) + W'*\dot\rho(x_{i-1})} \\
&\leq \frac1N \sum_{j\neq i, i\pm 1} \abs*{
		\int_{x_i}^{x_{i+1}} W''(x-x_j)\d x
		- \int_{x_{i-1}}^{x_i} W''(x-x_j)\d x }
		\spliteq
	+ \biggl\lvert \frac1N [W'(x_{i+1}-x_{i-1})+W'(x_{i+1}-x_i)]
	-\frac2N [W'(x_i-x_{i-1})+W'(x_i-x_{i+1})]
		\spliteq\qquad
	+\frac1N [W'(x_{i-1}-x_i)+W'(x_{i-1}-x_{i+1})] \biggr\rvert \\
&\leq \frac1N \sum_{j\neq i, i\pm 1}
		\int_0^1 \abs*{ W''\bigl(x_i+t(x_{i+1}-x_i)-x_j\bigr)(x_{i+1}-x_i)
			-W''\bigl(x_{i-1}+t(x_i-x_{i-1})-x_j\bigr)(x_i-x_{i-1}) } \d t
		\spliteq
	+\frac1N \abs{W'(x_{i+1}-x_{i-1})-W'(x_i-x_{i-1})}
	+\frac1N \abs{W'(x_i-x_{i+1})-W'(x_{i-1}-x_{i+1})}
		\spliteq
	+\frac1N \abs{W'(x_{i+1}-x_i)-W'(x_i-x_{i-1})}
	+\frac1N \abs{W'(x_{i-1}-x_i)-W'(x_i-x_{i+1})} \\
&\leq \frac1{N^2} \abs*{\frac{1}{\rho_i}-\frac{1}{\rho_{i+1}}}\sum_{j\neq i, i\pm 1}
		\int_0^1 \abs*{ W''\bigl(x_i+t(x_{i+1}-x_i)-x_j\bigr) }\d t
		\spliteq
	+\frac1{N^2\rho_i} \sum_{j\neq i, i\pm 1}
		\int_0^1 \int_{x_{i-1}+t(x_i-x_{i-1})}^{x_i+t(x_{i+1}-x_i)}
			\abs{W'''(y-x_j)} \d y \d t
		\spliteq
	+ \frac1N\norm{W''}_{L^\infty([x_0-x_N,x_N-x_0])}
		[(x_{i+1}-x_i) + (x_i-x_{i-1}) + 2(x_{i+1}-x_{i-1})] \\
&\leq \frac{\abs{\rho_{i+1}-\rho_i}}{N\rho_i} \sum_{j\neq i, i\pm 1}
		\int_{x_i}^{x_{i+1}} \abs{W''(y-x_j)} \d y
	+ \frac1{N^2\rho_i} \sum_{j\neq i, i\pm 1}
		\int_{x_{i-1}}^{x_{i+1}} \abs{W'''(y-x_j)} \d y
		\spliteq
	+ \frac3N\norm{W''}_{L^\infty([x_0-x_N,x_N-x_0])} (x_{i+1}-x_{i-1}),
\end{split}
\]
where in the last step, by Fubini, we used the estimate
\[
\begin{split}
\int_0^1 \int_{x_{i-1}+t(x_i-x_{i-1})}^{x_i+t(x_{i+1}-x_i)} \abs{W'''(y-x_j)} \d y \d t
&\leq \int_{x_{i-1}}^{x_{i+1}} \int_0^1 \abs{W'''(y-x_j)}
	\min\left\{\frac{y-x_{i-1}}{x_i-x_{i-1}},\frac{x_{i+1}-y}{x_{i+1}-x_i}\right\}
		\d t \d y \\
&\leq \int_{x_{i-1}}^{x_{i+1}} \abs{W'''(y-x_j)} \d y.
\end{split}
\]
Then \eqref{eq:second-derivative-bound-sampled} follows from \eqref{eq:second-derivative-V}, the previous computation and
\begin{align*}
\sum_{j\neq i, i\pm 1} \int_{x_i}^{x_{i+1}} \abs{W''(y-x_j)} \d y
&\leq \frac1{\rho_{i+1}} \norm{W''}_{L^\infty([x_0-x_N,x_N-x_0])}, &
x_{i+1}-x_{i-1}
&= \frac1N \left(\frac1{\rho_i}+\frac1{\rho_{i+1}}\right).
\end{align*}


Let us consider now the integrated scheme.
The interaction can be split into two sums of long range and short range interactions respectively
\[
\begin{split}
W'*\bar\rho(x_{i+1})&-2W'*\bar\rho(x_i)+W'*\bar\rho(x_{i-1}) \\
&= \int_{-\infty}^\infty [W'(x_{i+1}-x)-2W'(x_i-x)+W'(x_{i-1}-x)] \bar\rho(x) \d x \\
&= \sum_{j\neq i,i+1} \rho_j
	\int_{x_{j-1}}^{x_j} [W'(x_{i+1}-x)-2W'(x_i-x)+W'(x_{i-1}-x)]\d x
	\spliteq
	+\sum_{j=i,i+1} \rho_j
	\int_{x_{j-1}}^{x_j} [W'(x_{i+1}-x)-2W'(x_i-x)+W'(x_{i-1}-x)]\d x \\
&= \sum_{j\neq i,i+1} \rho_j \int_{x_{j-1}}^{x_j} \left(
		\int_{x_i}^{x_{i+1}} W''(y-x) \d y - \int_{x_{i-1}}^{x_i} W''(y-x) \d y
	\right) \d x
	\spliteq
	+\sum_{j=i,i+1} \rho_j
	\int_{x_{j-1}}^{x_j} [W'(x_{i+1}-x)-2W'(x_i-x)+W'(x_{i-1}-x)]\d x .
\end{split}
\]
With a computation similar to the one used for the sampled interaction, the first sum can be bounded in absolute value by
\[
\begin{split}
&\sum_{j\neq i,i+1} \rho_j \int_{x_{j-1}}^{x_j} \abs*{
		\int_{x_i}^{x_{i+1}} W''(y-x) \d y - \int_{x_{i-1}}^{x_i} W''(y-x) \d y
	} \d x \\
&\leq \frac{\abs{\rho_{i+1}-\rho_i}}{\rho_i} \sum_{j\neq i, i\pm 1}
		\int_{x_i}^{x_{i+1}} \int_{x_{j-1}}^{x_j} \abs{W''(y-x)} \rho_j \d x \d y
	+ \frac1{N\rho_i} \sum_{j\neq i, i\pm 1}
		\int_{x_{i-1}}^{x_{i+1}} \int_{x_{j-1}}^{x_j} \abs{W'''(y-x)} \rho_j \d x \d y \\
&\leq \frac{\abs{\rho_{i+1}-\rho_i}}{\rho_i}
		\norm{W''}_{L^\infty([x_0-x_N,x_N-x_0])} (x_{i+1}-x_i)
	+ \frac1{N\rho_i} \int_{x_{i-1}}^{x_{i+1}} \int_{x_0}^{x_N}
		\abs{W'''(y-x)}\bar\rho(x) \d x \d y \\
&\leq \frac{\abs{\rho_{i+1}-\rho_i}}{N\rho_i\rho_{i+1}}
		\norm{W''}_{L^\infty([x_0-x_N,x_N-x_0])}
	+ \frac1{N\rho_i}
		\int_{x_{i-1}}^{x_{i+1}} \abs{W'''}*\bar\rho(y) \d y .
\end{split}
\]
On the other hand, the second sum can be rewritten as
\[
\begin{split}
\Biggl\lvert \rho_i & \int_{x_{i-1}}^{x_i} [W'(x_{i+1}-x)-2W'(x_i-x)+W'(x_{i-1}-x)] \dx
	\spliteq
	+ \rho_{i+1} \int_{x_i}^{x_{i+1}} [W'(x_{i+1}-x)-2W'(x_i-x)+W'(x_{i-1}-x)] \dx \Biggr\rvert \\
&\leq \rho_i \int_{x_{i-1}}^{x_i} \abs{W'(x_{i+1}-x)-W'(x_i-x)} \dx
	+ \rho_{i+1} \int_{x_i}^{x_i+1} \abs{W'(x_i-x)-W'(x_{i-1}-x)} \dx
	\spliteq
	+ \abs*{ \rho_i \int_{x_{i-1}}^{x_i} [W'(x_i-x)-W'(x_{i-1}-x)] \dx
	- \rho_{i+1} \int_{x_i}^{x_{i+1}} [W'(x_{i+1}-x)-W'(x_i-x)] \dx } \\
&\leq \frac1N \norm{W''}_{L^\infty([x_0-x_N,x_N-x_0])} (x_{i+1}-x_{i-1})
	\spliteq
	+ \frac1N \Biggl\lvert \int_0^1 \Bigl\{
		\bigl[W'\bigl(x_i - x_{i-1}-t(x_i-x_{i-1})\bigr)
		-W'\bigl(x_{i-1} - x_{i-1}-t(x_i-x_{i-1})\bigr)\bigr]
	\spliteq\qquad\qquad
		-\bigl[W'\bigl(x_{i+1} - x_i-t(x_{i+1}-x_i)\bigr)
		-W'\bigl(x_i - x_i-t(x_{i+1}-x_i)\bigr)\bigr] \Bigr\} \dt \Biggr\rvert \\
&\leq \frac1N \norm{W''}_{L^\infty([x_0-x_N,x_N-x_0])} (x_{i+1}-x_{i-1})
	\spliteq
	+ \frac1N \int_0^1 \Bigl\{
	\abs{ W'\bigl((1-t)(x_i - x_{i-1})\bigr)
		- W'\bigl((1-t)(x_{i+1} - x_i)\bigr) }
	\spliteq\qquad\qquad
	+\abs{ W'\bigl(-t(x_i-x_{i-1})\bigr)
			-W'\bigl(-t(x_{i+1}-x_i)\bigr) } \Bigr\} \dt \\
&\leq \frac3N \norm{W''}_{L^\infty([x_0-x_N,x_N-x_0])} (x_{i+1}-x_{i-1})
\leq \frac3{N^2} \norm{W''}_{L^\infty([x_0-x_N,x_N-x_0])}
	\left(\frac1{\rho_i}+\frac1{\rho_{i+1}}\right) .
\end{split}
\]
In summary, we have that
\[
\begin{split}
\abs{W'*\bar\rho(x_{i+1})&-2W'*\bar\rho(x_i)+W'*\bar\rho(x_{i-1})} \\
&\leq \frac{\abs{\rho_{i+1}-\rho_i}}{N\rho_i\rho_{i+1}}
		\norm{W''}_{L^\infty([x_0-x_N,x_N-x_0])}
	+ \frac1{N\rho_i}
		\int_{x_{i-1}}^{x_{i+1}} \abs{W'''}*\bar\rho(y) \d y
	\spliteq
	+\frac3{N^2} \norm{W''}_{L^\infty([x_0-x_N,x_N-x_0])}
	\left(\frac1{\rho_i}+\frac1{\rho_{i+1}}\right) ,
\end{split}
\]
which, together with \eqref{eq:second-derivative-V}, concludes the proof of \eqref{eq:second-derivative-bound-integrated}.
\end{proof}

\begin{lemma}\label{lem:comparison}
Let $X=(x_0,\dots,x_N)$ be sorted particles and let $\bar U_i$ and $\dot U_i$ be defined as in \eqref{eq:ode-integrated-interaction} and \eqref{eq:ode-sampled-interaction}. Then
\begin{align}
\abs{\bar U_i-\dot U_i}
&\leq \frac1N \norm{W''}_{L^\infty([x_0-x_N,x_N-x_0])} (x_N-x_0),
	\label{eq:comparison-0}\\
\abs{(\bar U_i-\dot U_i) - (\bar U_{i-1}-\dot U_{i-1})}
&\leq \frac2N \bigl(\norm{W'''}_{L^1([x_0-x_N,x_N-x_0])}
	+\norm{W''}_{L^\infty([x_0-x_N,x_N-x_0])}\bigr) (x_i-x_{i-1}) .
	\label{eq:comparison-1}
\end{align}
\end{lemma}

\begin{proof}
We directly compute
\[
\begin{split}
&\hspace{-.25cm} \abs{\bar U_i-\dot U_i} \\
&\leq
	\sum_{j=1}^i \int_{x_{j-1}}^{x_j} \abs{W'(x_i-x)-W'(x_i-x_{j-1})} \rho_j\d x
	+\sum_{j=i+1}^N \int_{x_{j-1}}^{x_j} \abs{W'(x_i-x)-W'(x_i-x_j)} \rho_j\d x \\
&\leq \sum_{j=1}^i \norm{W''}_{L^\infty([0,x_N-x_0])} \rho_j (x_j-x_{j-1})^2
	+\sum_{j=i+1}^N \norm{W''}_{L^\infty([x_0-x_N,0])} \rho_j (x_j-x_{j-1})^2 \\
&\leq \norm{W''}_{L^\infty([x_0-x_N,x_N-x_0])} \frac1N \sum_{j=1}^N (x_j-x_{j-1}) ,
\end{split}
\]
which implies \eqref{eq:comparison-0}, and
\[
\begin{split}
&\abs{(\bar U_i-\dot U_i) - (\bar U_{i-1}-\dot U_{i-1})} \\
&\leq \sum_{j=1}^{i-1} \int_{x_{j-1}}^{x_j} \abs*{
	\bigl[W'(x_i-x)-W'(x_{i-1}-x)\bigr]
	-\bigl[W'(x_i-x_{j-1})-W'(x_{i-1}-x_{j-1})\bigr]
	} \rho_j \d x
		\spliteq
	+\sum_{j=i+1}^{N} \int_{x_{j-1}}^{x_j} \abs*{
	\bigl[W'(x_i-x)-W'(x_{i-1}-x)\bigr]
	-\bigl[W'(x_i-x_j)-W'(x_{i-1}-x_j)\bigr]
	} \rho_j \d x
		\spliteq
	+ \int_{x_{i-1}}^{x_i} \abs*{
		\bigl[W'(x_i-x)-W'(x_{i-1}-x)\bigr]
		-\bigl[W'(x_i-x_{i-1})-W'(x_{i-1}-x_i)\bigr]
	} \rho_i \d x \\
&\leq \sum_{j=1}^{i-1} \rho_j \int_{x_{j-1}}^{x_j} \int_{x_{i-1}}^{x_i} \int_{x_{j-1}}^x
	\abs{W'''(y-z)} \d z \d y \d x
	+\sum_{j=i+1}^N \rho_j \int_{x_{j-1}}^{x_j} \int_{x_{i-1}}^{x_i} \int_x^{x_j}
	\abs{W'''(y-z)} \d z \d y \d x
		\spliteq
	+ \int_{x_{i-1}}^{x_i} \Bigl[\abs{W'(x_i-x)-W'(x_i-x_{i-1})}
	+ \abs{W'(x_{i-1}-x)-W'(x_{i-1}-x_i)}\Bigr] \rho_i \d x \\
&\leq \sum_{j=1}^{i-1} \rho_j \int_{x_{i-1}}^{x_i} \int_{x_{j-1}}^{x_j}
	\abs{W'''(y-z)} (x_j-z) \d z \d y
	+\sum_{j=i+1}^N \rho_j \int_{x_{i-1}}^{x_i} \int_{x_{j-1}}^{x_j}
	\abs{W'''(y-z)} (z-x_{j-1}) \d z \d y
		\spliteq
	+\frac1N \Bigl(\norm{W''}_{L^\infty([0,x_N-x_0])}
		+\norm{W''}_{L^\infty([x_0-x_N,0])}\Bigr) (x_i-x_{i-1}) \\
&\leq \frac1N \int_{x_{i-1}}^{x_i} \int_{x_0}^{x_N} \abs{W'''(y-z)} \d z \d y
	+ \frac2N \norm{W''}_{L^\infty([x_0-x_N,x_N-x_0])} (x_i-x_{i-1}) ,
\end{split}
\]
which implies \eqref{eq:comparison-1}.
\end{proof}

Thanks to the above estimates on the velocity fields, we can now proceed to show uniform bounds independent of $N$ for both the support and the $L^\infty$ norm of piecewise constant densities $\bar\rho$ originating from particles solving \eqref{eq:ode-integrated-interaction} or \eqref{eq:ode-sampled-interaction}. The main tool is the Gronwall inequality and its generalization due to \cite{Bihari}, and \cite{ButlerRogers}. Since the latter is less know, for the reader's convenience we recall here the statement, taken from \cite{ButlerRogers}.
\begin{theorem}[Bihari, Butler-Rogers]
Let $x(t),a(t),b(t)$ be positive functions of $t$, bounded in $c\leq t\leq d$, let $k(t,s)$ be nonnegative, bounded on the triangular region $c\leq s\leq t\leq d$; assume further that $x(t)$ is measurable and $k(t,s)$ is a measurable function of $s$ for each $t$. Let $f(u),g(u)$ be positive functions for $u\geq0$, with $f$ strictly increasing and $g$ nondecreasing. Then defining
\begin{gather*}
A(t) = \sup_{c\leq s\leq t} a(s), \qquad\qquad
B(t) = \sup_{c\leq s\leq t} b(s), \\
K(t,s) = \sup_{s\leq \tau\leq t} k(\tau, s),
\end{gather*}
the inequality
\[
f\bigl(x(t)\bigr) \leq a(t) + b(t)\int_c^t k(t,s)g\bigl(x(s)\bigr) \d s,
\qquad c \leq t \leq d,
\]
implies the inequality
\[
x(t) \leq f^{-1}\oleft(\Omega^{-1}\oleft(
	\Omega\bigl(A(t)\bigr) + B(t) \int_c^t K(t,s) \d s \right)\right),
\qquad c \leq t \leq d' \leq d,
\]
where
\[
\Omega(u) = \int_\eps^u \frac{\d w}{g\bigl(f^{-1}(w)\bigr)},
\qquad \eps>0,u>0,
\]
and
\[
d' = \max\set*{c\leq t\leq d}
	{\Omega\bigl(A(t)\bigr) + B(t) \int_c^t K(t,s) \d s \leq \Omega\bigl(f(\infty)\bigr)}.
\]
\end{theorem}

\begin{proposition}[Uniform bound of the supports]\label{lem:supp-bound}
Let $v$, $V$, $W$ satisfy the assumptions \autoref{as:assumptions} and let $\bar\rho$ be the piecewise constant density associated to particles $X=(x_0,\dots,x_N)$ solving either \eqref{eq:ode-integrated-interaction} or \eqref{eq:ode-sampled-interaction}, with initial conditions $-S_0 \leq x_0(0) < x_N(0) \leq S_0$ for some $S_0>0$.

Then there exists an increasing function $S \in C\bigl([0,\infty);[0,\infty)\bigr)$ independent of $N$ such that
\[
-S(t) \leq x_0(t) \leq x_N(t) \leq S(t), \qquad \forall t\in[0,\infty),
\]
hence $\supp\bigl(\bar\rho(t)\bigr) \subseteq [-S(t),S(t)]$.
\end{proposition}

\begin{proof}
Consider the last particle $x_N$. We have
\[
\begin{split}
x_N'(t) &= v_N [V(t,x_N(t)) + W'*\bar\rho(t,x_N(t))] \\
&\leq V_+(t,x_N(t)) + \norm{W'}_{L^\infty([-x_N(t)+x_0(t),x_N(t)-x_0(t)])}.
\end{split}
\]
A similar estimate is valid for $x_0'(t)$, hence if we define $s(t)=\max\{x_N(t)_+,x_0(t)_-\}$, using \eqref{eq:mild-growth} we have
\[
s'(t) \leq F(t)\lambda\bigl(s(t)\bigr) + \norm{W'}_{L^\infty([-2s(t),2s(t)])} \leq 2F(t)\lambda\bigl(s(t)\bigr).
\]
By the assumptions on $F$ and $\lambda$ we can apply Bihari's estimate and deduce that
\[
s(t) \leq S(t) = \Lambda^{-1}(\Lambda(S_0) + A(t))
\]
where $A$ is the primitive of $F$.
\end{proof}

\begin{proposition}[Uniform bound of the density]\label{lem:max-princ}
Let $v$, $V$, $W$ satisfy the assumptions \autoref{as:assumptions} and let $\bar\rho$ be the piecewise constant density associated to particles $X=(x_0,\dots,x_N)$ solving either \eqref{eq:ode-integrated-interaction} or \eqref{eq:ode-sampled-interaction}, with initial conditions $-S_0 \leq x_0(0) < x_N(0) \leq S_0$ and $\bar\rho(0)\leq R_0$ for some $S_0,R_0>0$.

Then there exists an increasing function $R\in C\bigl([0,\infty);[0,\infty)\bigr)$ independent of $N$ such that
$\rho_i(t) \leq R(t)$ for every $t\in[0,\infty)$ and $i=1,\dots,N$, or equivalently $\bar\rho(t)\leq R(t)$.
\end{proposition}

\begin{proof}
At a fixed $t$, let $i$ be the index such that $\rho_i(t)=\max_j\rho_j(t)$. We can compute the derivative
\[
\begin{split}
\rho_i'(t) &= -N\rho_i(t)^2[x_i'(t)-x_{i-1}'(t)] \\
&= -N\rho_i^2[(v_i-v(\rho_i))U_i-(v_{i-1}-v(\rho_i))U_{i-1}]
	-N\rho_i^2v(\rho_i)(U_i-U_{i-1}).
\end{split}
\]
If $U_i<0$, then $v_i-v(\rho_i)=0$. If instead $U_i>0$, then $v_i-v(\rho_i)=v(\rho_{i+1})-v(\rho_i)\geq0$ because $\rho_i\geq\rho_{i+1}$. A similar consideration for $U_{i-1}$ leads to conclude that the first term is negative, thus $\rho_i'(t)\leq-N\rho_i^2v(\rho_i)(U_i-U_{i-1})$.

Let us assume that \eqref{eq:assumption-v-decay} holds. Thanks to \eqref{eq:lipschitz-bound} we have
\[
\begin{split}
&\rho_i'(t) \\
&\leq N\rho_i^2v(\rho_i)\bigl(
	\norm{V'}_{L^\infty([x_0,x_N])} + \norm{W''}_{L^\infty([x_0-x_N,x_N-x_0])}
	+ 2 \norm{W'}_{L^\infty([x_0-x_N,x_N-x_0])}\rho_i \bigr)
		(x_i-x_{i-1}) \\
&\leq \norm{v}_\infty \bigl(
	\norm{V'}_{L^\infty([x_0,x_N])} + \norm{W''}_{L^\infty([x_0-x_N,x_N-x_0])} \bigr) \rho_i
	+ 2\norm{W'}_{L^\infty([x_0-x_N,x_N-x_0])}\rho_i^2v(\rho_i) \\
&\leq 2\norm{v}_\infty F(t)G\bigl(2S(t)\bigr) \rho_i
	+ 2F(t)G\bigl(2S(t)\bigr) \rho_i^2v(\rho_i) .
\end{split}
\]
If we define $r(t)=\max_i\rho_i(t)=\norm{\bar\rho(t)}_\infty$, then from the previous estimate we get
\[
r'(t) \leq h(t)r(t) + k(t) r(t)^2v(r(t)),
\]
for some increasing functions $h,k$.
We are again in a position to apply Bihari's estimate and deduce
\[
r(t)\leq R(t) = \Omega^{-1}[\Omega(H(t))+K(t)],
\]
where $H$, $K$ and $\Omega$ are the primitives of $h$, $k$ and $1/(r^2v(r))$ respectively.

Let us now assume that \eqref{eq:assumption-W-repulsive} holds instead, with no longer any requirement about the decay of $v$. In this case, we proceed in a similar way but instead of directly applying the Lipschitz estimate \eqref{eq:lipschitz-bound} we must derive a better one-sided inequality in which we treat separately the interaction between the two consecutive particles, which has a sign thanks to the repulsive nature of $W$. In particular by \eqref{eq:assumption-W-repulsive} we have either $x_i-x_{i-1}>h$ or $W'(x_i-x_{i-1})-W'(x_{i-1}-x_i)\leq 0$. In the latter case, dealing with the sampled interaction, we find
\[
\begin{split}
U_i-U_{i-1}
&= V(x_i) - V(x_{i-1}) - \frac1N \sum_{j\neq i-1,i} [W'(x_i-x_j)-W'(x_{i-1}-x_j)]
	\spliteq
	- \frac1N[W'(x_i-x_{i-1})-W'(x_{i-1}-x_i)] \\
&\geq -\abs{V_i - V_{i-1}} - \frac1N \sum_{j\neq i-1,i} \abs{W'(x_i-x_j)-W'(x_{i-1}-x_j)} \\
&\geq - \bigl(
	\norm{V'}_{L^\infty([x_0,x_N])} + \norm{W''}_{L^\infty([x_0-x_N,x_N-x_0])}
	\bigr) (x_i-x_{i-1}) .
\end{split}
\]
Therefore, either $\rho_i(t) \leq \frac{1}{Nh}$, or we get the inequality
\[
\begin{split}
\rho_i'(t)
&\leq N\rho_i^2 v(\rho_i) \bigl(
	\norm{V'}_{L^\infty([x_0,x_N])} + \norm{W''}_{L^\infty([x_0-x_N,x_N-x_0])}
	\bigr) (x_i-x_{i-1}) \\
&\leq \norm{v}_\infty \bigl(
	\norm{V'}_{L^\infty([x_0,x_N])} + \norm{W''}_{L^\infty([x_0-x_N,x_N-x_0])}
	\bigr) \rho_i,
\end{split}
\]
from which we conclude by Gronwall instead of Bihari, using \autoref{as:assumptions} to estimate the norms.
The case with the integrated interaction is treated in the exact same way, by exploiting the sign of the term
\[
\int_{x_{i-1}}^{x_i} [W'(x_i-x)-W'(x_{i-1}-x)] \rho_i \d x
\]
when the particles $x_{i-1}$ and $x_i$ are near each other.
\end{proof}

\begin{remark}\label{rmk:discrete-derivatives}
Let $X$ be either $X^I$ or $X^S$ solving \eqref{eq:ode-integrated-interaction} or \eqref{eq:ode-sampled-interaction} respectively and let $U_i$ be either $\bar U_i$ or $\dot U_i$ accordingly. Then, as a consequence of \autoref{as:assumptions}, \autoref{lem:supp-bound} and \autoref{lem:max-princ}, the estimates in \autoref{lem:discrete-derivatives} take the form
\begin{align*}
\abs{U_i-U_{i-1}} &\leq 2 F(t)G\bigl(2S(t)\bigr) [1+R(t)] (x_i-x_{i-1}), \\
\abs{\dot U_{i+1}-2\dot U_i+\dot U_{i-1}}
&\leq \frac{3F(t)G\bigl(2S(t)\bigr)}{N\rho_i\rho_{i+1}} \left(
	\abs{\rho_{i+1}-\rho_i} + \frac{2R(t)}{N} \right)
	+ \frac{R(t)}{N\rho_i\rho_{i+1}} \int_{x_{i-1}}^{x_{i+1}} \abs{V''(x)}\d x
	\spliteq\notag
	+ \frac1{N^2\rho_i} \sum_{j\neq i, i\pm 1}
		\int_{x_{i-1}}^{x_{i+1}} \abs{W'''(y-x_j)} \d y , \\
\abs{\bar U_{i+1}-2\bar U_i+\bar U_{i-1}}
&\leq \frac{3F(t)G\bigl(2S(t)\bigr)}{N\rho_i\rho_{i+1}} \left(
	\abs{\rho_{i+1}-\rho_i} + \frac{2R(t)}{N} \right)
	+ \frac{R(t)}{N\rho_i\rho_{i+1}} \int_{x_{i-1}}^{x_{i+1}} \abs{V''(x)}\d x
	\spliteq\notag
	+ \frac1{N\rho_i}
		\int_{x_{i-1}}^{x_{i+1}} \abs{W'''}*\bar\rho(y) \d y .
\end{align*}
Similarly, the estimates in \autoref{lem:comparison} take the form
\begin{align*}
\abs{\bar U_i-\dot U_i}
&\leq \frac2N F(t)G\bigl(2S(t)\bigr)S(t), \\
\abs{(\bar U_i-\dot U_i) - (\bar U_{i-1}-\dot U_{i-1})}
&\leq \frac4N F(t)G\bigl(2S(t)\bigr) (x_i-x_{i-1}) .
\end{align*}
\end{remark}

The estimate for the total variation of a piecewise constant density relies on the following lemma, which provides the negativity of the error term appearing when replacing the mobilities of two consecutive particles with the one evaluated at the intermediate density. Together with the proof of the discrete entropy inequality \eqref{eq:discrete-entropy}, this is the only step where the monotonicity of $v$ is required.

\begin{lemma}\label{lem:negative-middle-term}
Letting $\sigma_i=\sign(\rho_{i+1}-\rho_i)$ and $\mu_i=\sigma_i-\sigma_{i-1}$, we have
\[
\mu_i \bigl[ \bigl(v_i-v(\rho_i)\bigr)U_i
	- \bigl(v_{i-1}-v(\rho_i)\bigr)U_{i-1} \bigr] \leq 0 .
\]
\end{lemma}

\begin{proof}
If $\rho_{i-1}<\rho_i<\rho_{i+1}$ or $\rho_{i-1}>\rho_i>\rho_{i+1}$, then $\mu_i=0$ and the inequality is satisfied. We need therefore to consider the cases $\rho_i\leq\rho_{i\pm1}$ and $\rho_i\geq\rho_{i\pm1}$. Let's investigate the former, as the latter is completely analogous. We have that $\sigma_i\leq0\leq\sigma_{i+1}$, hence $\mu_i\geq0$. If $U_i\geq0$, then $v_i=v(\rho_{i+1})\leq v(\rho_i)$, otherwise $v_i=v(\rho_i)$. In both cases we have that $\bigl(v_i-v(\rho_i)\bigr)U_i\leq0$. Similarly, $\bigl(v_{i-1}-v(\rho_i)\bigr)U_{i-1}\geq0$. The stated inequality is therefore satisfied.
\end{proof}

\begin{proposition}[Total variation growth]\label{prop:total-variation}
Given $N\in\setN$, let $\bar\rho$ be either
\begin{itemize}
\item the piecewise constant density associated to the particles $x^I=(x^I_0,\dots,x^I_N)$ solving \eqref{eq:ode-integrated-interaction}, driven by the total velocity field $v\bar U$;
\item the piecewise constant density associated to the particles $x^S=(x^S_0,\dots,x^S_N)$ solving \eqref{eq:ode-sampled-interaction}, driven by the total velocity field $v\dot U$;
\end{itemize}
with the assumptions of \autoref{lem:supp-bound}.
Then there exists an increasing function $B:[0,\infty)\to[0,\infty)$ independent of $N$ such that
\[
\TV\bigl(\bar\rho(t)\bigr) \leq B(t).
\]
\end{proposition}

\begin{proof}
Letting $\sigma_i=\sign(\rho_{i+1}-\rho_i)$, $\mu_i=\sigma_i-\sigma_{i-1}$ and using the fact that
\[
\rho_i' = -N\rho_i^2(x_i'-x_{i-1}'),
\]
since $\abs\plchldr$ is Lipschitz, by the chain rule we can compute
\[
\begin{split}
\frac\d{\d t}\TV\bigl(\bar\rho(t)\bigr)
&= \sum_{i=0}^N \sigma_i(\rho_{i+1}'-\rho_i')
 = -\sum_{i=1}^N \mu_i\rho_i'
 = \sum_{i=1}^N \mu_i N\rho_i^2(x_i'-x_{i-1}') \\
&= \sum_{i=1}^N \mu_i N\rho_i^2(v_iU_i-v_{i-1}U_{i-1}) \\
&= \sum_{i=1}^N \mu_i N\rho_i^2\bigl[ v(\rho_i)(U_i-U_{i-1})
	+ \bigl(v_i-v(\rho_i)\bigr)U_i - \bigl(v_{i-1}-v(\rho_i)\bigr)U_{i-1} \bigr] \\
&\leq \sum_{i=1}^N \mu_i N\rho_i^2 v(\rho_i)(U_i-U_{i-1})
 = \sum_{i=1}^N \mu_i \rho_i I_i ,
\end{split}
\]
where $I_i = N\rho_i v(\rho_i)(U_i-U_{i-1}) = N m(\rho_i)(U_i-U_{i-1})$ and the inequality is true because of \autoref{lem:negative-middle-term}.
Summing again by parts, we obtain
\begin{equation}\label{eq:TV-split}
\begin{split}
\frac\d{\d t}\TV\bigl(\bar\rho(t)\bigr)
&\leq \sum_{i=1}^N \mu_i \rho_i I_i
 = \sum_{i=1}^N (\sigma_i-\sigma_{i-1}) \rho_i I_i \\
&= \sum_{i=0}^N \sigma_i(\rho_iI_i-\rho_{i+1}I_{i+1}) \\
&= \sum_{i=0}^N \sigma_i(\rho_i-\rho_{i+1})I_{i+1}
	+\sum_{i=0}^N \sigma_i\rho_i(I_i-I_{i+1})
\end{split}
\end{equation}

Thanks to \eqref{eq:lipschitz-bound} with \autoref{rmk:discrete-derivatives}, we can bound
\[
\abs{I_i}
\leq 2\norm{v}_\infty N\rho_i F(t)G\bigl(2S(t)\bigr)[1+R(t)](x_i-{x_{i-1}})
\leq 2\norm{v}_\infty F(t)G\bigl(2S(t)\bigr)[1+R(t)] ,
\]
hence
\begin{equation}\label{eq:TV-term1}
\abs*{\sum_{i=0}^N \sigma_i (\rho_i-\rho_{i+1}) I_{i+1}}
\leq 2\norm{v}_\infty F(t)G\bigl(2S(t)\bigr)[1+R(t)] \TV\bigl(\bar\rho(t)\bigr).
\end{equation}
On the other hand, observe that
\[
\begin{split}
I_{i+1}-I_i
&= Nm(\rho_{i+1})(U_{i+1}-U_i)-Nm(\rho_i)(U_i-U_{i-1}) \\
&= Nm(\rho_{i+1})(U_{i+1}-2U_i+U_{i-1})
	+ N\bigl(m(\rho_{i+1})-m(\rho_i)\bigr)(U_i-U_{i-1}) ,
\end{split}
\]
therefore, for the integrated scheme, thanks to \eqref{eq:lipschitz-bound} and \eqref{eq:second-derivative-bound-integrated} with \autoref{rmk:discrete-derivatives}, we have
\begin{equation}\label{eq:TV-term2}
\begin{split}
\abs*{\sum_{i=0}^N \sigma_i\rho_i(I_i-I_{i+1})}
&\leq \sum_{i=1}^N \rho_i \abs{I_{i+1}-I_i} \\
&\leq N \sum_{i=1}^N \rho_i [ m(\rho_{i+1})\abs{\bar U_{i+1}-2\bar U_i+\bar U_{i-1}}
	+ \bigl(m(\rho_{i+1})-m(\rho_i)\bigr)\abs{\bar U_i-\bar U_{i-1}} ] \\
&\leq N \sum_{i=1}^N \rho_i \rho_{i+1} v(\rho_{i+1}) \Biggl\{
	\frac{3F(t)G\bigl(2S(t)\bigr)}{N\rho_i\rho_{i+1}} \left(
	\abs{\rho_{i+1}-\rho_i} + \frac{2R(t)}{N} \right)
		\spliteq\qquad
	+ \frac{R(t)}{N\rho_i\rho_{i+1}} \int_{x_{i-1}}^{x_{i+1}} \abs{V''(x)}\d x
	+ \frac1{N\rho_i}
		\int_{x_{i-1}}^{x_{i+1}} \abs{W'''}*\bar\rho(y) \d y \Biggr\}
	\spliteq
	+ N \sum_{i=1}^N \rho_i \norm{m'}_{L^\infty([0,R(t)])}
	\abs{\rho_{i+1}-\rho_i} 2 F(t) G\bigl(2S(t)\bigr)[1+R(t)](x_i-x_{i-1}) \\
&\leq 3F(t)G\bigl(2S(t)\bigr) \norm{v}_\infty \TV\bigl(\bar\rho(t)\bigr)
	+ 6F(t)G\bigl(2S(t)\bigr)R(t) \norm{v}_\infty
	\spliteq
	+ 2R(t) \norm{v}_\infty \norm{V''}_{L^1([-S(t),S(t)])}
	+ 2R(t) \norm{v}_\infty \norm{\abs{W'''}*\bar\rho}_{L^1([-S(t),S(t)])}
	\spliteq
	+ 2F(t) G\bigl(2S(t)\bigr)[1+R(t)]
	\bigl( \norm{v}_\infty + R(t)\norm{v'}_{L^\infty([0,R(t)])} \bigr)
		\TV\bigl(\bar\rho(t)\bigr) \\
&\leq 6F(t)G\bigl(2S(t)\bigr)R(t) \norm{v}_\infty
	+ 2F(t)G\bigl(S(t)\bigr)R(t) \norm{v}_\infty
	\spliteq
	+ 2R(t) \norm{v}_\infty \norm{W'''}_{L^1[-2S(t),2S(t)]}
	\spliteq 
	+5F(t) G\bigl(2S(t)\bigr)[1+R(t)]
		[\norm{v}_\infty + R(t)\norm{v'}_{L^\infty([0,R(t)])}] \TV\bigl(\bar\rho(t)\bigr) \\
&\leq 10F(t)G\bigl(2S(t)\bigr)R(t) \norm{v}_\infty
	\spliteq
	+5F(t) G\bigl(2S(t)\bigr)[1+R(t)]
		\bigl[\norm{v}_\infty + R(t)G\bigl(R(t)\bigr)\bigr] \TV\bigl(\bar\rho(t)\bigr) .
\end{split}
\end{equation}
For the sampled scheme, we can perform a similar computation using \eqref{eq:second-derivative-bound-sampled} with \autoref{rmk:discrete-derivatives}, where the only difference is in the term
\[
\begin{split}
N \sum_{i=1}^N & \rho_i \rho_{i+1} v(\rho_{i+1})
	\frac1{N^2\rho_i} \sum_{j\neq i, i\pm 1}
		\int_{x_{i-1}}^{x_{i+1}} \abs{W'''(y-x_j)} \d y \\
&\leq 2 R(t) \norm{v}_\infty \frac1N \sum_{j\neq i,i\pm1}
	\int_{x_0}^{x_N} \abs{W'''(y-x_j)} \d y \\
&\leq 2 R(t) F(t) G\bigl(2S(t)\bigr) \norm{v}_\infty,
\end{split}
\]
which leads to the same estimate as for the integrated scheme.


In conclusion, putting together \eqref{eq:TV-split}, \eqref{eq:TV-term1} and \eqref{eq:TV-term2}, we deduce
\[
\frac{\d}{\d t}\TV\bigl(\bar\rho(t)\bigr)
\leq P(t) + Q(t) \TV\bigl(\bar\rho(t)\bigr) ,
\]
for some finite functions $P$ and $Q$ independent of $N$, hence by Gronwall we deduce that $\TV\bigl(\bar\rho(t)\bigr) \leq B(t)$ for some increasing function $B:[0,\infty)\to[0,\infty)$.
\end{proof}

We conclude this section showing that the sequence of $\bar\rho^N$ constructed with either scheme \eqref{eq:ode-integrated-interaction} or \eqref{eq:ode-sampled-interaction} is compact in $L^1_\loc$.

\begin{theorem}\label{thm:compactness}
Let $v$, $V$ and $W$ satisfy the assumptions \autoref{as:assumptions} and let $\rho_0\in\Prob(\setR)\cap L^\infty(\setR)\cap BV(\setR)$ with the bound $\rho_0\leq R_0$ and compact support $\supp(\rho_0)\subseteq[-S_0,S_0]$.

Let $\bar\rho^N_0$ satisfy the initial conditions \eqref{eq:initial-conditions} and $\bar\rho^N$ be the piecewise constant density associated to the particles solving either \eqref{eq:ode-integrated-interaction} or \eqref{eq:ode-sampled-interaction}.

Then there exists a non negative probability density $\rho\in L^\infty_\loc\bigl([0,\infty)\times\setR\bigr)$ such that up to a subsequence $\bar\rho^N$ converges to $\rho$ in the strong topology of $L^p([0,T] \times \setR)$ for every $T>0$ and every $p\in[1,\infty)$. Moreover, for almost every $t\in[0,\infty)$, one has $\bar\rho^N(t,\plchldr)\to\rho(t,\plchldr)$ in $L^p(\setR)$ for every $p\in[1,\infty)$.
\end{theorem}

\begin{proof}
The proof relies on a very powerful generalization of Aubin-Lions Lemma as obtained in \cite[Theorem 2]{RossiSavare} which ensures that the sequence of $\bar\rho^N$ is strongly relatively compact in $L^1_\loc([0,\infty)\times\setR)$ as soon as 
\begin{gather}
\sup_N \int_0^T \bigl[ \abs{\supp\bigl(\bar\rho^N(t)\bigr)} + \TV\bigl(\bar\rho^N(t)\bigr) \bigr] \dt < \infty, \label{eq:uniform-bounds} \\
\sup_N W_1\bigl(\bar\rho^N(t,\plchldr), \bar\rho^N(s,\plchldr)\bigr) < C \abs{t-s}
	\quad\text{for any $s,t \in [0,T]$}, \label{eq:W1-lip}
\end{gather}
where $W_1$ denotes the 1-Wasserstein distance.

Notice that \eqref{eq:uniform-bounds} is an immediate consequence of \autoref{lem:supp-bound} and \autoref{prop:total-variation}.

Let us now focus on the proof of \eqref{eq:W1-lip}. The optimal transport map $\Xi$ from $\bar\rho^N(s,\plchldr)$ to $\bar\rho^N(t,\plchldr)$ is increasing and piecewise affine, mapping every interval $[x^N_{i-1}(s),x^N_i(s)]$ onto $[x^N_{i-1}(t),x^N_i(t)]$, given by the formula
\[
\Xi(x) = \sum_{i=1}^N \left(
	\frac{x-x^N_{i-1}(s)}{x^N_i(s)-x^N_{i-1}(s)} x^N_i(t)
	+\frac{x^N_i(s)-x}{x^N_i(s)-x^N_{i-1}(s)} x^N_{i-1}(t)
	\right) \bm1_{[x^N_{i-1}(s),x^N_i(s)]}(x).
\]
Notice in particular that for $x\in[x^N_{i-1}(s),x^N_i(s)]$ we have
\[
\begin{split}
\abs{\Xi(x)-x}
&= \abs*{
	\frac{x-x^N_{i-1}(s)}{x^N_i(s)-x^N_{i-1}(s)} [x^N_i(t)-x^N_i(s)]
	+\frac{x^N_i(s)-x}{x^N_i(s)-x^N_{i-1}(s)} [x^N_{i-1}(t)-x^N_{i-1}(s)]} \\
&\leq \abs{x^N_{i-1}(t)-x^N_{i-1}(s)} + \abs{x^N_i(t)-x^N_i(s)},
\end{split}
\]
therefore
\[
\begin{split}
W_1\bigl(\bar\rho^N(t,\plchldr), \bar\rho^N(s,\plchldr)\bigr)
&\leq \frac2N \sum_{i=0}^N \abs{x^N_i(t)-x^N_i(s)}
= \frac2N \sum_{i=0}^N \abs*{\int_s^t (x^N_i)'(\tau) \d\tau} \\
&\leq \frac{\norm{v}_\infty}N \sum_{i=0}^N \int_s^t \abs*{U_i(\tau)} \d\tau \\
&\leq \norm{v}_\infty \int_s^t \bigl(
	\norm{V}_{L^\infty([-S(\tau),S(\tau)])} + \norm{W'}_{L^\infty([-2S(\tau),2S(\tau)])}
	\bigr) \d\tau \\
&\leq 2 \norm{v}_\infty F(T) G\bigl(2S(T)\bigr) \abs{t-s} .
\end{split}
\]
%
%

Since \eqref{eq:uniform-bounds} and \eqref{eq:W1-lip} are verified, from Aubin-Lions Lemma we get that, up to a subsequence, $\bar\rho^N\to\rho$ in $L^1\bigl([0,T]\times\setR\bigr)$ for every $T>0$.

In addition, since $\bar\rho^N,\rho\leq R(T)$, for every $p\in[1,\infty)$ we have that
\[
\norm{\bar\rho^N-\rho}_{L^p([0,T]\times\setR)}^p
\leq \norm{\bar\rho^N-\rho}_{L^\infty([0,T]\times\setR)}^{p-1}
	\norm{\bar\rho^N-\rho}_{L^1([0,T]\times\setR)}
\leq [2R(T)]^{p-1} \norm{\bar\rho^N-\rho}_{L^1([0,T]\times\setR)},
\]
hence there is also strong convergence in $L^p\bigl([0,T]\times\setR\bigr)$.

By a standard argument, up to a further subsequence, for almost every $t\in[0,T]$ we have $\bar\rho^N(t,\plchldr)\to\rho(t,\plchldr)$ in $L^p(\setR)$.
\end{proof}

\begin{corollary}\label{cor:continuity}
With the assumptions of \autoref{thm:compactness}, we have that $\rho\in C\bigl([0,\infty);L^p(\setR)\bigr)$ for every $p\in[1,\infty)$.
\end{corollary}

\begin{proof}
From the equi-Lipschitzianity \eqref{eq:W1-lip}, we deduce that also $[0,T]\to \bigl(\Prob(\setR),W_1\bigr):t\mapsto\rho(t,\plchldr)$ is Lipschitz, hence $\rho(t,\plchldr)$ is continuous in time with respect to the weak convergence of measures.

Let us observe that for $p\in(1,\infty)$ and $s\to t$ we have $\rho(s)\weakto\rho(t)$ in $L^p(\setR)$. In fact, given $\phi\in L^{p'}(\setR)$ and $\eps>0$ we can find $\psi\in C_0(\setR)\cap L^{p'}(\setR)$ such that $\norm{\phi-\psi}_{L^{p'}(\setR)}<\eps$. But then
\begin{multline*}
\abs*{\int_\setR \rho(s,x)\phi(x)\dx-\int_\setR \rho(t,x)\phi(x)\dx} \\
\leq \abs*{\int_\setR [\rho(s,x)-\rho(t,x)]\psi(x)\dx}
	+ \int_\setR \abs{\rho(s,x)-\rho(t,x)}\cdot\abs{\phi(x)-\psi(x)}\dx,
\end{multline*}
and the first integral goes to $0$ by weak convergence whereas the second can be estimated by H\"older
\[
\norm{\rho(s)-\rho(t)}_{L^p(\setR)} \norm{\phi-\psi}_{L^{p'}(\setR)}
\leq \eps \norm{\rho(s)-\rho(t)}_{L^1(\setR)}^{1/p}
	\norm{\rho(s)-\rho(t)}_{L^\infty(\setR)}^{1/p'}
\leq \eps 2 R(t)^{1/p'} .
\]

Moreover, the functions $[0,T]\to\setR:t\to\norm{\bar\rho^N(t,\plchldr)}_{L^p(\setR)}$ are equi-Lipschitz because, thanks to \eqref{eq:lipschitz-bound} with \autoref{rmk:discrete-derivatives} and \autoref{prop:total-variation}, we have
\[
\begin{split}
\abs*{\frac{\d}{\dt} \norm{\bar\rho^N(t,\plchldr)}_{L^p(\setR)}^p}
&= \abs*{\frac{\d}{\dt} \sum_{i=1}^N \rho_i^p(x_i-x_{i-1})}
= \abs*{\frac{\d}{\dt} \frac1N \sum_{i=1}^N \rho_i^{p-1}}
= \abs*{\frac{p-1}N \sum_{i=1}^N \rho_i^{p-2}\rho_i'} \\
&\leq (p-1) \sum_{i=1}^N \rho_i^p \abs{x_i'-x_{i-1}'}
= (p-1) \sum_{i=1}^N \rho_i^p \abs{v_iU_i-v_{i-1}U_{i-1}} \\
&\leq (p-1) \sum_{i=1}^N \rho_i^p \bigl(
	v_i\abs{U_i-U_{i-1}}+\abs{v_i-v_{i-1}}\cdot \abs{U_{i-1}} \bigr) \\
&\leq (p-1) \sum_{i=1}^N \rho_i^p \norm{v}_\infty
	2F(T)G\bigl(2S(T)\bigr)[1+R(T)](x_i-x_{i-1})
	\spliteq
	+ (p-1) \sum_{i=1}^N \rho_i^p \norm{v'}_{L^\infty([0,R(T)])} \bigl(
			\abs{\rho_{i+1}-\rho_i} + \abs{\rho_i-\rho_{i-1}} \bigr)
		2F(t)G\bigl(2S(t)\bigr) \\
&\leq 2(p-1)R(T)^{p-1}F(T)G\bigl(2S(T)\bigr)[1+R(T)]\norm{v}_\infty
	\spliteq
	+ 4(p-1)R(T)^pF(T)G\bigl(2S(T)\bigr)G\bigl(R(T)\bigr)B(T),
\end{split}
\]
which is independent of $N$; therefore the limit function $[0,T]\to\setR:t\mapsto\norm{\rho(t,\plchldr)}_{L^p(\setR)}$ is Lipschitz as well.

It follows that, for $s\to t$, $\rho(s)$ converges to $\rho(t)$ weakly in $L^p(\setR)$ and the $L^p$ norms converge as well. Since $L^p(\setR)$ is uniformly convex for $p\in(1,\infty)$, we have that $\rho(s)\to\rho(t)$ strongly in $L^p(\setR)$.

By the boundedness of the support, this implies the strong convergence in $L^1(\setR)$ too.
\end{proof}

\subsection{Existence, uniqueness and convergence to the entropy solution}
\label{sec:entropy-solution}

The proof of \autoref{thm existence of entropy solution} relies on several ingredients. We begin by presenting the tools used to prove the existence. The first proposition states an approximate entropic inequality for the discrete densities, whereas the second allows to pass to the limit this inequality.

\begin{proposition}\label{prop:discrete-entropy-condition}
For every $N\in\setN$, let $\bar\rho^N$ be either
\begin{itemize}
\item the piecewise constant density associated to the particles $x^I=(x^I_0,\dots,x^I_N)$ solving \eqref{eq:ode-integrated-interaction}, driven by the total velocity field $v\bar U$;
\item the piecewise constant density associated to the particles $x^S=(x^S_0,\dots,x^S_N)$ solving \eqref{eq:ode-sampled-interaction}, driven by the total velocity field $v\dot U$.
\end{itemize}

Then for every constant $c\geq0$, $T>0$ and non-negative test function $\phi\in C^\infty_c\bigl((0,T)\times\setR;[0,\infty)\bigr)$ we have
\begin{multline}\label{eq:discrete-entropy}
\int_0^\infty \int_\setR \left\{ \abs{\bar\rho^N-c}\partial_t\phi
	+\sign(\bar\rho^N-c) \bigl[\bigl(m(\bar\rho^N)-m(c)\bigr)\bar U^N\partial_x\phi
	-m(c)\partial_x\bar U^N\phi\bigr] \right\} \d x \d t \\
\geq - \frac1N H(T) (\norm{\partial_x\phi}_\infty + \norm{\phi}_\infty),
\end{multline}
where
\[
\bar U^N(t,x) = V(t,x) - (W'*\bar\rho^N)(t,x)
\]
and $H:[0,\infty)\to[0,\infty)$ is an increasing function depending on $F,G,R,S,B$.

Moreover, $\bar\rho^N$ satisfies an approximate continuity equation: for every $\phi\in C^\infty_c\bigl((0,T)\times\setR\bigr)$ we have
\begin{equation}\label{eq:approx-cont-eq}
\abs*{
	\int_0^\infty\int_\setR \bigl[
	\bar\rho^N\partial_t\phi+m(\bar\rho^N)\bar U^N\partial_x\phi \bigr]\dx\dt }
= O_\phi\oleft(\frac1N\right) .
\end{equation}
\end{proposition}


\begin{proposition}\label{prop:entropic-limit}
Let $v$, $V$ and $W$ satisfy \autoref{as:assumptions}. For every $n\in\setN$ let $\bar\rho^N:[0,\infty)\times\setR\to[0,\infty)$ be a curve of probability measures such that  $\bar\rho^N(t)\leq R(t)$, $\supp\bar\rho^N(t)\subset[-S(t),S(t)]$, $\TV\bigl(\bar\rho^N(t)\bigr)\leq B(t)$, where $R,S,B:(0,\infty)\to[0,\infty)$ are increasing functions.

Assume that, for every $T>0$, $\bar\rho^N$ converges to $\rho$ in $L^1([0,T]\times\setR)$ and that $\bar\rho^N$ satisfies \eqref{eq:discrete-entropy}.

Then $\rho$ satisfies \eqref{eq:continuity-equation} and \eqref{eq:entropy-condition}, hence it is both a weak solution and an entropy solution of \eqref{eq:conservation-law} according to \autoref{def:entropy-solution}.
\end{proposition}

On the other hand, the uniqueness in \autoref{thm existence of entropy solution} follows from a variation of the arguments developed in \cite{KarlsenRisebro} and \cite{Kruzkov}. More precisely, adapting \cite[Theorem 1.3]{KarlsenRisebro} we deduce the following statement. Given the strong similarity with \cite[Theorem 1.3]{KarlsenRisebro}, we will not give a fully detailed proof of the \autoref{thm:KR}, but will instead just sketch the main steps, pointing out the main differences.

\newcommand{\citekr}{\cite[Theorem 1.3]{KarlsenRisebro}}
\begin{theorem}[\citekr]\label{thm:KR}
Consider two fluxes $\setR_t\times\setR^n_x\times\setR_u\to\setR$ of the form $p(u)P(t,x)$ and $q(u)Q(t,x)$, with $p,q\in\Lip_\loc\bigl([0,\infty)\bigr)$ and $P,Q\in L^1_\loc\bigl([0,\infty);W^{1,1}_\loc(\setR)\bigr)$.
Let then $\rho,\sigma \in L^\infty_\loc\bigl([0,\infty);\Prob(\setR)\cap BV(\setR)\bigr)$ be ``almost entropy solutions'' of
\begin{align*}
\partial_t \rho + \div_x\bigl(p(\rho)P(t,x)\bigr) &\simeq 0, &
\partial_t \sigma + \div_x\bigl(q(\sigma)Q(t,x)\bigr) &\simeq 0,
\end{align*}
in the sense that
\begin{align*}
\int_0^\infty \int_\setR \bigl\{ \abs{\rho-c}\partial_t\phi
	+\sign(\rho-c) \bigl[\bigl(p(\rho)-p(c)\bigr)P\div_x\phi
&	-p(c)\div_x P\phi\bigr] \bigr\} \d x \d t \\
&\geq -\frac1N H(T)(\norm{\partial_x\phi}_\infty+\norm{\phi}_\infty), \\
\int_0^\infty \int_\setR \bigl\{ \abs{\sigma-c}\partial_t\phi
	+\sign(\sigma-c) \bigl[\bigl(q(\sigma)-q(c)\bigr)Q\div_x\phi
&	-q(c)\div_x Q\phi\bigr] \bigr\} \d x \d t \\
&\geq -\frac1N H(T)(\norm{\partial_x\phi}_\infty+\norm{\phi}_\infty),
\end{align*}
for every constant $c\geq 0$ and non-negative test function $\phi\in C^\infty_c\bigl((0,T)\times\setR;[0,\infty)\bigr)$, with
\begin{align*}
\supp\bigl(\rho(t)\bigr) &\subseteq [-S(t),S(t)], &
\rho(t) &\leq R(t), &
\TV\bigl(\rho(t)\bigr) &\leq B(t), \\
\supp\bigl(\sigma(t)\bigr) &\subseteq [-S(t),S(t)], &
\sigma(t) &\leq R(t), &
\TV\bigl(\sigma(t)\bigr) &\leq B(t),
\end{align*}
where $R,S,B,H:[0,\infty)\to[0,\infty)$ are increasing functions.

Then 
\begin{equation}\label{eq:KR}
\begin{split}
\norm{\rho(t)-\sigma(t)}_{L^1(\setR)} \Bigr\rvert_{t_1}^{t_2}
&\leq \int_{t_1}^{t_2} \int_{-S(t)}^{S(t)} \Bigl\{
	\norm{p-q}_{L^\infty([0,R(t)])} \abs{\div_x P}
	+ \norm{p}_{L^\infty([0,R(t)])} \abs{\div_x(P-Q)}
	\spliteq\qquad\qquad
	+ \norm{P-Q}_{L^\infty([-S(t),S(t)])} \norm{p'}_{L^\infty([0,R(t)])} \abs{\div_x\rho}
	\spliteq\qquad\qquad
	+ \norm{P}_{L^\infty([-S(t),S(t)])} \norm{p'-q'}_{L^\infty([0,R(t)])} \abs{\div_x\rho}
	\Bigr\} \dx \dt
	\spliteq
	+ \frac{4H(T)}N .
\end{split}
\end{equation}
\end{theorem}

\begin{remark}\label{rmk:stability}
In the setting of \autoref{thm:KR}, assume that $p=q=m$ and the fields $P$ and $Q$ are given by $P=V-W'*\bar\rho^M$ and $Q=V-W'*\bar\rho^N$, where $v$, $V$ and $W$ satisfy \autoref{as:assumptions}, $M\geq N$ are positive integers, and $\bar\rho^M$ and $\bar\rho^N$ are the piecewise constant densities associated to particles $X^M$ and $X^N$ solving either \eqref{eq:ode-integrated-interaction} or \eqref{eq:ode-sampled-interaction}. In this case, $\bar\rho^M$ and $\bar\rho^N$ satisfy the assumptions of \autoref{thm:KR} thanks to \autoref{prop:discrete-entropy-condition}, therefore \eqref{eq:KR} becomes
\[
\begin{split}
&\hspace{-0.5cm}\norm{\bar\rho^M(t)-\bar\rho^N(t)}_{L^1(\setR)} \\
&\leq \norm{\bar\rho^M(0)-\bar\rho^N(0)}_{L^1(\setR)} + \frac{8H(T)}N
	+ \norm{m}_{L^\infty([0,R(T)])} \int_0^t \int_{-S(s)}^{S(s)}
		\abs{\div_x(P-Q)} \dx \d s
	\spliteq
	+ \norm{m'}_{L^\infty([0,R(T)])} \int_0^t
		\norm{P-Q}_{L^\infty([-S(s),S(s)])} \TV\bigl(\rho(s)\bigr) \d s \\
&\leq \norm{\bar\rho^M(0)-\bar\rho^N(0)}_{L^1(\setR)} + \frac{8H(T)}N
	\spliteq
	+ \norm{m}_{L^\infty([0,R(T)])}
		\bigl(\norm{W''}_{L^1([-2S(T),2S(T)])}+\norm{w}_{L^\infty([0,T])} \bigr)	
		\int_0^t \norm{\bar\rho^M(s)-\bar\rho^N(s)}_{L^1(\setR)} \d s
	\spliteq
	+ \norm{m'}_{L^\infty([0,R(T)])} \norm{W'}_{L^\infty([-2S(T),2S(T)])} B(T)
		\int_0^t \norm{\bar\rho^M(s)-\bar\rho^N(s)}_{L^1(\setR)} \d s \\
&\leq \norm{\bar\rho^M(0)-\bar\rho^N(0)}_{L^1(\setR)} + \frac{8H(T)}N
	+ A(T) \int_0^t \norm{\bar\rho^M(s)-\bar\rho^N(s)}_{L^1(\setR)} \d s,
\end{split}
\]
for some increasing function $A$, depending on $F,G,R,S,B$.
By Gronwall lemma, this implies
\begin{equation}\label{eq:cauchy}
\norm{\bar\rho^M(t)-\bar\rho^N(t)}_{L^1(\setR)}
\leq \left( \norm{\bar\rho^M(0)-\bar\rho^N(0)}_{L^1(\setR)} + \frac{8H(T)}N \right)
	\exp\bigl(t A(T)\bigr), \qquad \forall t\in(0,T).
\end{equation}
Since we have the convergence of the initial datum $\bar\rho^M(0),\bar\rho^N(0)\to\rho_0$,
this proves that $\bigl(\bar\rho^N(t)\bigr)_{N\in\setN}$ is a Cauchy sequence in $L^1(\setR)$, uniformly for $t\in[0,T]$. This provides an alternative way to deduce the existence of a limit density $\rho$ which does not rely on the compactness argument of \autoref{thm:compactness}.
\end{remark}

%

We now proceed to prove the propositions and theorems presented in this section.

\begin{proof}[Proof of Proposition \ref{prop:discrete-entropy-condition}]
Let us fix the constant $c$ and the test function $\phi$. In order to keep the notation light, we write $\bar U=\bar U^N$ and we omit the dependence on $(t,x)$ whenever possible. Applying the definition of the particle approximations $\bar\rho^N$ to the left hand side of \eqref{eq:discrete-entropy}, we reduce to study the non negativity of the following quantity
\begin{equation}\label{eq:entropy-terms}
\begin{split}
&\mathrel{\phantom{=}} \int_0^T \int_{(x_0,x_N)^c} \Bigl[
	c \partial_t \phi(t,x) + m(c) \partial_x \bigl(\bar U^N(t,x) \phi(t,x)\bigr)
	\Bigr] \d x \d t
	\spliteq
+ \sum_{i=1}^N \int_0^T \int_{x_{i-1}}^{x_i}
	\Bigl[ \abs{\rho_i-c} \partial_t \phi(t,x)
	- \sign(\rho_i-c) m(\rho_i) \partial_x \bar U^N(t,x) \phi(t,x) \Bigr] \d x \d t
	\spliteq
+ \sum_{i=1}^N \int_0^T \Bigl[
	\sign(\rho_i-c) \bigl(m(\rho_i) - m(c)\bigr) \bigl(\bar U^N(t,x_i)\phi(t,x_i)
	- \bar U^N(t,x_{i-1})\phi(t,x_{i-1}) \bigr) \Bigr] \d t.
\end{split}
\end{equation}
The boundary term can be dealt with in the following manner
\begin{equation}\label{eq:boundary-term}
\begin{split} 
\int_0^T \int_{(x_0,x_N)^c} & \Bigl[
	c \partial_t \phi + m(c) \partial_x (\bar U \phi) \Bigr] \d x \d t \\
&= c \int_0^T \left[ \frac{\d}{\d t} \int_{(x_0,x_N)^c} \phi \d x
	+ \phi(x_N)x'_N - \phi(x_0)x'_0 \right] \d t
	\spliteq
	+ m(c) \int_0^T \bigl( \bar U(x_0)\phi(x_0) - \bar U(x_N)\phi(x_N) \bigr) \d t \\
&= c \int_0^T \Bigl[ 
	\phi(x_N)\bigl(x'_N - v(c)\bar U(x_N)\bigr)
	- \phi(x_0) \bigl(x'_0 - v(c)\bar U(x_0)\bigr) \Bigr] \d t
\end{split}
\end{equation}
Integrating twice by parts we can rewrite
\begin{equation*}
\begin{split}
\sum_{i=1}^N \int_0^T \int_{x_{i-1}}^{x_i} \abs{\rho_i-c}\partial_t\phi \d x
&= \sum_{i=1}^N \int_0^T \Biggl\{
	\frac{\d}{\d t} \left[ \int_{x_{i-1}}^{x_i} |\rho_i - c| \phi\d x \right]\d t
	- \int_{x_{i-1}}^{x_i} \sign(\rho_i -c) \rho'_i \phi \d x
	\spliteq
 - \abs{\rho_i - c} \bigl( \phi(x_i) x'_i - \phi(x_{i-1})x'_{i-1} \bigr)
 \Biggr\} \d t \\
&= \sum_{i=1}^N \int_0^T \sign(\rho_i-c) \Bigl[\rho_i (x'_i-x'_{i-1}) \phi(\bar x_i)
	\spliteq
	- (\rho_i - c) \big( \phi(x_i) x'_i - \phi(x_{i-1}) x'_{i-1} \bigr)\Bigr] \d t \\
&= \sum_{i=1}^N \int_0^T \sign(\rho_i-c) \rho_i (x'_i - x'_{i-1}) \bigl(\phi(\bar{x}_i) - \phi(x_i)\bigr) \d t
	\spliteq
- \sum_{i=1}^N \int_0^T \sign(\rho_i-c) \rho_i x'_{i-1} \bigl(\phi(x_i) - \phi(x_{i-1})\bigr) \d t
	\spliteq
+ c \sum_{i=1}^N \int_0^T \sign(\rho_i-c) \bigl(x'_i \phi(x_i) - x'_{i-1} \phi(x_{i-1})\bigr) \d t.
\end{split}
\end{equation*}
where $\bar x_i\in(x_{i-1},x_i)$ is such that $\phi(\bar x_i)=\dashint_{x_{i-1}}^{x_i} \phi(x) dx$.
Thanks to the above computation, we can combine the first term of the second line and the third line of \eqref{eq:entropy-terms} to get
\begin{equation}\label{eq:entropic-some-terms}
\begin{split}
&\sum_{i=1}^N \int_0^T \left\{
	\int_{x_{i-1}}^{x_i} \abs{\rho_i-c} \partial_t \phi \d x 
	+ \sign(\rho_i-c) \bigl(m(\rho_i) - m(c)\bigr) \bigl(\bar U(x_i)\phi(x_i) - \bar U(x_{i-1})\phi(x_{i-1}) \bigr) \right\} \d t \\
&= c \sum_{i=1}^N \int_0^T \sign(\rho_i-c)\Bigl[\bigl(x'_i-v(c)\bar U(x_i)\bigr)\phi(x_i)
	- \bigl(x'_{i-1}-v(c)\bar U(x_{i-1})\bigr)\phi(x_{i-1})\Bigr] \d t
	\spliteq
	+ \sum_{i=1}^N \int_0^T \sign(\rho_i-c)\rho_i \Bigl[
	(x'_i-x'_{i-1})\bigl(\phi(x_i)-\phi(x_{i-1})\bigr)
	-x'_{i-1}\bigl(\phi(x_i)-\phi(x_{i-1})\bigr)
		\spliteq \qquad\qquad
	+v(\rho_i)\bar U(x_{i-1})\bigl(\phi(x_i)-\phi(x_{i-1})\bigr)
	+v(\rho_i)\bigl(\bar U(x_i)-\bar U(x_{i-1})\bigr)\phi(x_i)\Bigr] \d t \\
&= - c \sum_{i=1}^{N-1} \int_0^T [\sign(\rho_{i+1}-c) - \sign(\rho_i-c)]
	\phi(x_i) \bigl(x'_i - v(c)\bar U(x_i)\bigr) \d t
	\spliteq
+ c \int_0^T \Bigl[ \sign(\rho_N-c) \phi(x_N) \bigl(x'_N - v(c)\bar U(x_N)\bigr)
	- \sign(\rho_1-c) \phi(x_0) \bigl(x'_0 - v(c)\bar U(x_0)\bigr) \Bigr] \d t
	\spliteq
+\sum_{i=1}^N \int_0^T \sign(\rho_i-c) \biggl\{
	m(\rho_i) \phi(x_i) (\bar U(x_i) - \bar U(x_{i-1}))
	\spliteq \qquad\qquad
	+ \rho_i \Bigl[ (x'_i-x'_{i-1}) (\phi(\bar x_i) - \phi(x_i))
	- \bigl(\phi(x_i) - \phi(x_{i-1})\bigr)\bigl(x'_{i-1} - v(\rho_i)\bar U(x_{i-1})\bigr) \Bigr] \biggr\} \d t .
\end{split}
\end{equation}
Applying \eqref{eq:boundary-term} and \eqref{eq:entropic-some-terms} to the appropriate terms in \eqref{eq:entropy-terms}, we can rewrite the latter as $I + \II + \III$ where
\begin{align*}
I &= - c \sum_{i=1}^{N-1} \int_0^T [\sign(\rho_{i+1}-c) - \sign(\rho_i-c)]
	\phi(x_i) \bigl(x'_i - v(c)\bar U(x_i)\bigr) \d t, \\
\II &= c \int_0^T \Bigl[
	\phi(x_N) \bigl(x'_N-v(c)\bar U(x_N)\bigr)\bigl(1+\sign(\rho_N-c)\bigr)
		\spliteq \qquad
	- \phi(x_0) \bigl(x'_0-v(c)\bar U(x_0)\bigr)\bigl(1+\sign(\rho_1-c)\bigr) \Bigr]\d t, \\
\III &= \sum_{i=1}^N \int_0^T \sign(\rho_i-c) \biggl\{
	m(\rho_i) \biggl[ \phi(x_i) \bigl(\bar U(x_i) - \bar U(x_{i-1})\bigr)
		-\int_{x_{i-1}}^{x_i} \partial_x\bar U \phi \d x \biggr]
	\spliteq \qquad\qquad
	+ \rho_i \Bigl[ (x'_i-x'_{i-1}) (\phi(\bar x_i) - \phi(x_i))
	- \bigl(\phi(x_i) - \phi(x_{i-1})\bigr)\bigl(x'_{i-1} - v(\rho_i)\bar U(x_{i-1})\bigr) \Bigr] \biggr\} \d t .
\end{align*}

We make the important remark that up to this point we have only performed algebraic manipulations and never used the actual ordinary differential equation solved by the particles $X=(x_0,\dots,x_N)$. In particular, all computations carried out so far apply to both schemes.

We now proceed to discuss the case in which $\bar\rho^N$ is associated to the particles $X^I$ solving \eqref{eq:ode-integrated-interaction}: in particular, we substitute $x'_i=v_i\bar U(x_i)$ in the terms $I,\II,\III$. The validity of \eqref{eq:discrete-entropy} follows once we prove that $I$ and $\II$ are non negative terms, while $\III$ goes to zero as $1/N$ when $N \to \infty$.

Observe, that the contribution of $\sign(\rho_{i+1}-c) - \sign(\rho_i-c)$ is non trivial only when $\rho_i \leq c \leq \rho_{i+1}$ or $\rho_i \geq c \geq \rho_{i+1}$ and, similarly, the contribution of $1+\sign(\rho_1-c)$ and $1+\sign(\rho_N-c)$ is non zero only when $\rho_1,\rho_N \geq c$. Then the non negativity of $I$ and $\II$ is an immediate consequence of the monotonicity of $v$.

Let us now focus on the term $\III$.
The first line can be estimated as
\[
\begin{split}
&\hspace{-.25cm} \abs*{\sum_{i=1}^N \int_0^T \sign(\rho_i-c)
	m(\rho_i) \biggl[ \phi(x_i) \bigl(\bar U(x_i) - \bar U(x_{i-1})\bigr)
		-\int_{x_{i-1}}^{x_i} \partial_x\bar U \phi \d x \biggr] \d t} \\
&\leq \norm{v}_\infty \int_0^T \sum_{i=1}^N \rho_i \abs*{\int_{x_{i-1}}^{x_i}
	\partial_x \bar U \bigl(\phi(x_i)-\phi(x)\bigr)\d x} \d t\\
&\leq \norm{v}_\infty \norm{\partial_x\phi}_\infty \int_0^T \rho_i(x_i-x_{i-1}) \sum_{x=1}^N
	\int_{x_{i-1}}^{x_i} \abs{\partial_x \bar U} \d x \d t \\
&\leq \frac1N \norm{v}_\infty \norm{\partial_x\phi}_\infty
	\int_0^T \int_{-S(t)}^{S(t)} \abs{\partial_x \bar U}\d x \d t \\
&\leq \frac1N \norm{v}_\infty \norm{\partial_x\phi}_\infty \int_0^T
	\norm{V' - (W''+w(t)\delta_0)*\bar\rho}_{L^1([-S(t),S(t)])} \d t \\
&\leq \frac1N \norm{v}_\infty \norm{\partial_x\phi}_\infty \int_0^T \left(
	\norm{V'}_{L^1([-S(t),S(t)])}
	+ \norm{W''}_{L^1([-2S(t),2S(t)])}
	+ \abs{w(t)}\right) \d t \\
&\leq \frac1N \norm{v}_\infty \norm{\partial_x\phi}_\infty \int_0^T
	F(t) \bigl[1+ G\bigl(S(t)\bigr)+G\bigl(2S(t)\bigr)\bigr] \d t \\
&\leq \frac1N \norm{v}_\infty \norm{\partial_x\phi}_\infty
	\norm{F}_{L^1([0,T])} \bigl[1+2G\bigl(2S(T)\bigr)\bigr].
\end{split}
\]
To estimate the second line of $\III$, first observe that thanks to \eqref{eq:lipschitz-bound} and \autoref{rmk:discrete-derivatives} we have
\[
\begin{split}
\abs{x'_i-x'_{i-1}}
&= \abs{v_i\bar U_i-v_{i-1}\bar U_{i-1}} \\
&\leq v(\rho_i) \abs{\bar U_i-\bar U_{i-1}} + \abs{v_i-v(\rho_i)}\cdot\abs{\bar U_i}
	+ \abs{v_{i-1}-v(\rho_i)}\cdot\abs{\bar U_{i-1}} \\
&\leq 2\norm{v}_\infty F(t)G\bigl(2S(t)\bigr)[1+R(t)] (x_i-x_{i-1})
		\spliteq
	+ \norm{v'}_{L^\infty([0,R(t)])} \norm{V-W'*\bar\rho}_{L^\infty([-S(t),S(t)])}
	\bigl(\abs{\rho_{i+1}-\rho_i}+\abs{\rho_i-\rho_{i-1}}\bigr) \\
&\leq 2\norm{v}_\infty F(t)G\bigl(2S(t)\bigr)[1+R(t)] (x_i-x_{i-1})
		\spliteq
	+ 2G\bigl(R(t)\bigr) F(t)G\bigl(2S(t)\bigr)
		\bigl(\abs{\rho_{i+1}-\rho_i}+\abs{\rho_i-\rho_{i-1}}\bigr)
\end{split}
\]
hence
\[
\begin{split}
&\abs*{\sum_{i=1}^N \int_0^T \sign(\rho_i-c)
	\rho_i \Bigl[ (x'_i-x'_{i-1}) (\phi(\bar x_i) - \phi(x_i))
	- \bigl(\phi(x_i) - \phi(x_{i-1})\bigr)\bigl(x'_{i-1} - v(\rho_i)\bar U(x_{i-1})\bigr) \Bigr] \d t} \\
&\leq \frac{\norm{\partial_x\phi}_\infty}N \sum_{i=1}^N \int_0^T \bigl(
	\abs{x'_i-x'_{i-1}} + \abs{v_{i-1}-v(\rho_i)}\cdot\abs{\bar U(x_{i-1})} \bigr)\d t \\
&\leq \frac{\norm{\partial_x\phi}_\infty}N \Biggl[
	2 \norm{v}_\infty \int_0^T F(t)G\bigl(2S(t)\bigr)[1+R(t)]S(t)\d t
	+8 \int_0^T F(t) G\bigl(R(t)\bigr) G\bigl(2S(t)\bigr) \TV\bigl(\bar\rho^N(t)\bigr) \d t
	\Biggr] \\
&\leq \frac2N \norm{\partial_x\phi}_\infty (\norm{v}_\infty+4) G\bigl(2S(T)\bigr)
	[1+R(T)+G\bigl(R(T)\bigr)][S(T)+B(T)] \norm{F}_{L^1([0,T])}.
\end{split}
\]
Putting together these last two estimates one sees that $\abs{\III} \leq \frac1N H(T)\norm{\partial_x\phi}_\infty$ for some increasing function $H:[0,\infty)\to[0,\infty)$ which can be expressed in terms of $F,G,R,S,B$.

Let us now move on to the case in which $\bar\rho^N$ is associated to the particles $X^S$ solving \eqref{eq:ode-sampled-interaction}. Adding and subtracting $\dot U$ where appropriate in the terms $I$, $\II$ and $\III$, we can rewrite \eqref{eq:entropy-terms} as $\tilde I + \tilde\II + \tilde\III + \tilde\IV$ where
\begin{align*}
\tilde I &= - c \sum_{i=1}^{N-1} \int_0^T [\sign(\rho_{i+1}-c) - \sign(\rho_i-c)]
	\phi(x_i) \bigl(x'_i - v(c)\dot U_i\bigr) \d t, \\
\tilde\II &= c \int_0^T \Bigl[
	\phi(x_N) \bigl(x'_N-v(c)\dot U_N\bigr)\bigl(1+\sign(\rho_N-c)\bigr)
		\spliteq \qquad
	- \phi(x_0) \bigl(x'_0-v(c)\dot U_0\bigr)\bigl(1+\sign(\rho_1-c)\bigr) \Bigr]\d t, \\
\tilde\III &= \sum_{i=1}^N \int_0^T \sign(\rho_i-c) \biggl\{
	m(\rho_i) \biggl[ \phi(x_i) \bigl(\bar U_i - \bar U_{i-1}\bigr)
		-\int_{x_{i-1}}^{x_i} \partial_x\bar U \phi \d x \biggr]
	\spliteq \qquad\qquad
	+ \rho_i \Bigl[ (x'_i-x'_{i-1}) (\phi(\bar x_i) - \phi(x_i))
	- \bigl(\phi(x_i) - \phi(x_{i-1})\bigr)\bigl(x'_{i-1} - v(\rho_i)\dot U_{i-1}\bigr) \Bigr] \biggr\} \d t, \\
\tilde\IV &= m(c) \sum_{i=1}^{N-1} \int_0^T [\sign(\rho_{i+1}-c) - \sign(\rho_i-c)]
	\phi(x_i) \bigl(\bar U_i-\dot U_i\bigr) \d t
		\spliteq
	- m(c) \int_0^T \Bigl[
	\phi(x_N) \bigl(\bar U_N-\dot U_N\bigr)\bigl(1+\sign(\rho_N-c)\bigr)
	- \phi(x_0) \bigl(\bar U_0-\dot U_0\bigr)\bigl(1+\sign(\rho_1-c)\bigr) \Bigr]\d t
		\spliteq
	+ \sum_{i=1}^N \int_0^T \sign(\rho_i-c) m(\rho_i) \bigl(\phi(x_i) - \phi(x_{i-1})\bigr)
		\bigl(\bar U_{i-1} - \dot U_{i-1}\bigr) \d t.
\end{align*}
Since in this case $x'_i=v_i\dot U_i$, the non-negativity of $\tilde I$ and $\tilde\II$ and the fact that $\abs{\tilde\III}\leq\frac1N H(T)\norm{\partial_x\phi}_\infty$ follow by the same arguments as before. What is left to show is that the remainder term $\tilde\IV$ goes to $0$ as well at the same rate $1/N$. Thanks to \autoref{lem:comparison} with \autoref{rmk:discrete-derivatives} we can estimate
\[
\begin{split}
\abs{\tilde\IV}
&\leq m(c) \sum_{i=1}^N \int_0^T \abs*{ \sign(\rho_i-c)
	\Bigl[\phi(x_i)\bigl(\bar U_i-\dot U_i\bigr)
	-\phi(x_{i-1})\bigl(\bar U_{i-1}-\dot U_{i-1}\bigr)\Bigr] } \d t
		\spliteq
	+m(c)\int_0^T \abs*{ \phi(x_N)(\bar U_N-\dot U_N) - \phi(x_0)(\bar U_0-\dot U_0)}\d t
		\spliteq
	+\sum_{i=1}^N \int_0^T \abs*{
		\sign(\rho_i-c) m(\rho_i) \bigl(\phi(x_i) - \phi(x_{i-1})\bigr)
		\bigl(\bar U_{i-1} - \dot U_{i-1}\bigr) } \d t \\
&\leq \sum_{i=1}^N \int_0^T
	[m(c)+m(\rho_i)] \cdot \abs{\phi(x_i)-\phi(x_{i-1})}
		\cdot \abs{\bar U_{i-1}-\dot U_{i-1}} \d t
		\spliteq
	+ m(c) \sum_{i=1}^N \int_0^T \phi(x_i)
		\abs{(\bar U_i-\dot U_i) - (\bar U_{i-1}-\dot U_{i-1})} \d t
		\spliteq
	+m(c)\int_0^T \Bigl[
		\phi(x_N)\abs{\bar U_N-\dot U_N} + \phi(x_0)\abs{\bar U_0-\dot U_0} \Bigr] \d t \\
&\leq \frac1N \sum_{i=1}^N \int_0^T \bigl[m(c)+R(t)\norm{v}_{L^\infty}\bigr]
	\norm{\partial_x\phi}_\infty (x_i-x_{i-1}) F(t) G\bigl(2S(t)\bigr) S(t) \d t
		\spliteq
	+ \frac4N m(c) \norm{\phi}_{L^\infty} \sum_{i=1}^N
		\int_0^T F(t)G\bigl(2S(t)\bigr) (x_i-x_{i-1}) \d t
		\spliteq
	+ \frac2N m(c) \norm{\phi}_{L^\infty} \int_0^T F(t)G\bigl(2S(t)\bigr) S(t) \d t \\
&\leq \frac1N \bigl(\norm{\partial_x\phi}_\infty+6\norm{\phi}_{L^\infty}\bigr)
	\bigl(m(c)+R(T)\norm{v}_{L^\infty}\bigr)
	G\bigl(2S(T)\bigr) \bigl(S(T)+S(T)^2\bigr) \norm{F}_{L^1([0,T])} .
\end{split}
\]
This concludes the proof of \eqref{eq:discrete-entropy}.

Let us now check that $\bar\rho^N$ satisfies the approximate continuity equation \eqref{eq:approx-cont-eq}. Fix $\phi\in C^\infty_c\bigl((0,T)\times\setR;[0,\infty)\bigr)$. Taking $c=0$ we get from \eqref{eq:discrete-entropy}
\begin{equation}\label{eq:cont-eq-side1}
\int_0^T\int_\setR \left\{
	\bar\rho^N \partial_t\phi + m(\bar\rho^N)\bar U^N\partial_x\phi \right\} \dx\dt
\geq O_\phi\oleft(\frac1N\right),
\end{equation}
whereas taking $c=R(T)$ and recalling that $\bar\rho^N(t,\plchldr)<R(T)$ for $t\leq T$ we get
\[
\int_0^T\int_\setR \left\{
	[R(T)-\bar\rho^N]\partial_t\phi
	- m(\bar\rho^N)\bar U^N\partial_x\phi
	+ m\bigl(R(T)\bigr)\partial_x(\bar U^N\phi) \right\} \dx\dt
\geq O_\phi\oleft(\frac1N\right),
\]
which implies
\begin{equation}\label{eq:cont-eq-side2}
\int_0^T\int_\setR \left\{
	-\bar\rho^N\partial_t\phi
	- m(\bar\rho^N)\bar U^N\partial_x\phi \right\} \dx\dt
\geq O_\phi\oleft(\frac1N\right)
\end{equation}
because the terms with $R(T)$ and $m\bigl(R(T)\bigr)$ give a zero contribution to the integral. Combining \eqref{eq:cont-eq-side1} and \eqref{eq:cont-eq-side2} we deduce the validity of \eqref{eq:approx-cont-eq} for non-negative test functions. This requirement can then be dropped by writing an arbitrary $\phi\in C^\infty_c\bigl((0,T)\times\setR\bigr)$ as a difference $\phi=\phi_1-\phi_2$ with $\phi_1,\phi_2\in C^\infty_c\bigl((0,T)\times\setR;[0,\infty)\bigr)$.
\end{proof}

In the proof of \autoref{prop:entropic-limit} we reuse as much as possible the argument from \cite{DiFrancescoFagioliRadici}, focusing the attention primarily on the term \eqref{eq:limit2} which is more problematic to pass to the limit.

\begin{proof}[Proof of \autoref{prop:entropic-limit}]
Up to a subsequence, we may assume that $\bar\rho^N\to\rho$ almost everywhere.

Let us first show that $\rho$ solves the continuity equation \eqref{eq:continuity-equation}.
For a fixed test function $\phi\in C^\infty_c\bigl((0,T)\times\setR\bigr)$, thanks to \autoref{prop:discrete-entropy-condition} we know that $\bar\rho^N$ solves the approximate continuity equation \eqref{eq:approx-cont-eq}, therefore we just need to pass it to the limit. The term $\bar\rho^N\partial_t\phi$ is obvious. For the other one we use the fact that
\[
\begin{split}
\int_0^T & \int_\setR \abs{m(\bar\rho^N)\bar U^N-m(\rho)U} \dx \dt \\
&\leq \int_0^T \int_{-S(t)}^{S(t)} \abs{m(\bar\rho^N)(\bar U^N-U)} \dx \dt
	+ \int_0^T \int_{-S(t)}^{S(t)} \abs{[m(\bar\rho^N)-m(\rho)]U} \dx \dt \\
&\leq \norm{v}_\infty R(T) \int_0^T \int_{-S(t)}^{S(t)} \abs{W'*(\bar\rho^N-\rho)} \dx \dt
	+ \norm{m'}_{L^\infty([0,R(T)])}
		\int_0^T \int_{-S(t)}^{S(t)} \abs{\bar\rho^N-\rho}\cdot\abs{U} \dx \dt \\
&\leq \norm{v}_\infty R(T) \int_0^T
	\norm{W'}_{L^1([-2S(t),2S(t)])} \norm{\bar\rho^N-\rho}_{L^1(\setR)} \dt
		\spliteq
	+ \bigl[ \norm{v}_\infty + R(T)\norm{v'}_{L^\infty([0,R(T)])} \bigr]
		\int_0^T \norm{V-W'*\rho}_{L^\infty([-S(t),S(t)])}
		\int_\setR \abs{\bar\rho^N-\rho} \dx \dt \\
&\leq 4 \norm{v}_\infty R(T) F(T) G\bigl(2S(T)\bigr) S(T)
	\norm{\bar\rho^N-\rho}_{L^1([0,T]\times\setR)}
		\spliteq
	+ \bigl[ \norm{v}_\infty + R(T)G\bigl(R(T)\bigr) \bigr]
		F(T)\bigl[G\bigl(S(T)\bigr)+G\bigl(2S(T)\bigr)\bigr]
		\norm{\bar\rho^N-\rho}_{L^1([0,T]\times\setR)},
\end{split}
\]
which goes to $0$ because $\norm{\bar\rho^N-\rho}_{L^1([0,T]\times\setR)}\to0$.

Let us now move on to proving that $\rho$ solves the entropy inequality \eqref{eq:entropy-condition}.
We split the left hand side of \eqref{eq:discrete-entropy} into three terms:
\begin{gather*}
\int_0^T\int_\setR \abs{\bar\rho^N -c} \partial_t \phi \d x \d t, \\
\int_0^T\int_\setR \sign(\bar\rho^N -c)
	\bigl(m(\bar\rho^N) - m(c)\bigr)\bar U^N\partial_x \phi \d x \d t, \\
\int_0^T\int_\setR \sign(\bar\rho^N -c) m(c) \phi \partial_x \bar U^N \d x \d t.
\end{gather*}
The first two can be shown to pass to the limit in a similar fashion as in \cite[Lemma 4.4]{DiFrancescoFagioliRadici}, in particular exploiting the fact that $r\mapsto\sign(r-c)\bigl(m(r)-m(c)\bigr)$ is a Lipschitz function.

Let us now deal with the third integral. Notice that the sought convergence is implied by
\begin{gather}
\lim_{N\to\infty} \int_0^T\int_\setR \sign(\bar\rho^N -c) m(c) \phi
	\partial_x \bigl(\bar U^N - U\bigr)\d x \d t = 0,
	\label{eq:limit1} \\
\lim_{N\to\infty} \int_0^T\int_\setR \bigl(\sign(\bar\rho^N -c)-\sign(\rho -c)\bigr)
	m(c) \phi \partial_x U\d x \d t = 0 .
	\label{eq:limit2}
\end{gather}

To prove the limit \eqref{eq:limit1}, assume w.l.o.g.\ that $\supp(\phi)\subset(0,T)\times (-a,a)$. Then the limit follows from the estimate
\[
\begin{split}
\abs*{\int_0^T\int_\setR \sign(\bar\rho^N -c) m(c) \phi
	\partial_x \bigl(\bar U^N - U\bigr)\d x \d t}
&\leq \int_0^T\int_\setR \abs*{\partial_x (\bar U^N - U)} \phi \d x \d t \\
&\leq \norm{\phi}_\infty \int_0^T \norm{\partial_x (\bar U^N - U)(t,\plchldr)}_{L^1((-a,a))} \dt
\end{split}
\]
and the fact that for every $t\in(0,T)$ by Young inequality we have
\[
\begin{split}
\norm{\partial_x (\bar U^N - U)(t,\plchldr)}_{L^1((-a,a))}
&= \norm{\partial_x [W'*(\bar\rho^N-\rho)]}_{L^1((-a,a))} \\
&\leq \norm{W''*(\bar\rho^N-\rho)}_{L^1((-a,a))}
	+ \abs{w(t)}\cdot\norm{\bar\rho^N-\rho}_{L^1(\setR)} \\
&\leq \bigl(\norm{W''}_{L^1([-a-S(t),a+S(t)])} + \abs{w(t)}\bigr)
	\norm{\bar\rho^N-\rho}_{L^1(\setR)} \\
&\leq F(T) \bigl[1+G\bigl(a+S(T)\bigr)\bigr] \norm{\bar\rho^N-\rho}_{L^1(\setR)},
\end{split}
\]
which goes to zero because $\norm{\bar\rho^N-\rho}_{L^1(\setR)}\to0$.

To prove \eqref{eq:limit2}, we first need to study the behavior of $\abs{\partial_x U}$ on the set where $\rho=c$.
Let $Z=\set{(t,x)\in (0,\infty)\times\setR}{\rho(t,x)=c}$ and consider an open set $\Omega\supset Z$, to be optimized later. Consider the family of squares in $(0,T)\times\setR$
\[
\mathcal{F} = \set{Q(z,r)}{z\in Z,\ Q(z,r)\subset\Omega}
\]
covering $Z$, where $Q(z,r)=\set{w=(w_t,w_x)}{\abs{w_t-z_t}<r,\abs{w_x-z_x}<r}\subset(0,T)\times\setR$. By Besicovitch covering theorem there exists $M\in\setN$ and a countable subfamily $\mathcal{G}\subset\mathcal{F}$ such that
\begin{equation}\label{eq:besicovitch}
\bm{1}_Z \leq \sum_{Q\in\mathcal{G}} \bm{1}_Q \leq M \bm{1}_\Omega.
\end{equation}
Let us now focus the attention on a single square $Q=I\times A\subset(0,T)\times\setR$. For positive $\chi\in C^1_c(I)$ and $\eta\in C^1_c(A)$, consider the test function $\phi(t,x)=\chi(t)\eta(x)$. Since $\rho$ solves \eqref{eq:continuity-equation}, we have
\[
\int_I\int_A \bigl[
	\rho\partial_t\chi\eta + m(\rho)U\chi\partial_x\eta \bigr] \dx\dt = 0.
\]

We know that $\rho(t,\plchldr)\in BV(\setR)$ for every $t$, but from the equation \eqref{eq:conservation-law} we deduce that also $\partial_t\rho$ is a measure, therefore $\rho\in BV\bigl((0,T)\times\setR\bigr)$ and, by Leibniz and the chain rule, also $m(\rho)U$ belongs to $BV\bigl((0,T)\times\setR\bigr)$. Hence, integrating by parts and looking only at the absolutely continuous parts of the distributional derivatives, we have
\[
\int_I\int_A \bigl[
	(\partial_t\rho)^a + \partial_x\bigl(m(\rho)U\bigr)^a \bigr] \chi\eta \dx\dt
= \int_I\int_A \bigl[
	(\partial_t\rho)^a + m'(\rho)(\partial_x\rho)^aU
		+ m(\rho)\partial_xU \bigr] \chi\eta \dx\dt
= 0.
\]
By \cite[Proposition 3.92]{AFP} we have that $(\partial_t\rho)^a=(\partial_x\rho)^a=0$ a.e.\ on $Z$, hence
\[
m(c) \iint_{Q\cap Z} \partial_xU \chi\eta \dx\dt
= - \iint_{Q\setminus Z} \bigl[ m(\rho)\partial_xU
	+(\partial_t\rho)^a + m'(\rho)(\partial_x\rho)^aU\bigr] \chi\eta \dx\dt.
\]
By a standard argument (considering linear combinations of terms $\chi\eta$) we deduce
\[
m(c) \iint_{Q\cap Z} \abs{\partial_xU} \dx\dt
\leq \iint_{Q\setminus Z} \abs{ m(\rho)\partial_xU
	+(\partial_t\rho)^a + m'(\rho)(\partial_x\rho)^aU } \dx\dt.
\]
Summing this inequality over $Q\in\mathcal{G}$ and recalling \eqref{eq:besicovitch} therefore we get
\[
\begin{split}
m(c) \iint_Z \abs{\partial_xU} \dx\dt
&= m(c) \iint_{[0,\infty)\times\setR} \abs{\partial_xU}\bm1_Z \bm1_Z \dx\dt \\
&\leq m(c) \sum_{Q\in\mathcal{G}} \iint_{[0,\infty)\times\setR}
	\abs{\partial_xU}\bm1_Z \bm1_Q \dx\dt \\
&= m(c) \sum_{Q\in\mathcal{G}} \iint_{Q\cap Z} \abs{\partial_xU} \dx\dt \\
&\leq \sum_{Q\in\mathcal{G}} \iint_{Q\setminus Z} \abs{ m(\rho)\partial_xU
	+(\partial_t\rho)^a + m'(\rho)(\partial_x\rho)^aU } \dx\dt \\
&= \sum_{Q\in\mathcal{G}} \iint_{[0,\infty)\times\setR} \abs{ m(\rho)\partial_xU
	+(\partial_t\rho)^a + m'(\rho)(\partial_x\rho)^aU } \bm1_{Z^c}\bm1_Q \dx\dt \\
&\leq M \iint_{[0,\infty)\times\setR} \abs{ m(\rho)\partial_xU
	+(\partial_t\rho)^a + m'(\rho)(\partial_x\rho)^aU } \bm1_{Z^c}\bm1_\Omega \dx\dt \\
&= M \iint_{\Omega\setminus Z} \abs{ m(\rho)\partial_xU
	+(\partial_t\rho)^a + m'(\rho)(\partial_x\rho)^aU } \dx\dt.
\end{split}
\]
By the absolute continuity of the integral, letting $\Omega$ decrease towards $Z$ we deduce
\[
m(c) \iint_Z \abs{\partial_xU} \dx\dt = 0.
\]
We can split the integral \eqref{eq:limit2} in $Z$ and $Z^c$. The integral in $Z$ is identically zero thanks to the previous computation. On the other hand, in $Z^c$ we have that $\sign(\bar\rho^N -c)-\sign(\rho -c)\to0$ a.e.\ because $\bar\rho^N\to\rho$ a.e.\ and $r\mapsto\sign(r-c)$ is continuous for $r\neq c$. Therefore by Lebesgue dominated convergence (a domination is $2m(c)\phi\abs{\partial_xU}\in L^1$) we have that
\[
\lim_{N\to\infty} \iint_{Z^c} \bigl(\sign(\bar\rho^N -c)-\sign(\rho -c)\bigr)
	m(c) \phi \partial_x U\d x \d t = 0 .
\]
This concludes the proof of \eqref{eq:limit2}.
\end{proof}

\begin{proof}[Proof of \autoref{thm:KR}]
We follow very closely the proof of \cite[Theorem 1.3]{KarlsenRisebro}, with the crucial difference that our vector fields $P,Q$ depend on time and there is the error term in the right hand side of the entropy inequality.

We fix a smooth, compactly supported, non-negative test function $\phi\in C^\infty_c\bigl(\setR\times\setR^n;[0,\infty)\bigr)$ and define
\[
\tilde\phi(t,x,\tau,y)
= \phi\oleft(\frac{t+\tau}{2},\frac{x+y}{2}\right)
	\delta_h\oleft(\frac{t-\tau}{2}\right)
	\omega_h\oleft(\frac{x-y}{2}\right) ,
\]
where $\delta\in C^\infty_c\bigl(\setR;[0,\infty)\bigr)$ and $\omega\in C^\infty_c\bigl(\setR^n;[0,\infty)\bigr)$ are two kernels such that $\int_\setR \delta(t)\d t=1$ and $\int_{\setR^n}\omega(x)\d x=1$ and
\begin{align*}
\delta_h(t) &= h^{-1}\delta(t/h), & \omega_h(x) &= h^{-n}\omega(x/h^n).
\end{align*}


We write the entropy inequality for $\rho(t,x)$ plugging the constant $\sigma(\tau,y)$ and integrate with respect to $\tau$ and $y$; we then do the converse treatment to the entropy inequality for $\sigma(\tau,y)$. The resulting inequalities are as follows:
\begin{align*}
\iiiint \bigl\{ \abs{\rho-\sigma}\partial_t\tilde\phi + \sign(\rho-\sigma) \bigl[
	\bigl(p(\rho)-p(\sigma)\bigr)P(t,x) \cdot \nabla_x\tilde\phi
	&- p(\sigma) \div_xP(t,x) \tilde\phi \bigr] \bigr\} \d t \d x \d \tau \d y \\
&\geq -\frac1N H(T)(\norm{\partial_x\tilde\phi}_\infty+\norm{\tilde\phi}_\infty), \\
\iiiint \bigl\{ \abs{\rho-\sigma}\partial_\tau\tilde\phi + \sign(\rho-\sigma) \bigl[
	\bigl(q(\rho)-q(\sigma)\bigr)Q(\tau,y) \cdot \nabla_y\tilde\phi
	&+ q(\rho) \div_yQ(\tau,y) \tilde\phi \bigr] \bigr\} \d t \d x \d \tau \d y \\
&\geq -\frac1N H(T)(\norm{\partial_y\tilde\phi}_\infty+\norm{\tilde\phi}_\infty),
\end{align*}
where we used the implicit notation $\rho=\rho(t,x)$, $\sigma=\sigma(\tau,y)$, $\tilde\phi=\tilde\phi(t,x,\tau,y)$ to shorten the formulas.

The proof now proceeds as in \cite{KarlsenRisebro} by summing the two inequalities, and rearranging the terms in the same manner. Notice in particular that the only manipulations performed are integrations by parts exclusively in space, so the time dependence does not interfere in these computations.

With the choice of a test function $\phi_k(t,x)=\psi_k(x)\chi_k(t)$ with $\psi_k\in C^\infty_c\bigl(\setR;[0,1]\bigr)$, $\psi_k\to1$, $\abs{\psi_k'}\leq 1$, and $\chi_k$ that approximates $\bm1_{[t_1,t_2]}$ for some fixed $0\leq t_1<t_2\leq T$, noticing that crucially the right hand side does not depend on $\partial_t\phi$, we can pass to the limit for $h,k\to\infty$ as in \cite{KarlsenRisebro} and arrive at the conclusion \eqref{eq:KR},
observing that all the integrals in space are localized to $[-S(T),S(T)]$ because the densities vanish outside.
\end{proof}

Finally, we can combine the results of this section to give a short proof of our main theorem.

\begin{proof}[Proof of \autoref{thm existence of entropy solution}]
We first deal with the existence.
Let $\bar\rho^N : [0,T] \times \setR \to \setR$ be the piecewise constant density associated to particles $X^N$ solving either \eqref{eq:ode-integrated-interaction} or \eqref{eq:ode-sampled-interaction}. Thanks to \autoref{thm:compactness} and \autoref{cor:continuity}, up to a subsequence, we have that $\bar\rho^N$ converges to a density $\rho\in L^\infty_\loc\bigl([0,\infty)\times\setR\bigr) \cap C\bigl([0,\infty);L^1(\setR)\bigr)$. Moreover, since $\bar\rho^N$ enjoys the uniform estimates stated in \autoref{lem:supp-bound}, \autoref{lem:max-princ} and \autoref{prop:total-variation}, passing to the limit and using the lower semicontinuity of the total variation one deduces \eqref{eq:solution-estimates}.

Combining \autoref{prop:discrete-entropy-condition} and \autoref{prop:entropic-limit}, we further obtain that $\rho$ satisfies the entropy inequality \eqref{eq:entropy-condition}, thus it is an entropy solution of \eqref{eq:conservation-law} according to \autoref{def:entropy-solution}.

Regarding the uniqueness, assume that $\rho$ and $\sigma$ are two entropy solutions. Thanks to \autoref{thm:KR}, with a similar argument as in \autoref{rmk:stability} which does not even involve the errors terms in the entropy inequality, one deduces with an application of Gronwall lemma the following stability estimate
\[
\norm{\rho(t)-\sigma(t)}_{L^1(\setR)}
\leq \norm{\rho(0)-\sigma(0)}_{L^1(\setR)}
	\exp\bigl(t A(T)\bigr), \qquad \forall t\in(0,T)
\]
for some increasing function $A:[0,\infty)\to[0,\infty)$. This leads to the uniqueness under the assumptions that the two initial conditions are the same.

The convergence rate \eqref{eq:convergence-rate} is deduced with a similar computation as in \autoref{rmk:stability} applied to the pair of densities $\rho$ and $\bar\rho^N$ (or, alternatively, sending $M\to\infty$ in \eqref{eq:cauchy}).
\end{proof}

\section{Numerical simulations}\label{sec:numerics}

\subsection{Implementation}

In this section we describe the numerical implementation of the two particle schemes introduced above and show some simulations of the evolution of the particles and the associated piecewise constant densities. The actual implementation is in the Julia programming language, using the \texttt{DifferentialEquation.jl} \cite{rackauckas2017differentialequations} package to solve the system of ODE describing the evolution of the particles. The code presented in this section is included and further developed in the package \texttt{ConservationLawsParticles.jl} \cite{ConservationLawsParticles}.

Assuming that one has at disposal the functions
\begin{minted}{julia}
V(t::Real, x::Real) :: Real
Wprime(t::Real, r::Real) :: Real
mobility(rho::Real) :: Real
\end{minted}
defining the external velocity field, the auto-interaction and the mobility term respectively, the resulting velocity of the particles according to the sampled scheme can be computed in-place with
\begin{minted}{julia}
function velocity(dx, x, p, t)
    R = pwc_density(x)
    for i in 1:length(x)
        v = (-sum(Wprime(t, x[i] - x[j]) for j in 1:i-1)
             -sum(Wprime(t, x[i] - x[j]) for j in i+1:length(x)))
        v /= length(x) - 1
        v += V(t, x[i])
        mob = mobility(v < 0 ? R[i] : R[i+1])
        dx[i] = v * mob
    end
end

function pwc_density(x::AbstractVector{<:Real})
    len = length(x)
    R = Array{float(eltype(x))}(undef, len+1)
    R[1] = R[end] = 0
    d = 1 / (len-1)
    for i in 2:len
        R[i] = d / (x[i] - x[i-1])
    end
    R
end
\end{minted}

A solution to the ODEs can be found by initializing an array \verb|x0| to the initial position of the particles whose piecewise density approximates the desired initial distribution of mass in the sense of \eqref{eq:initial-conditions}. The system of ODEs can then be solved with
\begin{minted}{julia}
using DifferentialEquations
tspan = (0., 3.)
prob = ODEProblem(velocity, x0, tspan)
sol = solve(prob, BS5(); abstol=1e-7, reltol=1e-7);
\end{minted}

The solution can then be visually inspected by either showing the trajectories of the particles or by displaying the density at different times.

\begin{figure}[h!]
\includegraphics[width=0.5\textwidth]{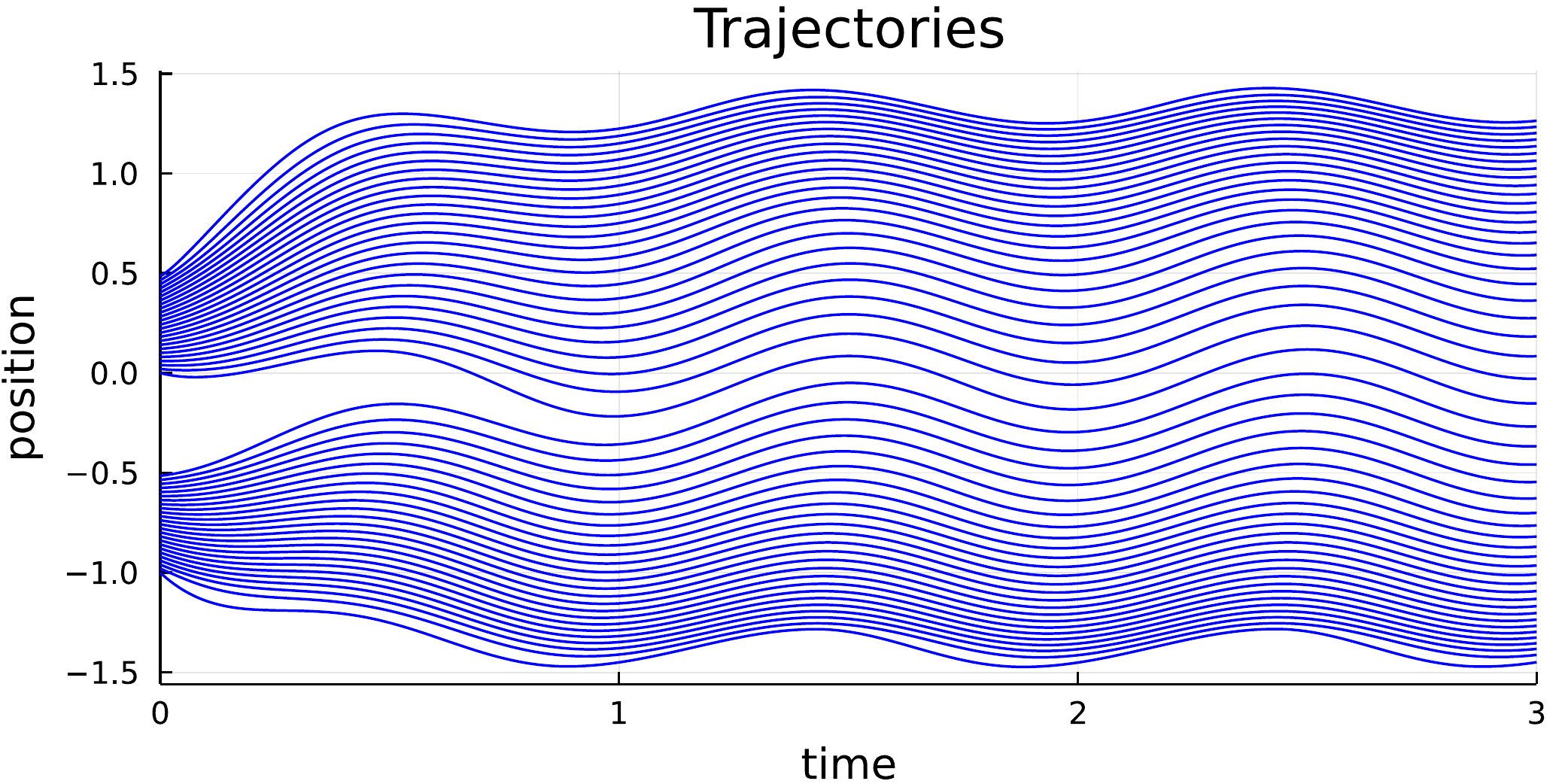}%
\includegraphics[width=0.5\textwidth]{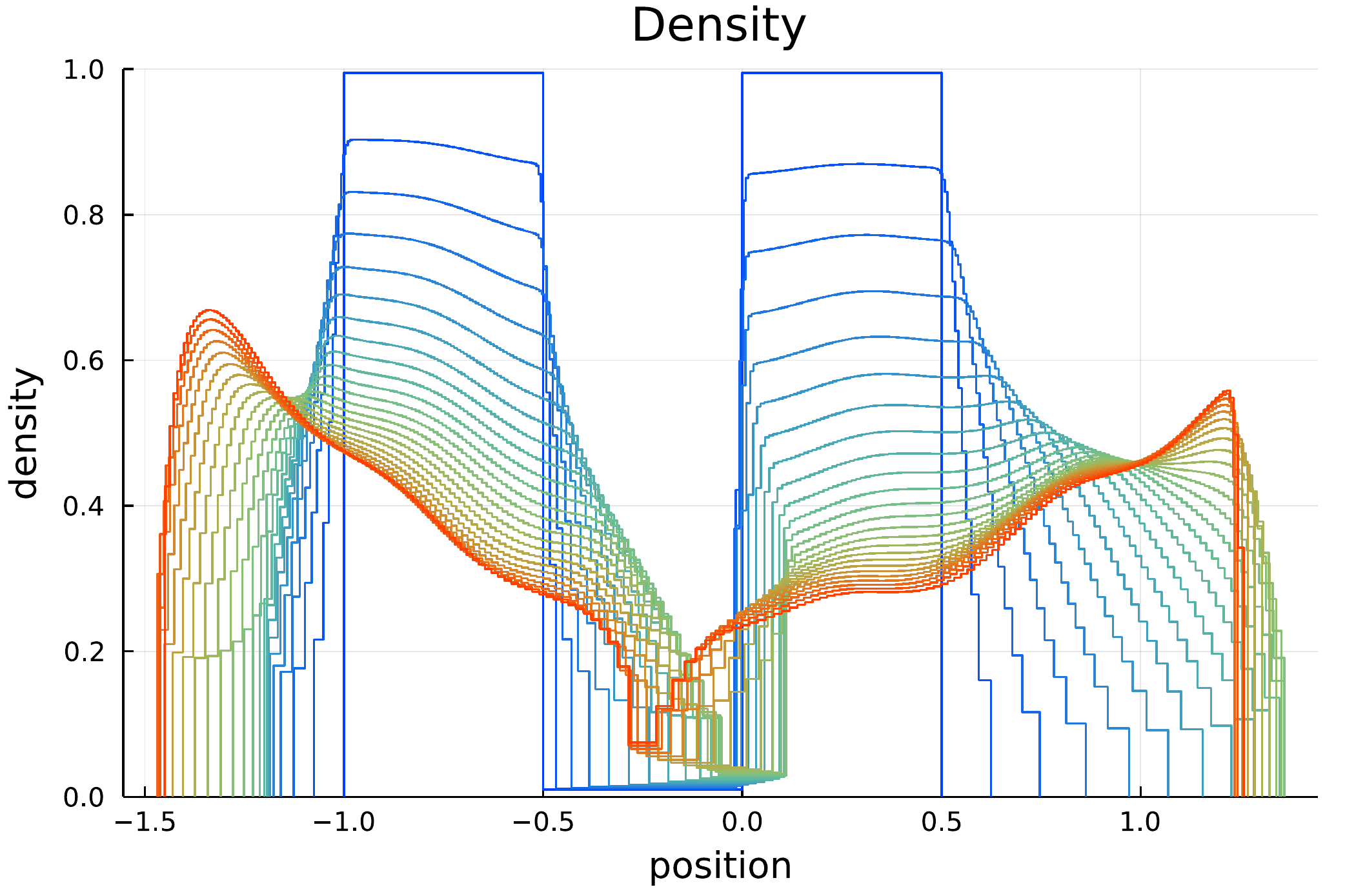}%
\caption{plot of the trajectories (left) and densities (right) for the example.}\label{fig:traj-dens}
\end{figure}

As an example, consider the time-dependent problem defined by
\begin{minted}{julia}
V(t, x) = -x^3 + 0.02sin(12x) + sinpi(2t)
Wprime(t, x) = -5sign(x) / (abs(x) + 1)
mob(rho) = 1 / (rho + 1)
model = SampledModel((V,), ((Wprime,),), (mob,))
\end{minted}
We can visualize the trajectories of the particles and the densities at different times, as shown in \autoref{fig:traj-dens}.

The code computing the velocities of the particles can be generalized to more than one species. To make the code more manageable, we introduce two types \verb|Model| and \verb|IntegratedModel| which describe the whole system of $N$ species (with external velocities \verb|Vs|, mutual interactions \verb|Wprimes|/\verb|Ws| and mobilities \verb|mobilities|) and auxiliary functions which computes the interaction effect of a species on each individual particle. The two types correspond to the sampled and integrated scheme.
\begin{minted}{julia}
mutable struct Model{
    N,
    TVs         <: Tuple{Vararg{Any,N}},
    TWprimes    <: Tuple{Vararg{Tuple{Vararg{Any,N}},N}},
    Tmobilities <: Tuple{Vararg{Any,N}},
}
    Vs::TVs
    Wprimes::TWprimes
    mobilities::Tmobilities
end

mutable struct IntegratedModel{
    N,
    TVs         <: Tuple{Vararg{Any,N}},
    TWs         <: Tuple{Vararg{Tuple{Vararg{Any,N}},N}},
    Tmobilities <: Tuple{Vararg{Any,N}},
}
    Vs::TVs
    Ws::TWs
    mobilities::Tmobilities
end

function sampled_interaction(t::Real, Wprime, ys::AbstractVector{<:Real}, x::Real)
    sum(Wprime(t, x - y) for y in ys) / (length(ys) - 1)
end

function integrated_interaction(t::Real, W, dens_diff::AbstractVector{<:Real},
                                ys::AbstractVector{<:Real}, x::Real)
    sum(i -> dens_diff[i] * W(t, x - ys[i]), eachindex(ys))
end
\end{minted}

The total velocity describing the evolution of the full system is then computed by two methods of the function \verb|velocities| which call either \verb|sampled_interaction| or \verb|integrated_interaction|.
\begin{minted}{julia}
function velocities(
  dx::ArrayPartition{F, T},
  x::ArrayPartition{F, T},
  p::Model{N, TVs, TWprimes, Tmobilities},
  t
) where {
  F,
  T <: Tuple{Vararg{AbstractVector{<:Real}}},
  N, TVs, TWprimes, Tmobilities
}
  dens = pwc_densities(x.x...)
  for spec in 1:N
    for i in 1:length(x.x[spec])
      v::F = p.Vs[spec](t, x.x[spec][i])
      for other in 1:N
        v -= sampled_interaction(t, p.Wprimes[spec][other],
                                 x.x[other], x.x[spec][i])
      end
      if v < 0
        mob = p.mobilities[spec](dens[spec][:, 1, i]...)
      else
        mob = p.mobilities[spec](dens[spec][:, 2, i]...)
      end
      dx.x[spec][i] = v * mob
    end
  end
end

function velocities(
  dx::ArrayPartition{F, T},
  x::ArrayPartition{F, T},
  p::IntegratedModel{N, TVs, TWs, Tmobilities},
  t
) where {
  F,
  T <: Tuple{Vararg{AbstractVector{<:Real}}},
  N, TVs, TWs, Tmobilities
}
  dens = pwc_densities(x.x...)
  dens_diff = Vector{Vector{F}}(undef, N)
  for s in 1:N
    dens_diff[s] = dens[s][s, 1, :] - dens[s][s, 2, :]
  end
  for spec in 1:N
    for i in 1:length(x.x[spec])
      v::F = p.Vs[spec](t, x.x[spec][i])
      for other in 1:N
        v -= integrated_interaction(t, p.Ws[spec][other], dens_diff[other],
                                    x.x[other], x.x[spec][i])
      end
      if v < 0
        mob = p.mobilities[spec](dens[spec][:, 1, i]...)
      else
        mob = p.mobilities[spec](dens[spec][:, 2, i]...)
      end
      dx.x[spec][i] = v * mob
    end
  end
end
\end{minted}

The multi-species version of the densities computation is a bit more involved. It scans the particles left to right and computes the densities of the various species on both sides of each particle.
\begin{minted}{julia}
function pwc_densities(xs::Vararg{AbstractVector{<:Real}, N}) where N
    N::Int
    T = promote_type(float.(eltype.(xs))...)
    len = length.(xs)
    dens = map(x -> new_undef_densities(T, N, length(x)), xs)::NTuple{N, Array{T, 3}}
    ind = ones(Int, N) # indices of current particles being examined
    ds = zeros(T, N) # current densities
    fs::NTuple{N, T} = 1 ./ (len .- 1) # normalization factors
    i_min = Int[0]
    @inbounds while true
        n_mins::Int = 0
        for i in 1:N
            if ind[i] <= len[i]
                if n_mins == 0 || xs[i][ind[i]] < xs[i_min[1]][ind[i_min[1]]]
                    i_min[1] = i
                    n_mins = 1
                elseif xs[i][ind[i]] <= xs[i_min[1]][ind[i_min[1]]]
                    n_mins += 1
                    if n_mins <= length(i_min)
                        i_min[n_mins] = i
                    else
                        push!(i_min, i)
                    end
                end
            end
        end
        if n_mins == 0
            break
        end
        i_mins = @view i_min[1:n_mins]
        for i in i_mins, s in 1:N
            dens[i][s, 1, ind[i]] = ds[s]
        end
        for i in i_mins
            if ind[i] < length(xs[i])
                ds[i] = fs[i] / (xs[i][ind[i]+1] - xs[i][ind[i]])
            else
                ds[i] = zero(T)
            end
        end
        for i in i_mins, s in 1:N
            dens[i][s, 2, ind[i]] = ds[s]
        end
        for i in i_mins
            ind[i] += 1
        end
    end
    dens
end
\end{minted}

This multi-species implementation can be demonstrated in an example such as
\begin{minted}{julia}
V1(t, x) = - (x - 1)^3 / 15 + 0.05sin(12x)
V2(t, x) = - (x + 1 - 4sin(5t))^3 / 5
W_attr(t, x) = 5log(abs(x) + 1)
W_rep(t, x) = -5log(abs(x) + 1)
mob1(rho, sigma) = max(1 - rho - 0.5sigma, 0)
mob2(rho, sigma) = max(2 - sigma - 0.5rho, 0)

model = IntegratedModel(
    (V1, V2),
    ((W_rep, W_attr), (W_attr, W_rep)),
    (mob1, mob2))

x0 = ArrayPartition(
    gaussian_particles(3, 75) .- 2,
    gaussian_particles(3, 75) .+ 2)

tspan = (0., 5.)
prob = ODEProblem(velocities, x0, tspan, model)

@time sol = solve(prob, BS5(); abstol=1e-6, reltol=1e-6);
\end{minted}

\begin{figure}[h!]
\includegraphics[width=0.5\textwidth]{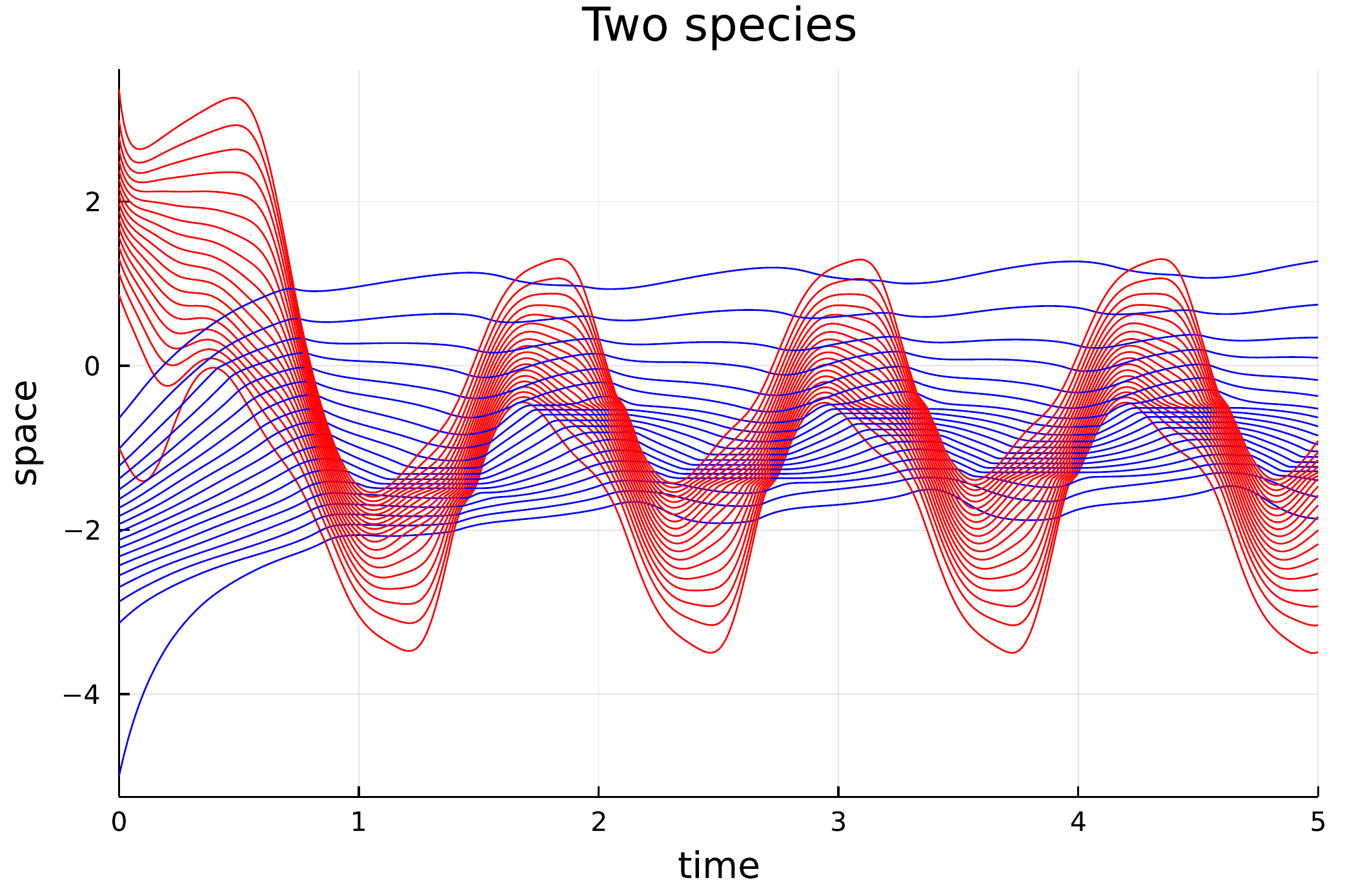}%
\includegraphics[width=0.5\textwidth]{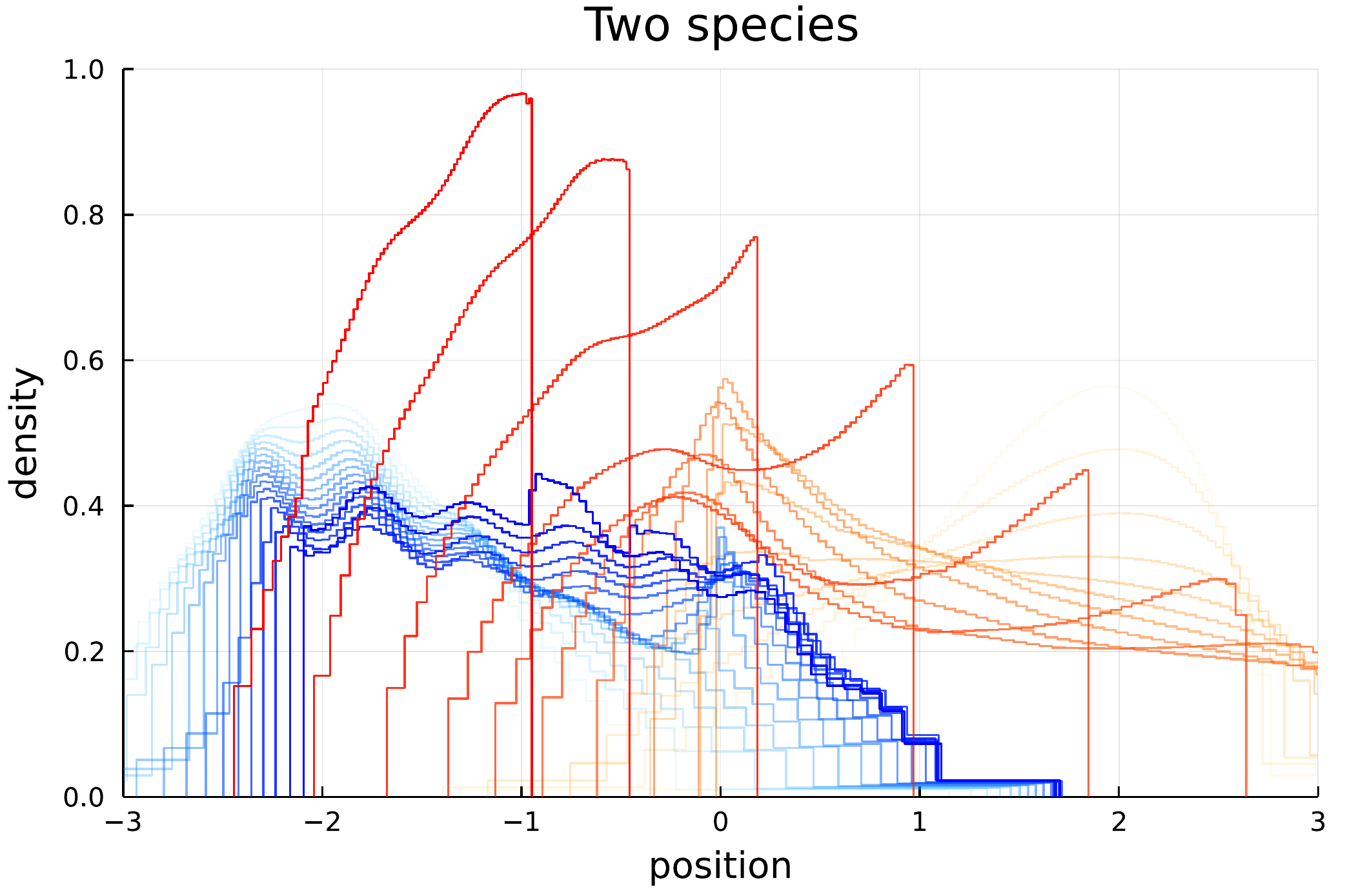}
\caption{trajectories and densities for the multi-species example.}\label{fig:multi-traj-dens}
\end{figure}

We want to point out that a possible downside of our piecewise constant construction of the density is that it has intrinsically lower resolution where the density is low. This is because the particles are necessarily more spread apart, hence they cannot capture finer details of the density. This might be problematic because these lower resolution regions can influence adversely the resolution of the other species, even where the latter have higher densities and hence are more detailed. The reason is that with a decaying mutual interaction the dominant terms are contributed by the local particles of the lower quality density, which produce a velocity field of lower quality. This effect can be observed as staircase-like artifacts which can develop in the densities. Over long periods of time, it is observed that these can lead to oscillations, which are further amplified if the sampled scheme is used. With our approach the only simple way to limit these artifacts seems to be to increase the total number of particles.

\begin{figure}[h!]
\centering
\includegraphics[width=0.5\textwidth]{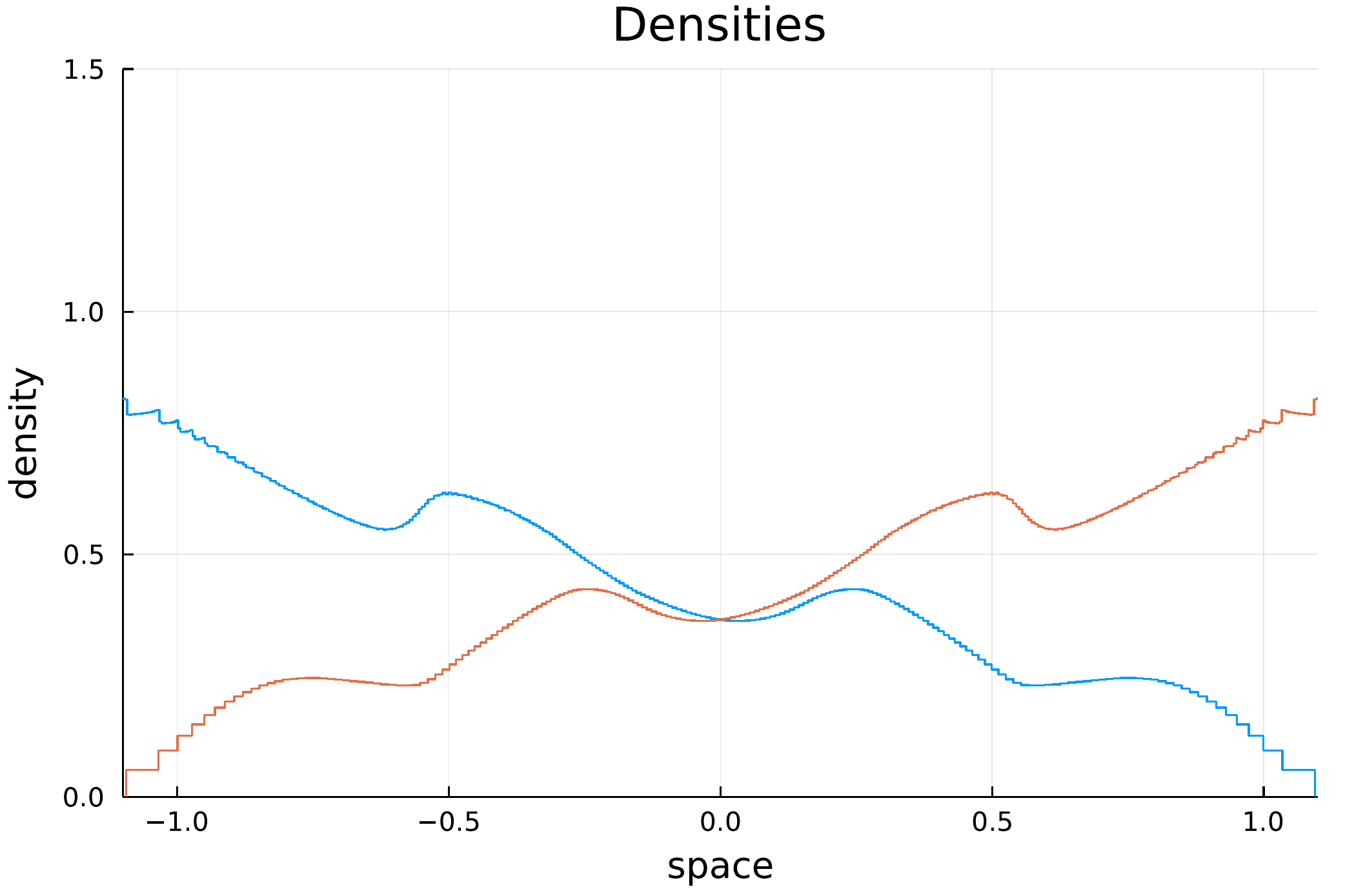}
\caption{example of staircase artifacts. Observe that near $x=-1$, even if the blue density is high and has a good resolution, it exhibits some staircase artifacts which are caused by the local interaction with the lower resolution orange density.}
\end{figure}

\begin{figure}[h!]
\centering
\includegraphics[width=0.5\textwidth]{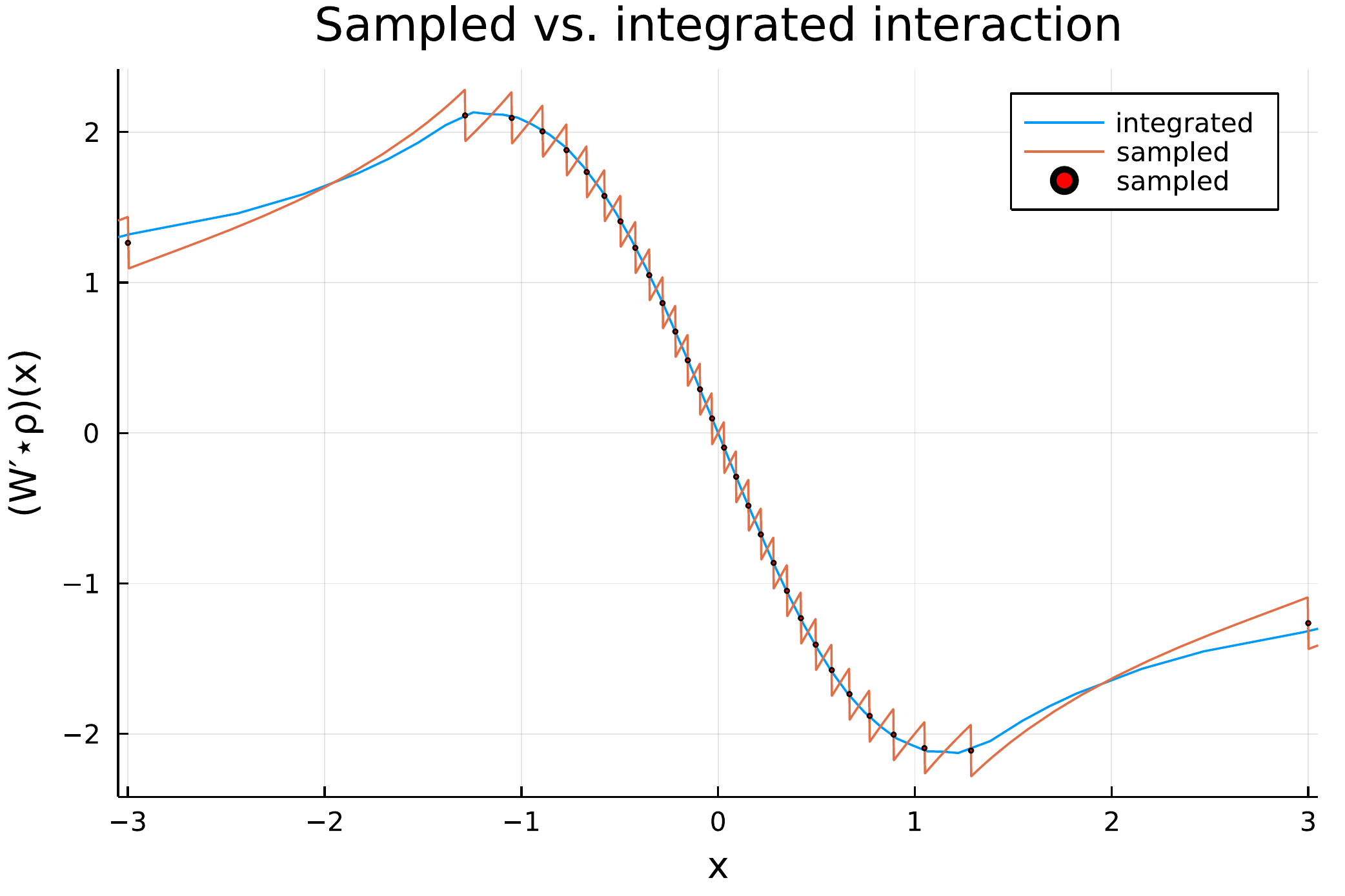}
\caption{sampled interaction velocity field (in blue) and integrated interaction velocity field. The particles follow a Gaussian distribution, the interaction is attractive. As it can be seen, the sampled interaction is discontinuous in correspondence of every particle.}
\label{fig:samp-int}
\end{figure}

When dealing with a multi-species problem, the usage of the sampled interaction is inadequate, since the resulting velocity field is discontinuous in correspondence of every particle, as can be visualized in the example \autoref{fig:samp-int}. This is problematic because the particles of two different species can become arbitrarily close, hence the effects of the discontinuity are very relevant and can be observed: for instance, particles of different species can become locked together and might not flow past each other. In contrast, the integrated scheme produces a velocity field which is continuous and suitable to be evaluated at any desired position, hence it does not cause this locking or jump behavior between particles.

\begin{figure}
\includegraphics[width=0.5\textwidth]{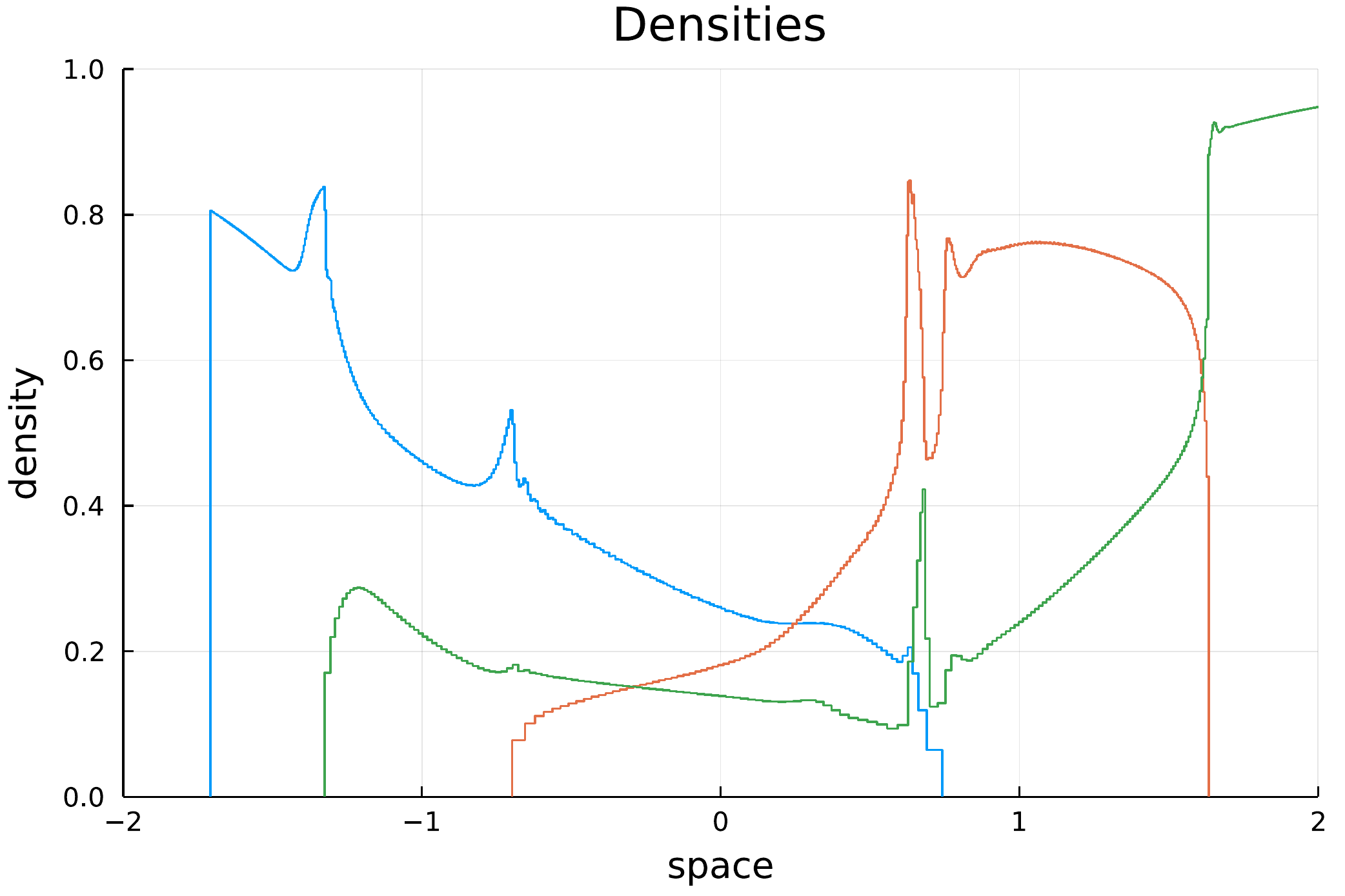}%
\includegraphics[width=0.5\textwidth]{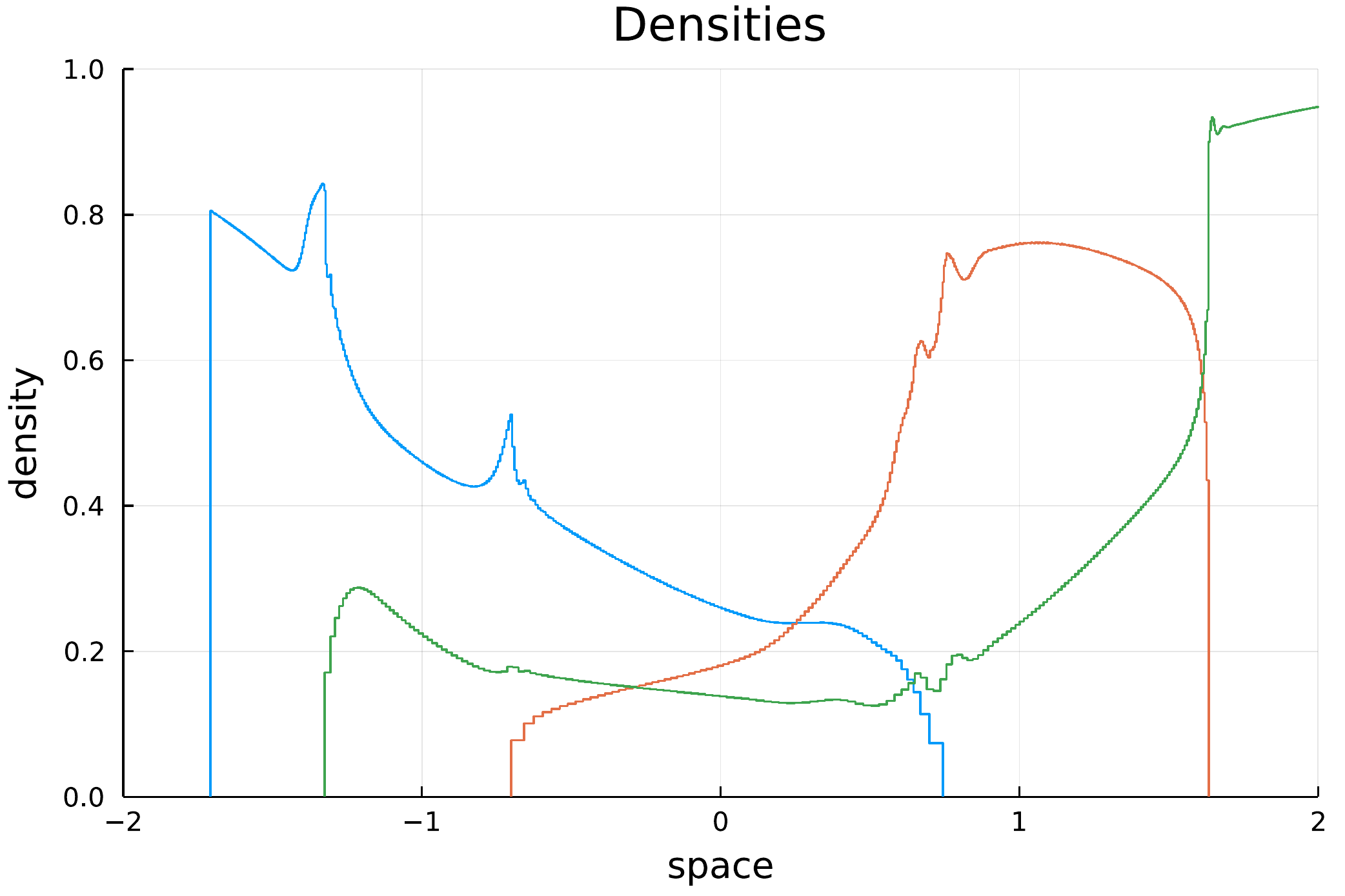}
\caption{example of artifacts in a model with three species. Notice that near $x=1$ there is the development of oscillations caused by the sampled interaction between the orange and the green species. These oscillation amplify and the solution becomes completely untrustworthy after this time. A way to remedy is to use the integrated interaction scheme, which is better behaved. The third picture shows the solution at the same instant in time solved with the same number of particles according to the integrated scheme.}
\end{figure}

Since in general the integrated scheme is not computationally more expensive than the sampled scheme, but the quality of the solution can be much higher, it is to be generally preferred, especially for multi-species models. Even in the single species case, the better regularity of the integrated scheme sometimes can allow the numerical integrator to perform longer steps in time, leading to higher performance.

\subsection{Examples}

In this section we show some peculiar examples that demonstrate the feasibility of the numerical computation with both schemes.

\begin{example}
We begin with an example with a single species, a compactly supported mobility which enforces $\rho\leq1$ and time-dependent potentials. The initial density is discretized with $100$ equally spaced particles in $[-1,-1/2]$ and $100$ particles in $[0,1/2]$. In the plot only the trajectory of one every 8 particles is displayed.
\begin{align*}
V(t, x) &= -(x - \sin(3t))^3, & v(\rho) &= (1 - \rho)_+, \\
W(t, x) &= - 5\sin(4t)^2 \log(\abs{x} + 1), & \rho_0 &= \bm1_{[-1,-1/2]} + \bm1_{[0,1/2]} .
\end{align*}
\begin{center}
\includegraphics[width=0.75\textwidth]{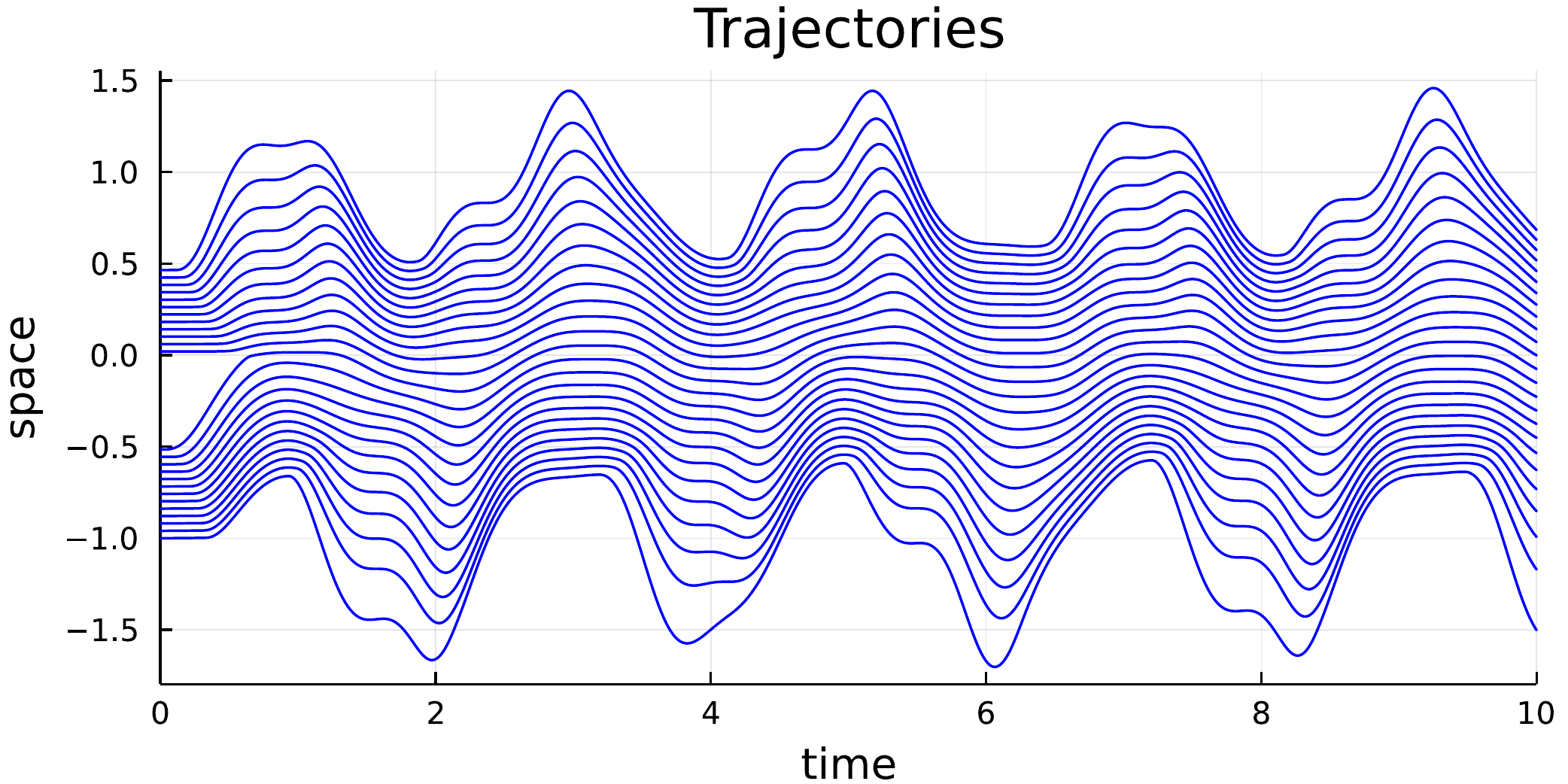}
\end{center}
\end{example}

\begin{example}
In this example we consider a strictly positive mobility $v(\rho) = 1 / (1 + \rho)$ which satisfies \eqref{eq:assumption-v-decay}. Consistently with \autoref{lem:max-princ}, the density does not blow up in finite time, but can grow arbitrarily large as time goes on. We plot the trajectories, showing the particles collapsing to the origin, and the density at different times on a logarithmic scale. The potentials are
\begin{align*}
V(t,x) &= 0, & W(t, x) = 5 \log(\abs{x} + 1),
\end{align*}
while the initial density is a perturbation of $\frac12\bm1_{[-1,1]}$ given by particles at position
\[
x_i = \frac{2i-N}{N} + 0.01\sin\oleft(20\frac{2i-N}{N}\right).
\]
\includegraphics[width=0.5\textwidth]{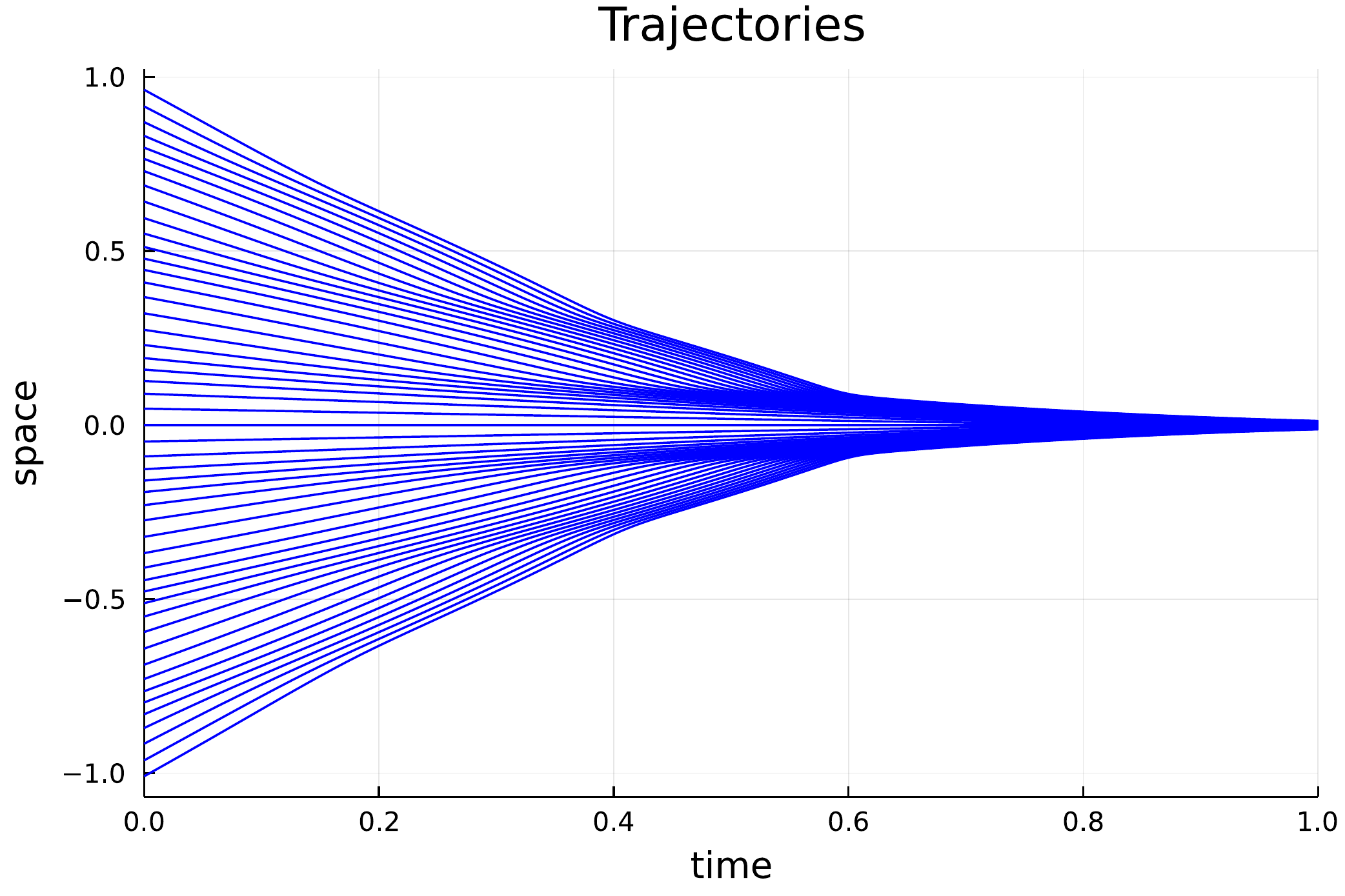}%
\includegraphics[width=0.5\textwidth]{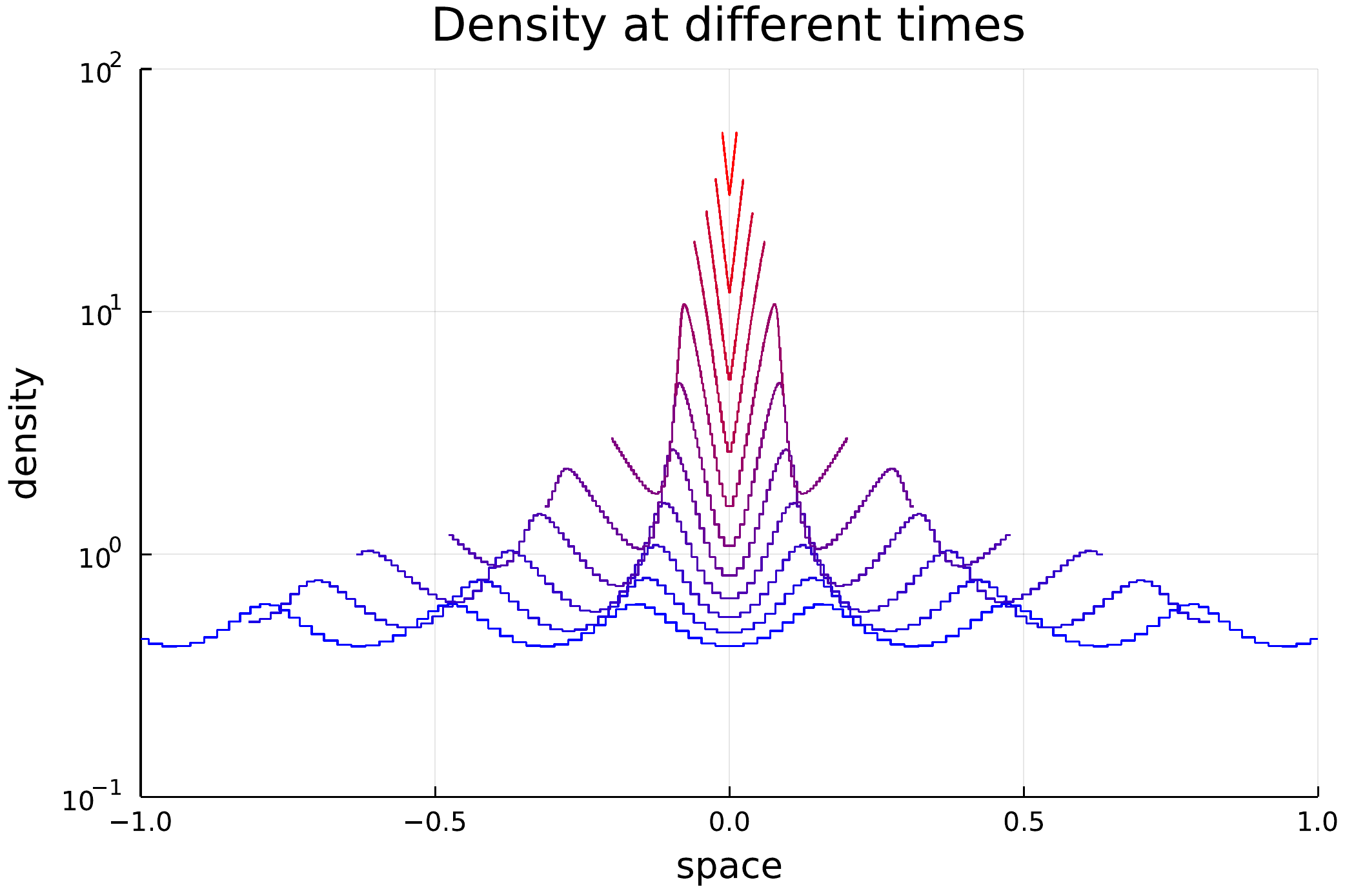}
\end{example}

\begin{example}
We compare the sampled and integrated schemes with time-dependent potentials, non-vaninishing mobility and $30$ particles.
\begin{align*}
V(t, x) &= -(x - sin(3t))^3, & v(\rho) &= 1 / (\rho + 1)^2, \\
W(t, x) &= \sin(2\pi x) / (2\pi x), & \rho_0 &= \bm1_{[-1,-1/2]}+\bm1_{[0,1/2]}.
\end{align*}
The two sets of trajectories seem to coincide to a much higher degree than suggested by \eqref{eq:comparison-0}, because that follows from a pessimistic estimate applied to each interval $[x_{i-1},x_i]$, whereas in reality compensations between consecutive intervals can occur. For example, in very regular situations with monotone potentials one could expect $\dot U-\bar U$ to be of order $1/N^2$ instead of $1/N$.
In the following picture we kept the number of particles extremely low otherwise the trajectories would have been indistinguishable.
\begin{center}
\includegraphics[width=0.75\textwidth]{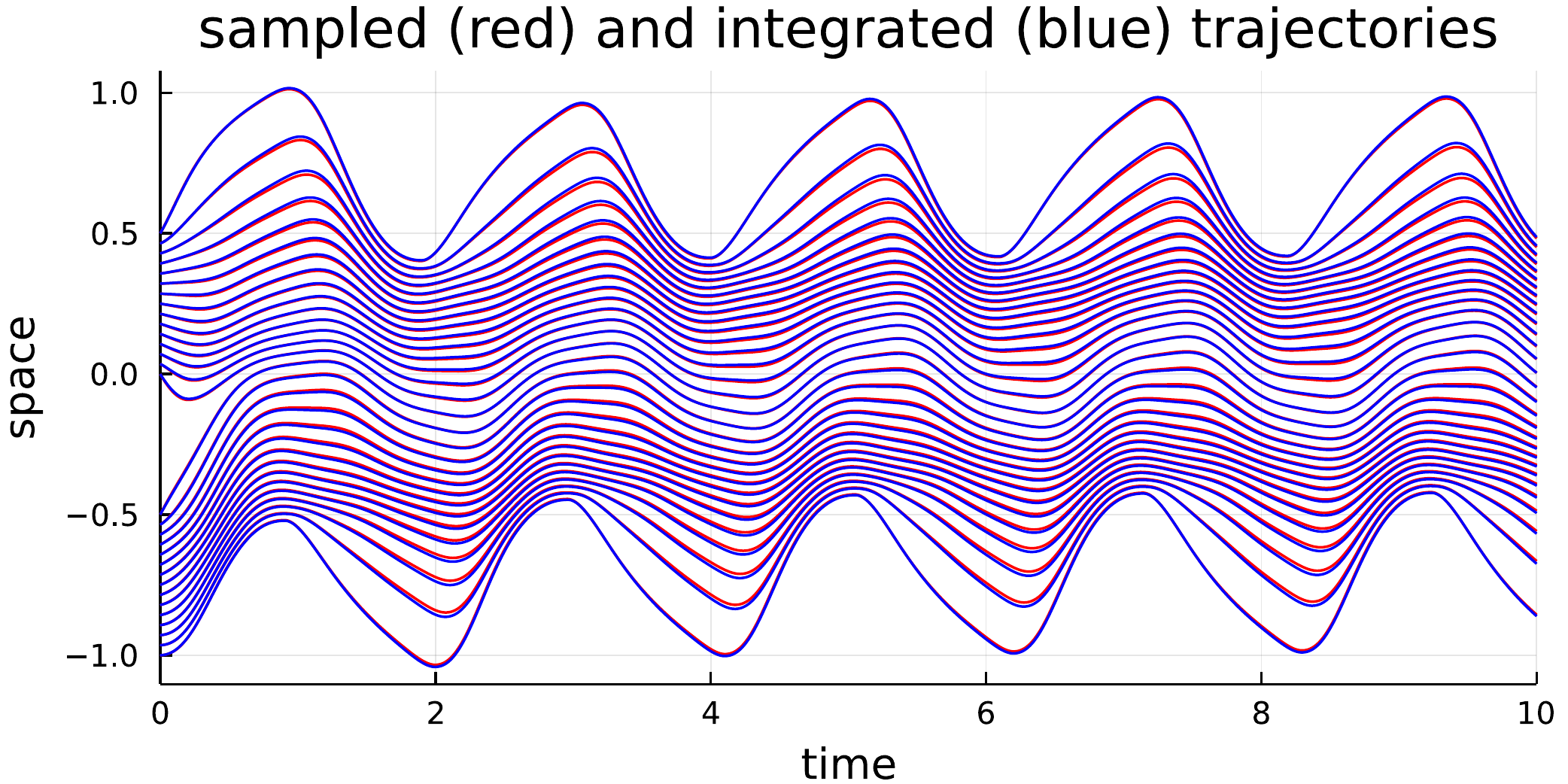}
\end{center}
\end{example}

\begin{example}
In the same setting as the previous example, we now compare the quality of the solution with different number of particles. From \eqref{eq:convergence-rate} we expect a convergence rate of order $1/N$. In the following picture we plot the trajectories of 31 particles in red, and one every 8 particles out of a total of 241 particles in blue.
\begin{center}
\includegraphics[width=0.75\textwidth]{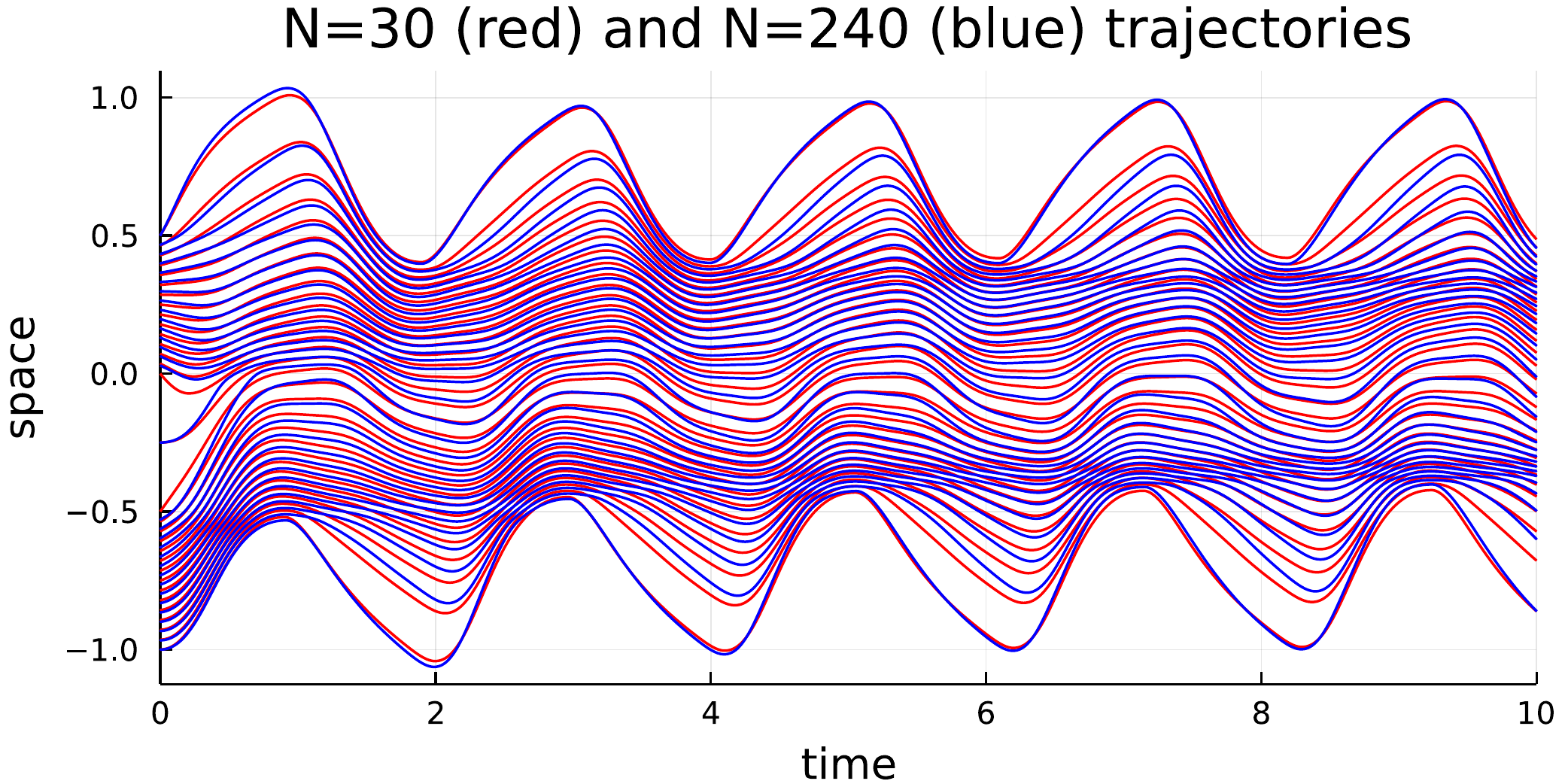}
\end{center}
\end{example}

\begin{example}
We propose an example in which the entropy solution is explicit (in particular it is stationary) and allows us to compare our new scheme with the one described in \cite{DiFrancescoFagioliRadici}. In order to fall inside the applicability domain of their scheme, we need to consider $V=0$, a regular attractive $W$ and a compactly supported $v$. We consider
\begin{align*}
V(t,x) &= 0, & v(\rho) &= (1-\rho)_+, \\
W(t,x) &= \log(\abs{x}+1)\frac{\abs{x}}{\abs{x}+1}, & \rho_0 &= \bm1_{[-2,-3/2]} + \bm_{[3/2,2]}.
\end{align*}
In the limit as $N\to\infty$, both schemes recover the correct stationary solution $\rho(t)=\rho_0$. However, while our new scheme correctly captures it with any number $N$, the old scheme (pictured in red in the left column) merges in finite time the two peaks of the density into a single bump, with a time of collapse which becomes increasingly large as $N\to\infty$.
\begin{center}
\includegraphics[width=0.5\textwidth]{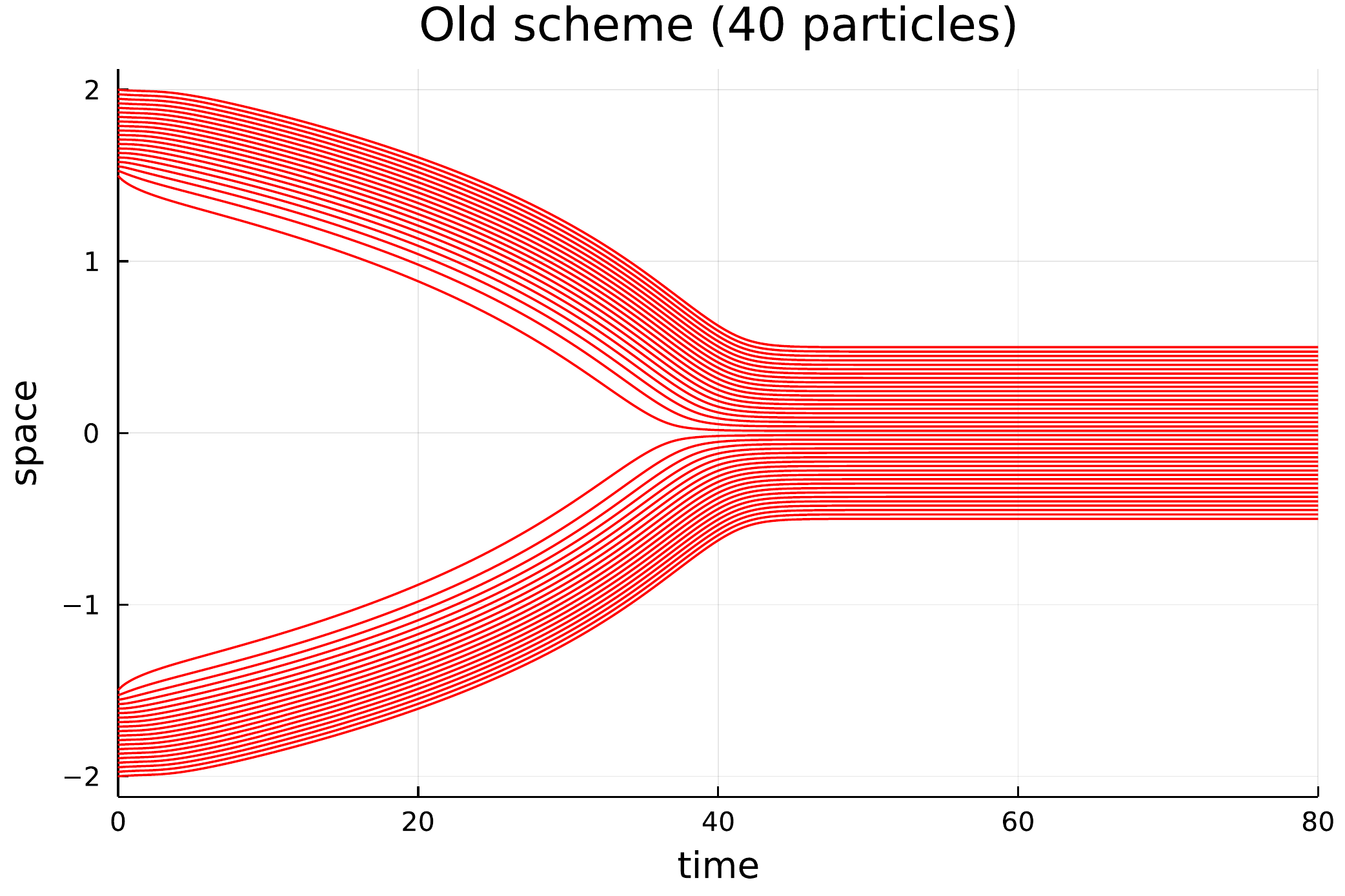}%
\includegraphics[width=0.5\textwidth]{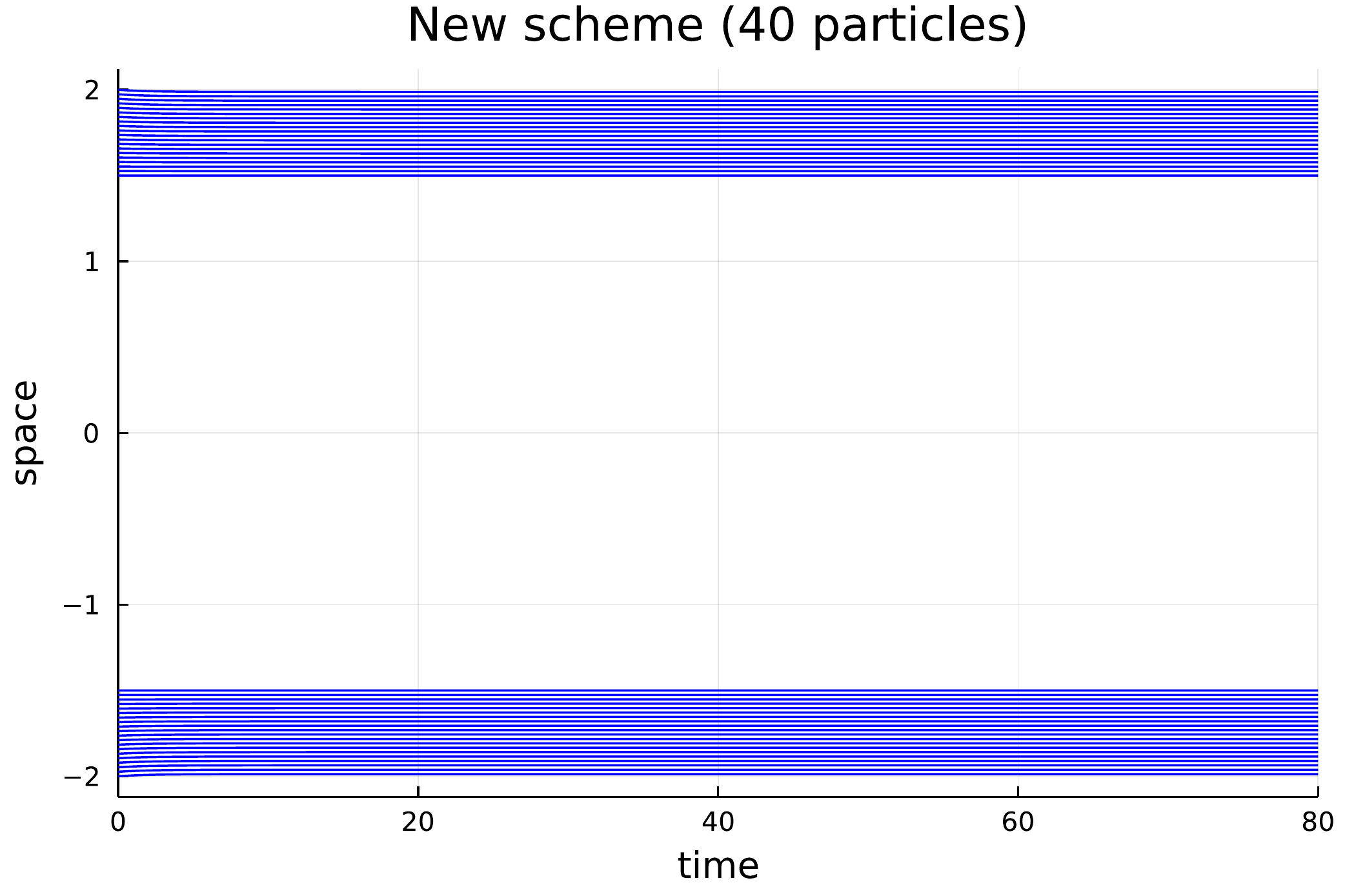}\\
\includegraphics[width=0.5\textwidth]{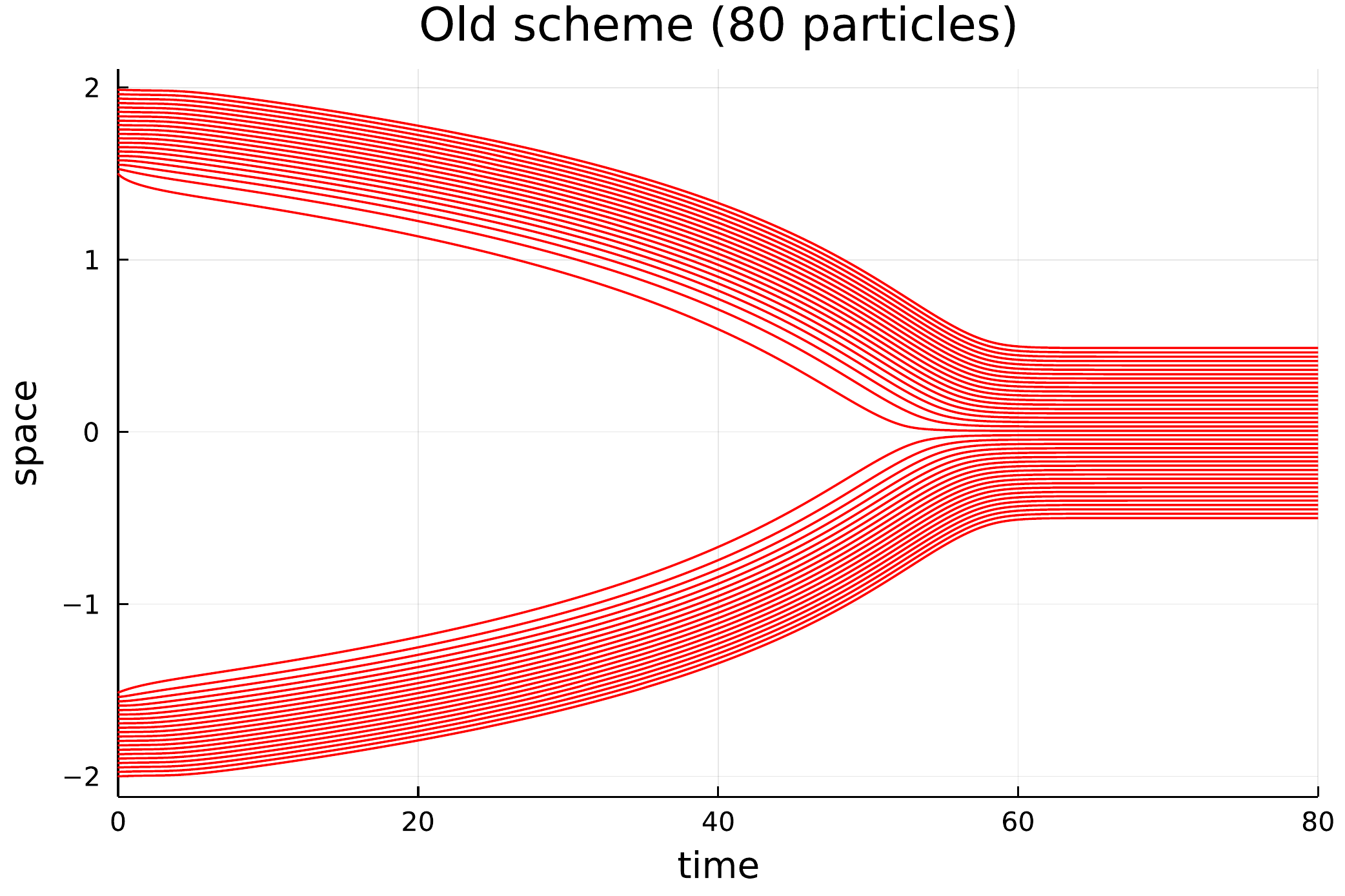}%
\includegraphics[width=0.5\textwidth]{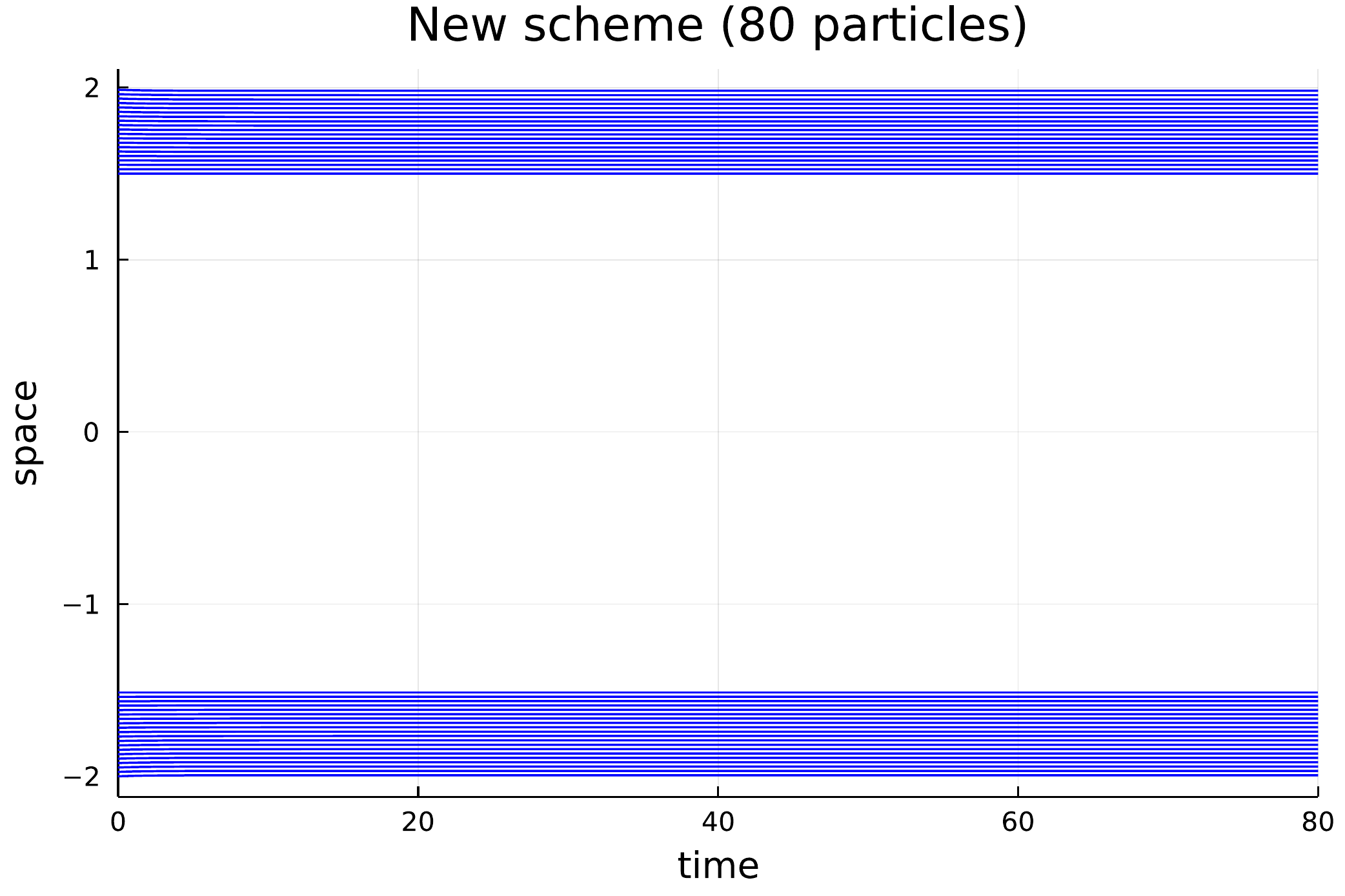}
\end{center}
\end{example}

\begin{example}
In this last example we showcase an application of the integrated scheme to a multi-species model of 2 populations that cross each other in opposite directions. With the notation of \eqref{eq:multi-species}, we use $S=2$ and
\begin{align*}
V_1(t,x) &= 2, &
W_{1,1}(t,x) &= W_{2,2}(t,x) = 2 (e^{\abs{x}/4} + e^{-2\abs{x}}), \\
V_2(t,x) &= -2, &
W_{1,2}(t,x) &= W_{2,1}(t,x) = -2 \log(\abs{x} + 1), \\
v_1(\rho_1,\rho_2) &= (2-\rho_1-\rho_2/2)_+, &
v_2(\rho_1,\rho_2) &= (2-\rho_2-\rho_1/2)_+, \\
\rho_1(0) &= \bm1_{[-2,-3/2]}+\bm1_{[-1,-1/2]}, &
\rho_2(0) &= \bm1_{[1/2,3/2]}.
\end{align*}
This model describes a density $\rho_1$ that wants to travel to the right, a density $\rho_2$ that wants to travel to the left; the mutual interaction is repulsive, whereas the self interaction is repulsive at small scales and attractive at long range.
\begin{center}
\includegraphics[width=0.75\textwidth]{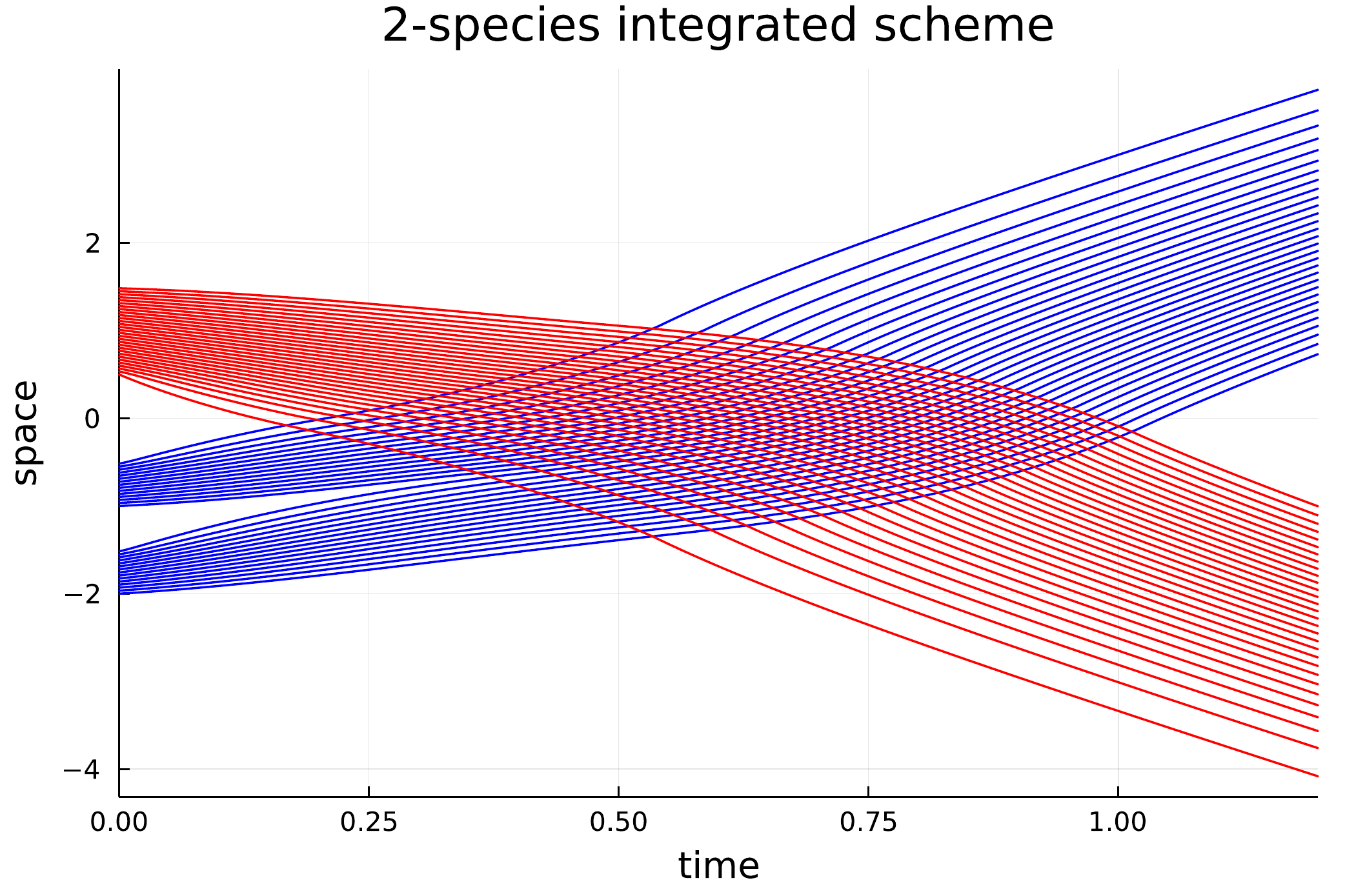}
\end{center}
\end{example}

\phantomsection
\addcontentsline{toc}{section}{\refname}
\printbibliography

\end{document}